\documentclass[11pt,english]{article}
\usepackage[T1]{fontenc}
\usepackage[latin9]{inputenc}
\usepackage{geometry}
\usepackage{color}
\usepackage{array}
\usepackage{rotating}
\usepackage{units}
\usepackage{multirow}
\usepackage{amsmath}
\usepackage{amssymb}
\usepackage{graphicx}
\usepackage{setspace}
\makeatletter

\providecommand{\tabularnewline}{\\}
\newcommand{\lyxdot}{.}

\makeatother

\usepackage{babel}
\begin{document}

\title{Sparse Linear Models and Two-Stage Estimation in High-Dimensional
Settings with Possibly Many Endogenous Regressors%
\thanks{First and foremost, I thank James Powell and Martin Wainwright for
useful suggestions and helpful comments. I am also grateful to Ron
Berman, Elena Manresa, Minjung Park, Demian Pouzo, Miguel Villas-Boas,
and other participants at the UC Berkeley econometric seminar. All
errors are my own. This work was supported by Haas School of Business
at Berkeley. %
}}

\author{Ying Zhu%
\thanks{Haas School of Business, UC Berkeley. 2220 Piedmont Ave., Berkeley,
CA 94720. ying\_zhu@haas.berkeley.edu. Tel: 406-465-0498. Fax: 510-643-4255 %
}}
\maketitle
\begin{abstract}
This paper explores the validity of the two-stage estimation procedure
for sparse linear models in high-dimensional settings with possibly
many endogenous regressors. In particular, the number of endogenous
regressors in the main equation and the instruments in the first-stage
equations can grow with and exceed the sample size $n$. The analysis
concerns the \textit{exact sparsity} case, i.e., the maximum number
of \textcolor{black}{non-zero} components in the vectors of parameters
in the first-stage equations, $k_{1}$, and the number of non-zero
components in the vector of parameters in the second-stage equation,
$k_{2}$, are allowed to grow with $n$ but slowly compared to $n$.
I consider the high-dimensional version of the two-stage least square
estimator where one obtains the fitted regressors from the first-stage
regression by a least square estimator with $l_{1}$- regularization
(the Lasso or Dantzig selector) when the first-stage regression concerns
a large number of instruments relative to $n$, and then construct
a similar estimator using these fitted regressors in the second-stage
regression.\textcolor{black}{{} }The main theoretical results of this
paper are non-asymptotic bounds from which I establish sufficient scaling conditions on the sample size 
for estimation consistency in $l_{2}-$norm and variable-selection
consistency (i.e., the two-stage high-dimensional estimators correctly
select the non-zero coefficients in the main equation with high probability).\textcolor{black}{{}
A technical issue regarding the so-called {}``restricted eigenvalue
(RE) condition'' for estimation consistency and the {}``mutual incoherence
(MI) condition'' for selection consistency arises in the two-stage
estimation procedure from allowing the number of regressors
in the main equation to exceed $n$ and this paper provides analysis
to verify these RE and MI conditions.} Depending on the underlying
assumptions that are imposed, the upper bounds on the $l_{2}-$error
and the sample size required to obtain these consistency results differ
by factors involving $k_{1}$ and/or $k_{2}$. Simulations are conducted
to gain insight on the finite sample performance of the high-dimensional
two-stage estimator.

\newpage{}
\end{abstract}

\section{Introduction }

The objective of this paper is consistent estimation and selection
of regression coefficients in models with a large number of endogenous
regressors. Consider the linear model 
\begin{equation}
y_{i}=\mathbf{x}_{i}^{T}\beta^{*}+\epsilon_{i}=\sum_{j=1}^{p}x_{ij}\beta_{j}^{*}+\epsilon_{i},\quad i=1,...,n
\end{equation}
where $\epsilon_{i}$ is a zero-mean random error possibly correlated
with $\mathbf{x}_{i}$ and $\beta^{*}$ is an unknown vector of parameters
of our main interests. The $j^{th}$ component of $\beta^{*}$ is
denoted by $\beta_{j}^{*}$. A component in the $p-$dimensional vector
$\mathbf{x}_{i}$ is said to be \textit{endogenous} if it is correlated
with $\epsilon_{i}$ (i.e., $\mathbb{E}(\mathbf{x}_{i}\epsilon_{i})\neq\mathbf{0}$)
and \textit{exogenous} otherwise (i.e., $\mathbb{E}(\mathbf{x}_{i}\epsilon_{i})=\mathbf{0}$).
Without loss of generality, I will assume all regressors are endogenous
throughout the rest of this paper for notational convenience (a modification
to allow mix of endogenous and exogenous regressors is trivial.).
When endogenous regressors are present, the classical least squares
estimator will be inconsistent for $\beta^{*}$ (i.e., $\hat{\beta}_{OLS}\overset{p}{\nrightarrow}\beta^{*}$)
even when the dimension $p$ of $\beta^{*}$ is small relative to
the sample size $n$. The classical solution to this problem of endogenous
regressors supposes that there is some $L$-dimensional vector of
instrumental variables, denoted by $\mathbf{z}_{i}$, which is observable
and satisfies $\mathbb{E}(\mathbf{z}_{i}\epsilon_{i})=\mathbf{0}$
for all $i$. In particular, the two-step estimation procedures including
the two-stage least square (2SLS) estimation and the control function
approach play an important role in accounting for endogeneity that
comes from individual choice or market equilibrium (e.g., Wooldrige,
2002). Consider the following {}``first stage'' equations for the
components of $\mathbf{x}_{i}$ 
\begin{equation}
x_{ij}=\mathbf{z}_{ij}^{T}\pi_{j}^{*}+\eta_{ij}=\sum_{l=1}^{d_{j}}z_{ijl}\pi_{jl}^{*}+\eta_{ij},\qquad i=1,....,n,\; j=1,...,p.
\end{equation}
For each $j=1,...,p$, $\mathbf{z}_{ij}$ is a $d_{j}\times1$ vector
of instrumental variables, and $\eta_{ij}$ a zero-mean random error
which is uncorrelated with $\mathbf{z}_{ij}$, and $\pi_{j}^{*}$
is an unknown vector of nuisance parameters. I will refer to the equation
in (1) as the main equation (or second-stage equation) and the equations
in (2) as the first-stage equations. Throughout the rest of this paper,
I will impose the following assumption. Without loss of generality,
this assumption implies a triangular simultaneous equations model
structure.\\
\\
\textbf{Assumption 1.1}: The data \textbf{$\{y_{i},\,\mathbf{x}_{i},\,\mathbf{z}_{i}\}_{i=1}^{n}$
}are \textit{i.i.d.} with finite second moments; $\mathbb{E}(\mathbf{z}_{ij}\epsilon_{i})=\mathbb{E}(\mathbf{z}_{ij}\eta_{ij})=\mathbf{0}$
for all $j=1,...,p$ and $\mathbb{E}(\mathbf{z}_{ij}\eta_{ij^{'}})=\mathbf{0}$
for all $j\neq j^{'}$. \\
\\
Statistical estimation and variable selection in the high-dimensional
setting is concerned with models in which the dimension of the parameters
of interests is larger than the sample size. In the past decade, a
tremendous increase of research activities in this field has been
facilitated by the advances in data collection technology. In the
literature on high-dimensional sparse linear regression models, a
great deal of attention has been given to the $l_{1}-$penalized least
squares. In particular, the Lasso and the Dantzig selector are the
most studied techniques (see, e.g., Tibshirani, 1996; Cand�s and Tao,
2007; Bickel, Ritov, and Tsybakov, 2009; Belloni, Chernozhukov, and
Wang, 2011; Belloni and Chernozhukov, 2011b; Loh and Wainwright, 2012;
Negahban, Ravikumar, Wainwright, and Yu, 2012). Variable selection
when the dimension of the problem is larger than the sample size has
also been studied in the likelihood method setting with penalty functions
other than the $l_{1}-$norm (see, e.g., Fan and Li, 2001; Fan and
Lv, 2011). Lecture notes by Koltchinskii (2011), as well as recent
books by B�hlmann and van de Geer (2011) and Wainwright (2014) have
given a more comprehensive introduction to high-dimensional statistics. 

The Lasso procedure is a combination of the residual sum of squares
and a $l_{1}-$regularization defined by the following program 
\[
\min_{\beta\in\mathbb{R}^{p}}\left\{ \frac{1}{2n}|y-X\beta|_{2}^{2}+\lambda_{n}|\beta|_{1}\right\} .
\]
Denote the minimizer to the above program by $\hat{\beta}_{Las}$.
A necessary and sufficient condition of $\hat{\beta}_{Las}$ is that
$0$ belongs to the subdifferential of the convex function $\beta\mapsto\frac{1}{2n}|y-X\beta|_{2}^{2}+\lambda_{n}|\beta|_{1}$.
This implies that the Lasso solution $\hat{\beta}_{Las}$ satisfies
the constraint 
\[
\left|\frac{1}{2n}X^{T}(y-X\hat{\beta}_{Las})\right|_{\infty}\leq\lambda_{n}.
\]
The Dantzig selector of the linear regression function is defined
as a vector having the smallest $l_{1}-$norm among all $\beta$ satisfying
the above constraint, i.e., 
\[
\hat{\beta}_{Dan}=\arg\min\left\{ |\beta|_{1}:\left|\frac{1}{2n}X^{T}(y-X\beta)\right|_{\infty}\leq\lambda_{n}\right\} .
\]
Recently, these $l_{1}-$penalized techniques have been applied in
a number of economics studies. Caner (2009) studies a Lasso-type GMM
estimator. Rosenbaum and Tsybakov (2010) study the high-dimensional
errors-in-variables problem where the non-random regressors are observed
with additive error and they present an application to hedge fund
portfolio replication. Lecture notes by Belloni and Chernozhukov (2011b)
discuss the $l_{1}-$based penalization methods with various econometric
problems including earning regressions and instrumental selection
in Angrist and Krueger data (1991). Belloni and Chernozhukov (2011a)
study the $l_{1}$-penalized quantile regression and illustrate its
use on an international economic growth application. Belloni, Chen,
Chernozhukov, and Hansen (2012) estimate the optimal instruments using
the Lasso and in an empirical example dealing with the effect of judicial
eminent domain decisions on economic outcomes, they find the Lasso-based
instrumental variable estimator outperforms an intuitive benchmark.
Belloni, Chernozhukov, and Hansen (2012) propose robust methods for
inference on the effect of a treatment variable on a scalar outcome
in the presence of very many controls with an application to abortion
and crime. Fan, Lv, and Li (2011) review the literature on sparse
high-dimensional econometric models including the vector autoregressive
model for measuring the effects of monetary policy, panel data model
for forecasting home price, and volatility matrix estimation in finance.
Their discussion is not restricted to $l_{1}-$based regularization
methods.

High dimensionality arises in the triangular simultaneous equations
structure (1) and (2) when the dimension $p$ of $\beta^{*}$ is large
relative to the sample size $n$ (namely, $p\gg n$) or when the dimension
$d_{j}$ of $\pi_{j}^{*}$ is large relative to the sample size $n$
(namely, $d_{j}\gg n$) for at least one $j$. In this paper, I consider
the scenario where the number of non-zero coefficients in $\beta^{*}$
and $\pi_{j}^{*}$ is small relative to $n$ (i.e., $\beta^{*}$ and
$\pi_{j}^{*}$for $j=1,...,p$ are \textit{exactly} \textit{sparse}).
The case where $d_{j}\gg n$ for at least one $j$ but $p\ll n$ has
been considered by Belloni and Chernozhukov (2011b), where they showed
the instruments selected by the Lasso technique in the first-stage
regression can produce an efficient estimator with a small bias at
the same time. To the best of my knowledge, the case where $p\gg n$
and $d_{j}\ll n$ for all $j$, or the case where $p\gg n$ and $d_{j}\gg n$
for at least one $j$ in the context of triangular simultaneous equations
with two-stage estimation has not been studied in the literature.
In both cases, one can still use the ideas of the 2SLS estimation
together with the Lasso technique. For instance, in the case where
$p\gg n$ and $d_{j}\ll n$ for all $j$, one can obtain the fitted
regressors by a standard least square estimation on each of the first-stage
equations separately as usual and then apply a Lasso-type technique
with these fitted regressors in the second-stage regression. Similarly,
in the case where $p\gg n$ and $d_{j}\gg n$ for all $j$, one can
obtain the fitted regressors by performing a regression with a Lasso-type
estimator on each of the first-stage equations separately and then
apply another Lasso-type estimator with these fitted regressors in
the second-stage regression. 

Compared to existing two-stage techniques which limit the number of
regressors entering the first-stage equations or the second-stage
equation or both, the two-stage estimation procedures with $l_{1}-$regularization
in both stages are more flexible and particularly powerful for applications
in which the vector of parameters of interests is sparse and there
is lack of information about the relevant explanatory variables and
instruments. In terms of practical implementations, these above-mentioned
high-dimensional two-stage estimation procedures are intuitive and
can be easily implemented using existing software packages for the
standard Lasso-type technique for linear models without endogeneity.
In analyzing the statistical properties of these estimators, the extension
from models with a few endogenous regressors to models with many endogenous
regressors ($p\gg n$) in the context of triangular simultaneous equations
with two-stage estimation is not obvious. This paper aims to explore
the validity of these two-step estimation procedures for the triangular
simultaneous linear equation models in the high-dimensional setting
under the sparsity scenario. 

In the presence of endogenous regressors, the direct implementation
of the Lasso or Dantzig selector fails as sparsity of coefficients
in equation (1) does not correspond to sparsity of linear projection
coefficients. The linear instrumental variable model with a single
or a few endogenous regressors and many instruments has been studied
in the econometrics literature on high dimensional models. For example,
Belloni and Chernozhukov (2011b) consider the following triangular
simultaneous equation model: 
\begin{eqnarray*}
y_{i} & = & \theta_{0}+\theta_{1}x_{1i}+\mathbf{x}_{2i}^{T}\gamma+\epsilon_{i}\\
x_{1i} & = & \mathbf{z}_{i}^{T}\beta+\mathbf{x}_{2i}^{T}\delta+\eta{}_{i},
\end{eqnarray*}
with $\mathbb{E}(\epsilon_{i}\vert\mathbf{x}_{2i},\,\mathbf{z}_{i})=\mathbb{E}(\eta_{i}\vert\mathbf{x}_{2i},\,\mathbf{z}_{i})=0.$
Here $y_{i}$, $x_{1i}$, and $\mathbf{x}_{2i}$ denote wage, education
(the endogenous regressor), and a vector of other explanatory variables
(the exogenous regressors) respectively, and $\mathbf{z}_{i}$ denotes
a vector of instrumental variables that have direct effect on education
but are uncorrelated with the unobservables (i.e., $\epsilon_{i}$)
such as innate abilities in the wage equation. 

In many applications, the number of endogenous regressors is also
large relative to the sample size. One example concerns the nonparametric
regression model with endogenous explanatory variables. Consider the
model $y_{i}=f(x_{i})+\epsilon_{i}$ where $\epsilon_{i}\sim\mathcal{N}(0,\,\sigma^{2})$
and $f(\cdot)$ is an unknown function of interest. Assume $\mathbb{E}(\epsilon_{i}\vert X_{i})\neq0$
for all $i$. Suppose we want to approximate $f(x_{i})$ by linear
combinations of some set of basis functions, i.e., $f(x_{i})=\sum_{j=1}^{p}\beta_{j}\phi_{j}(x_{i})$,
where $\{\phi_{1},...,\phi_{p}\}$ are some known functions. Then,
we end up with a linear regression model with many endogenous regressors. 

Empirical examples of many endogenous regressors can be found in hedonic
price regressions of consumer products (e.g., personal computers,
automobiles, pharmaceutical drugs, residential housing, etc.) sold
within a market (say, market \textit{i}) or by a firm (say, firm \textit{i}).
There are two major issues with using firm \textit{i}'s (or, market
\textit{i}'s) product characteristics as the explanatory variables.
First, the number of explanatory variables formed by the characteristics
(and the transformations of these characteristics) of products such
as personal computers, automobiles, and residential houses can be
very large. For example, in the study of hedonic price index analysis
in personal computers, the data considered by Benkard and Bajari involved
65 product characteristics (Benkard and Bajari, 2005). Together with
the various transformations of these characteristics, the number of
the potential regressors can be very large. On the other hand, it
is plausible that only a few of these variables matter to the underlying
prices but which variables constitute the relevant regressors are
unknown to the researchers. Housing data also tends to exhibit a similar
high-dimensional but sparse pattern in terms of the underlying explanatory
variables (e.g., Lin and Zhang, 2006; Ravikumar, et. al, 2009). Second,
firm \textit{i}'s product characteristics are likely to be endogenous
because just like price, product characteristics are typically choice
variables of firms, and it is possible that they are correlated with
unobserved components of price (Ackerberg and Crawford, 2009). An
alternative is to use other firms' (other markets') product characteristics
as the instruments for firm \textit{i}'s (market \textit{i}'s) product
characteristics. In demand estimation literature, this type of instruments
are sometimes referred to as BLP instruments, e.g., Berry, et. al.,
1995 (respectively, Hausman instruments, e.g., Nevo, 2001). 

Another empirical example of many endogenous regressors concerns the
study of network or community influence. For example, Manresa (2013)
looks at how a firm's production output is influenced by the investment
of other firms. As a future extension, she suggests an alternative
model that looks at the network influence in terms of the output of
the other firms rather than their investment: 
\[
y_{it}=\alpha_{i}+\zeta_{t}+\mathbf{x}_{it}^{T}\theta+\sum_{j\in\{1,...,n\},\, j\neq i}\beta_{ji}y_{jt}+\epsilon_{it},\qquad i=1,...,n,\: t=1,...,T
\]
$\mathbf{x}_{it}$ denotes a vector of exogenous regressors specific
to firm \textit{i} (e.g., investment) at period \textit{t}. $\alpha_{i}$
and $\zeta_{t}$ are the fixed effects of firm \textit{i} and period
\textit{t,} respectively. Notice that $y_{jt}$, the output of other
firms enters the right-hand-side of the equations above as additional
regressors and $\beta_{ji}$, $j=1,...,n$, and $j\neq i$ are interpreted
as the network influence arising from other firms' output on firm
\textit{i}'s output. Furthermore, the influence\textit{ }on firm \textit{i
}from firm \textit{j }is allowed to differ from the influence on firm
\textit{j} from firm \textit{i}. Endogeneity arises from the simultaneity
of the output variables when $\textrm{cov}(\epsilon_{it},\,\epsilon_{jt})\neq0$
(e.g., presence of unobserved network characteristics that are common
to all firms). As a result, the number of endogenous regressors in
the model above is of the order $O(n)$, which exceeds the number
of periods \textit{T }in the application considered by Manresa (2013). 

The case of many endogenous regressors and many instrumental variables
has been studied in the context of Generalized Method of Moments by
Fan and Liao (2011), and Gautier and Tsybakov (2011). Fan and Liao
show that the penalized GMM and penalized empirical likelihood are
consistent in both estimation and selection. Gautier and Tsybakov
propose a new estimation procedure called the Self Tuning Instrumental
Variables (STIV) estimator based on the moment conditions $\mathbb{E}(\mathbf{z}_{i}\epsilon_{i})=\mathbf{0}$.
They discuss the STIV procedure with estimated linear projection type
instruments, akin to the 2SLS procedure, and find it works successfully
in simulation. Gautier and Tsybakov also speculate on the rate of
convergence for this type of two-stage estimation procedures when
both stage equations are in the high-dimensional settings. As will
be shown in the subsequent section, the results in this paper partially
confirm their conjecture. 

In the low-dimensional setting, the properties of the 2SLS and GMM
estimators are well-understood. However, it is unclear how the regularized
2SLS procedures compare to the regularized GMM procedures in the high-dimensional
and sparse setting, so it is important to study these regularized
two-stage high-dimensional estimation procedures in depth. Another
important contribution of this paper is to introduce a set of assumptions
that are suitable for showing estimation consistency and selection
consistency of the two-step type of high-dimensional estimators. When
endogeneity is absent from model (1), there is a well-developed theory
on what conditions on the design matrix $X\in\mathbb{R}^{n\times p}$
are sufficient (sufficient and necessary) for an $l_{1}-$based regularized
estimator to consistently estimate (select) $\beta^{*}$. In some
situations one can impose these conditions directly as an assumption
on the underlying design matrix. However, when employing a regularized
2SLS estimator in the context of triangular simultaneous linear equation
models in the high-dimensional setting, namely, (1) and (2), there
is no guarantee that the random matrix $\hat{X}^{T}\hat{X}$ (with
$\hat{X}$ obtained from regressing $X$ on the instrumental variables)
would automatically satisfy these previously established conditions
for estimation or selection consistency. This paper explicitly proves
that these conditions for estimation consistency indeed hold for $\hat{X}^{T}\hat{X}$
with high probability under a broad class of sub-Gaussian design matrices
formed by the instrumental variables allowing for arbitrary correlations
among the covariates. It also establishes the sample size required
for $\hat{X}^{T}\hat{X}$ to satisfy these conditions. Furthermore,
with an additional stronger assumption on the structure of the design
matrices formed by the instrumental variables, this paper shows $\hat{X}^{T}\hat{X}$
also satisfies the conditions for selection consistency under a stronger
sample size requirement. In summary, the aims of this paper, as mentioned
earlier, are to provide a theoretical justification that has not been
given in literature for these regularized 2SLS procedures in the high-dimensional
setting. 

I begin in Section 2 with background on the standard Lasso theory
of high-dimensional estimation techniques as well as basic definitions
and notation used in this paper. Results regarding the estimation
consistency and selection consistency of the high-dimensional 2SLS
procedure under the sparsity scenario are established in Section 3.
Section 4 presents simulation results. Section 5 concludes this paper
and discusses future extensions. All the proofs are collected in the
appendix (Section 6).

\section{Background, notation and definitions}

\textbf{Notation}. For the convenience of the reader, I summarize
here notations to be used throughout this paper. The $l_{q}$ norm
of a vector $v\in m\times1$ is denoted by $|v|_{q}$, $1\leq q\leq\infty$
where $|v|_{q}:=\left(\sum_{i=1}^{m}|v_{i}|^{q}\right)^{1/q}$ when
$1\leq q<\infty$ and $|v|_{q}:=\max_{i=1,...,m}|v_{i}|$ when $q=\infty$.
For a matrix $A\in\mathbb{R}^{m\times m}$, write $|A|_{\infty}:=\max_{i,j}|a_{ij}|$
to be the elementwise $l_{\infty}$- norm of $A$. The $l_{2}$-operator
norm, or spectral norm of the matrix $A$ corresponds to its maximum
singular value; i.e., it is defined as $||A||_{2}:=\sup_{v\in S^{m-1}}|Av|_{2}$,
where $S^{m-1}=\{v\in\mathbb{R}^{m}\,\vert\,|v|_{2}=1\}$. The $l_{\infty}$
matrix norm (maximum absolute row sum) of $A$ is denoted by $||A||_{\infty}:=\max_{i}\sum_{j}|a_{ij}|$
(note the difference between $|A|_{\infty}$ and $||A||_{\infty}$).
I make use of the bound $||A||_{\infty}\leq\sqrt{m}||A||_{2}$ for
any symmetric matrix $A\in\mathbb{R}^{m\times m}$. For a matrix $\Sigma$,
denote its minimum eigenvalue and maximum eigenvalue by $\lambda_{\min}(\Sigma)$
and $\lambda_{\max}(\Sigma)$, respectively. For functions $f(n)$
and $g(n)$, write $f(n)\succsim g(n)$ to mean that $f(n)\geq cg(n)$
for a universal constant $c\in(0,\,\infty)$ and similarly, $f(n)\precsim g(n)$
to mean that $f(n)\leq c^{'}g(n)$ for a universal constant $c^{'}\in(0,\,\infty)$.
$f(n)\asymp g(n)$ when $f(n)\succsim g(n)$ and $f(n)\precsim g(n)$
hold simultaneously. For some integer $s\in\{1,\,2,...,m\}$, the
$l_{0}$-ball of radius $s$ is given by $\mathbb{B}_{0}^{m}(s):=\{v\in\mathbb{R}^{m}\,\vert\,|v|_{0}\leq s\}$
where $|v|_{0}:=\sum_{i=1}^{m}1\{v_{i}\neq0\}$. Similarly, the $l_{2}$-ball
of radius $r$ is given by $\mathbb{B}_{2}^{m}(r):=\{v\in\mathbb{R}^{m}\,\vert\,|v|_{2}\leq r\}$.
Also, write $\mathbb{K}(s,\, m):=\mathbb{B}_{0}^{m}(s)\cap\mathbb{B}_{2}^{m}(1)$
and $\mathbb{K}^{2}(s,\, m):=\mathbb{K}(s,\, m)\times\mathbb{K}(s,\, m)$.
For a vector $v\in\mathbb{R}^{p}$, let $J(v)=\{j\in\{1,...,p\}\,\vert\, v_{j}\neq0\}$
be its support, i.e., the set of indices corresponding to its non-zero
components $v_{j}$. The cardinality of a set $J\subseteq\{1,...,p\}$
is denoted by $|J|$.\\
\\
I will begin with a brief review of the case where all components
in $X$ in (1) are \textit{exogenous}. Assume the number of regressors
$p$ in equation (1) grows with and exceeds the sample size $n$.
Let us focus on the class of models where $\beta^{*}$ has at most
$k$ non-zero parameters, where $k$ is also allowed to increase to
infinity with $n$ but slowly compared to $n$. Consider the following
Lasso program: 
\[
\hat{\beta}_{Las}\in\arg\min_{\beta\in\mathbb{R}^{p}}\left\{ \frac{1}{2n}|y-X\beta|_{2}^{2}+\lambda_{n}|\beta|_{1}\right\} ,
\]
where $\lambda_{n}>0$ is some tuning parameter. Alternatively, we
can consider a constrained version of the Lasso 
\[
\hat{\beta}_{Las}\in\arg\min_{\beta\in\mathbb{R}^{p}}\frac{1}{2n}|Y-X\beta|_{2}^{2}\hspace{1em}\textrm{such that }|\beta|_{1}\leq R.
\]
By Lagrangian duality theory, the above two programs are equivalent.
For example, for any choice of radius $R>0$ in the constrained variant
of the Lasso, there is a tuning parameter $\lambda_{n}(R)\geq0$ such
that solving the Lagrangian form of the Lasso is equivalent to solving
the constrained version.

Consider the constrained Lasso program above with radius $R=|\beta^{*}|_{1}$.
With this setting, the true parameter vector $\beta^{*}$ is feasible
for the problem. By definition, the estimate $\hat{\beta}_{Las}$
minimizes the quadratic loss function $\mathcal{L}(\beta;\,(y,\, X))=\frac{1}{2n}|y-X\beta|_{2}^{2}$
over the $l_{1}-$ball of radius $R$. As $n$ increases, we expect
that $\beta^{*}$ should become a near-minimizer of the same loss,
so that $\mathcal{L}(\hat{\beta}_{Las};\,(y,\, X))\approx\mathcal{L}(\beta^{*};\,(y,\, X))$.
But when does closeness in the loss imply that the error vector $v:=\hat{\beta}_{Las}-\beta^{*}$
is also small? The link between the excess loss $\mathcal{L}(\hat{\beta}_{Las})-\mathcal{L}(\beta^{*})$
and the size of the error $v=\hat{\beta}_{Las}-\beta^{*}$ is the
Hessian of the loss function, $\nabla^{2}L(\beta)=\frac{1}{n}X^{T}X$,
which captures the curvature of the loss function. In the low-dimensional
setting where $p<n$, as long as $\textrm{rank}(X)=p$, we are guaranteed
that the Hessian matrix, $\hat{\Sigma}=\frac{1}{n}X^{T}X$, of the
loss function is positive definite, i.e., $v^{T}\hat{\Sigma}v\geq\delta>0$
for $v\in\mathbb{R}^{p}\backslash\{0\}$. In the high-dimensional
setting with $p>n$, the Hessian is a $p\times p$ matrix with rank
at most $n$, so that it is impossible to guarantee that it has a
positive curvature in all directions. 

The restricted eigenvalue (RE) condition is one of the plausible ways
to relax the stringency of the uniform curvature condition. The RE
condition assumes that the Hessian matrix, $\hat{\Sigma}=\frac{1}{n}X^{T}X$,
of the loss function is positive definite on a restricted set (the
choice of this set is associated with the $l_{1}-$penalty and to
be explained shortly). In the high-dimensional setting, it is well-known
that the RE condition is a sufficient condition for $l_{q}$- consistency
of the Lasso estimate $\hat{\beta}_{Las}$ (see, e.g., Bickel, et.
al., 2009; Meinshausen and Yu, 2009; Raskutti et al., 2010; B�hlmann
and van de Geer, 2011; Loh and Wainwright, 2012; Negahban, et. al.,
2012). In this paper, I will use the following definition (see, Negahban,
et. al., 2012; Wainwright, 2014).\textbf{}\\
\textbf{}\\
\textbf{Definition 1} (RE): The matrix $X\in\mathbb{R}^{n\times p}$
satisfies the RE condition over a subset $S\subseteq\{1,\,2,...,p\}$
with parameter $(\delta,\,\gamma)$ if 
\begin{equation}
\frac{\frac{1}{n}|Xv|_{2}^{2}}{|v|_{2}^{2}}\geq\delta>0\qquad\textrm{for all }v\in\mathbb{C}(S;\,\gamma)\backslash\{\mathbf{0}\},
\end{equation}
where 
\[
\mathbb{C}(S;\,\gamma):=\left\{ v\in\mathbb{R}^{p}\,\vert\,|v_{S^{c}}|{}_{1}\leq\gamma|v_{S}|{}_{1}\right\} \quad\textrm{for some constant \ensuremath{\gamma}}\geq1
\]
with $v_{S}$ denoting the vector in $\mathbb{R}^{p}$ that has the
same coordinates as $v$ on $S$ and zero coordinates on the complement
$S^{c}$ of $S$. \\
\\
When the unknown vector $\beta^{*}\in\mathbb{R}^{p}$ is exactly sparse,
a natural choice of $S$ is the support set of $\beta^{*}$, i.e.,
$J(\beta^{*})$. RE is a weaker condition than other restrictions
in the literature including the pairwise incoherence condition (Donoho,
2006; Gautier and Tsybakov, 2011, Proposition 4.2) and the restricted
isometry property (Cand�s and Tao, 2007). As shown by Bickel et al.,
2009, the restricted isometry property implies the RE condition but
not vice versa. Additionally, Raskutti et al., 2010 give examples
of matrix families for which the RE condition holds, but the restricted
isometry constants tend to infinity as $(n,\,|S|)$ grow. Furthermore,
they show that even when a matrix exhibits a high amount of dependency
among the covariates, it might still satisfy RE. To be more precise,
they show that, if $X\in\mathbb{R}^{n\times p}$ is formed by independently
sampling each row $X_{i}\sim N(0,\,\Sigma)$, then there are strictly
positive constants $(\kappa_{1},\,\kappa_{2})$, depending only on
the positive definite matrix $\Sigma$, such that 
\[
\frac{|Xv|_{2}^{2}}{n}\geq\kappa_{1}|v|_{2}^{2}-\kappa_{2}\frac{\log p}{n}|v|_{1}^{2},\qquad\textrm{for all }v\in\mathbb{R}^{p},
\]
with probability at least $1-c_{1}\exp(-c_{2}n)$ for some universal
constants $c_{1}$ and $c_{2}$. The bound above ensures the RE condition
holds with $\delta=\frac{\kappa_{1}}{2}$ and $\gamma=3$ as long
as $n>32\frac{\kappa_{2}}{\kappa_{1}}k\log p$. To see this, note
that for any $v\in\mathbb{C}(J(\beta^{*}),\,3)$, we have $|v|_{1}^{2}\leq16|v_{J(\beta^{*})}|_{1}^{2}\leq16k|v_{J(\beta^{*})}|_{2}^{2}$.
Given the lower bound above, for any $v\in\mathbb{C}(J(\beta^{*});\,3)$,
we have the lower bound 
\[
\frac{|Xv|_{2}^{2}}{n}\geq\left(\kappa_{1}-16\kappa_{2}\frac{k\log p}{n}\right)|v|_{2}^{2}\geq\frac{\kappa_{1}}{2}|v|_{2}^{2},
\]
where the final inequality follows as long as $n>32(\frac{\kappa_{2}}{\kappa_{1}})^{2}k\log p$.
An appropriate choice of the tuning parameter $\lambda_{n}$ in the
Lasso program ensures $\hat{v}:=\hat{\beta}_{Las}-\beta^{*}\in\mathbb{C}(J(\beta^{*});\,3)$.
This fact can be formalized in the following proposition. \\
\\
\textbf{Proposition 2.1}. For the linear model $y_{i}=\mathbf{x}_{i}\beta^{*}+\epsilon_{i}$
where $\mathbb{E}(\mathbf{x}_{i}\epsilon_{i})=\mathbf{0}$, if we
solve the Lasso program with parameter $\lambda_{n}\geq\frac{2}{n}||X^{T}\epsilon||_{\infty}>0$,
then the error $\hat{v}:=\hat{\beta}_{Las}-\beta^{*}\in\mathbb{C}(J(\beta^{*}),\,3)$.
\\
\\
\textbf{Proof}. Define the Lagrangian $L(\beta;\,\lambda_{n})=\frac{1}{2n}|y-X\beta|_{2}^{2}+\lambda_{n}|\beta|_{1}$.
Since $\hat{\beta}_{Las}$ is optimal, we have 
\[
L(\hat{\beta};\,\lambda_{n})\leq L(\beta^{*};\,\lambda_{n})=\frac{1}{2n}|\epsilon|_{2}^{2}+\lambda_{n}|\beta^{*}|_{1},
\]
Some algebraic manipulation of the \textit{basic inequality} above
yields 
\begin{eqnarray*}
0 & \leq & \frac{1}{2n}|X\hat{v}|_{2}^{2}\leq\frac{1}{n}\epsilon X\hat{v}+\lambda_{n}\left\{ |\beta_{J(\beta^{*})}^{*}|_{1}-|(\beta_{J(\beta^{*})}^{*}+\hat{v}_{J(\beta^{*})},\,\hat{v}_{J(\beta^{*})^{c}})|_{1}\right\} \\
 & \leq & |\hat{v}|_{1}|\frac{1}{n}X^{T}\epsilon|_{\infty}+\lambda_{n}\left\{ |\hat{v}_{J(\beta^{*})}|_{1}-|\hat{v}_{J(\beta^{*})^{c}}|_{1}\right\} \\
 & \leq & \frac{\lambda_{n}}{2}\left\{ 3|\hat{v}_{J(\beta^{*})}|_{1}-|\hat{v}_{J(\beta^{*})^{c}}|_{1}\right\} ,
\end{eqnarray*}
where the last line applies the assumption on $\lambda_{n}$. $\square$

Rudelson and Zhou (2011) as well as Loh and Wainwright (2012) extend
this type of RE analysis from the case of Gaussian designs to the
case of sub-Gaussian designs. The sub-Gaussian assumption says that
the explanatory variables need to be drawn from distributions with
well-behaved tails like Gaussian. In contrast to the Gaussian assumption,
sub-Gaussian variables constitute a more general family of distributions.
In this paper, I make use of the following definition for a sub-Gaussian
matrix. \\
\textbf{}\\
\textbf{Definition 2}: A random variable $X$ with mean $\mu=\mathbb{E}[X]$
is sub-Gaussian if there is a positive number $\sigma$ such that
\[
\mathbb{E}[\exp(t(X-\mu))]\leq\exp(\sigma^{2}t^{2}/2)\qquad\textrm{for all}\, t\in\mathbb{R},
\]
and a random matrix $A\in\mathbb{R}^{n\times p}$ is sub-Gaussian
with parameters $(\Sigma_{A},\,\sigma_{A}^{2})$ if (a) each row $A_{i}^{T}\in\mathbb{R}^{p}$
is sampled independently from a zero-mean distribution with covariance
$\Sigma_{A}$, (b) for any unit vector $u\in\mathbb{R}^{p}$, the
random variable $u^{T}A_{i}^{T}$ is sub-Gaussian with parameter at
most $\sigma_{A}^{2}$.\\
\textbf{}\\
For example, if $A\in\mathbb{R}^{n\times p}$ is formed by independently
sampling each row $A_{i}\sim N(0,\,\Sigma_{A})$, then the resulting
matrix\textbf{ $A\in\mathbb{R}^{n\times p}$ }is a sub-Gaussian matrix
with parameters \textbf{$(\Sigma_{A},\,||\Sigma_{A}||_{2})$}, recalling
\textbf{$||\Sigma_{A}||_{2}$ }denotes the spectral norm of $\Sigma_{A}$.

\section{High-dimensional 2SLS estimation}

Suppose from performing a first-stage regression on each of the equations
in (2) separately, we obtain estimates $\hat{\pi}_{j}$ and let $\hat{\mathbf{x}}_{j}:=Z_{j}\hat{\pi}_{j}$
for $j=1,...,p$. Denote the fitted regressors from the first-stage
estimation by $\hat{X}$, where $\hat{X}=\left(\hat{\mathbf{x}}_{1},...,\hat{\mathbf{x}}_{p}\right)$.
For the second-stage regression, consider the following Lasso program:
\begin{equation}
\hat{\beta}_{H2SLS}\in\textrm{argmin}_{\beta\in\mathbb{R}^{p}}:\:\frac{1}{2n}|y-\hat{X}\beta|_{2}^{2}+\lambda_{n}|\beta|_{1}.
\end{equation}
The following is a standard assumption in the literature on sparsity
for high-dimensional linear models. \\
\\
\textbf{Assumption 3.1}: The numbers of regressors $p(=p_{n})$ and
$d_{j}(=d_{jn})$ for every $j=1,...,p$ in (1) and (2) can grow with
and exceed the sample size $n$.\textbf{ }The number of non-zero components
in $\pi_{j}^{*}$ is at most $k_{1}(=k_{1n})$ for all $j=1,...,p$,
and the number of non-zero components in $\beta^{*}$ is at most $k_{2}(=k_{2n})$.
Both $k_{1}$ and $k_{2}$ can increase to infinity with $n$ but
slowly compared to $n$.\\
\\
I first present a general bound on the statistical error measured
by the quantity $|\hat{\beta}_{H2SLS}-\beta^{*}|_{2}$. \\
\\
\textbf{Lemma 3.1} (General upper bound on the $l_{2}-$error). Let
$\hat{\Gamma}=\frac{1}{n}\hat{X}^{T}\hat{X}$ and $e=(X-\hat{X})\beta^{*}+\boldsymbol{\eta}\beta^{*}+\epsilon$.
Suppose the random matrix $\hat{\Gamma}$ satisfies the RE condition
(3) with $\gamma=3$ and the vector $\beta^{*}$ is supported on a
subset $J(\beta^{*})\subseteq\{1,\,2,\,...p\}$ with its cardinality
$|J(\beta^{*})|\leq k_{2}$. If a solution $\hat{\beta}_{H2SLS}$,
defined in (4) has $\lambda_{n}$ satisfying 
\[
\lambda_{n}\geq2|\frac{1}{n}\hat{X}^{T}e|_{\infty}>0,
\]
for any given $n$, then there is a constant $c>0$ such that 
\[
|\hat{\beta}_{H2SLS}-\beta^{*}|_{2}\leq\frac{c}{\delta}\sqrt{k_{2}}\lambda_{n}.
\]

The proof for Lemma 3.1 is provided in Section 6.1.

Notice that the choice of $\lambda_{n}$ in Lemma 3.1 depends on unknown
quantities and therefore Lemma 3.1 does not provide guidance to practical
implementation. Rather, it should be viewed as an intermediate lemma
for proving consistency of the two-stage estimator later on. In the
appendix (Section 6) we show that the term $|\frac{1}{n}\hat{X}^{T}e|_{\infty}$
can be bounded from above and the order of the resulting upper bound
can be used to set the tuning parameter $\lambda_{n}$. In order to
apply Lemma 3.1 to prove consistency, we need to show (i) $\hat{\Gamma}=\frac{1}{n}\hat{X}^{T}\hat{X}$
satisfies the RE condition (3) with $\gamma=3$ and (ii) the term
$|\frac{1}{n}\hat{X}^{T}e|_{\infty}\precsim f(k_{1},\, k_{2},\, d_{1},...,d_{p},\, p,\, n)$
with high probability, and then we can show 
\[
|\hat{\beta}_{H2SLS}-\beta^{*}|_{2}\precsim\sqrt{k_{2}}f(k_{1},\, k_{2},\, d_{1},...,d_{p},\, p,\, n)
\]
by choosing $\lambda_{n}\asymp f(k_{1},\, k_{2},\, d_{1},...,d_{p},\, p,\, n)$.
The assumption $\sqrt{k_{2}}f(k_{1},\, k_{2},\, d_{1},...d_{p},\, p,\, n)=o(1)$
will therefore imply the $l_{2}$-consistency of $\hat{\beta}_{H2SLS}$.
Applying Lemma 3.1 to the triangular simultaneous equations model
(1) and (2) requires additional work to establish conditions (i) and
(ii) discussed above, which depends on the specific first-stage estimator
for $\hat{X}$. It is worth mentioning that, while in many situations
one can impose the RE condition as an assumption on the design matrix
(e.g., Belloni and Chernozhukov, 2011b; Belloni, Chen, Chernozhukov,
and Hansen, 2012) in analyzing the consistency property of the Lasso,
appropriate analysis is needed in this paper to verify that $\frac{1}{n}\hat{X}^{T}\hat{X}$
satisfies the RE condition because $\hat{X}$ is obtained from a first-stage
estimation and there is no guarantee that the random matrix $\frac{1}{n}\hat{X}^{T}\hat{X}$
would automatically satisfy the RE condition. To the best of my knowledge,
previous literature has not dealt with this issue directly. Consequently,
the RE analysis introduced in this paper is particularly useful for
analyzing the statistical properties of the two-step type of high-dimensional
estimators in the simultaneous equations model context. As discussed
previously, this paper focuses on the case where $p\gg n$ and $d_{j}\ll n$
for all $j$ and the case where $p\gg n$ and $d_{j}\gg n$ for at
least one $j$. The following two subsections present results concerning
estimation consistency and variable-selection consistency for the
exact sparsity case.

\subsection{Estimation consistency for the sparsity case }

To derive the non-asymptotic bounds and asymptotic properties (i.e.,
estimation consistency and selection consistency) for $\hat{\beta}_{H2SLS}$,
I impose the following regularity conditions.\\
\textbf{}\\
\textbf{Assumption 3.2}: The error terms $\epsilon$ and $\eta_{j}$
for $j=1,...,p$ are \textit{i.i.d.} zero-mean sub-Gaussian vectors
with parameters $\sigma_{\epsilon}^{2}$ and $\sigma_{\eta}^{2}$,
respectively.\textbf{ }The random matrix\textbf{ $Z_{j}\in\mathbb{R}^{n\times d_{j}}$
}is\textbf{ }sub-Gaussian with parameters $(\Sigma_{Z_{j}},\,\sigma_{Z}^{2})$
for $j=1,...,p$.\textbf{}\\
\textbf{}\\
\textbf{Assumption 3.3}: For every $j=1,...,p$, $\mathbf{x}_{j}^{*}:=Z_{j}\pi_{j}^{*}$.
The matrix $X^{*}\in\mathbb{R}^{n\times p}$ is sub-Gaussian with
parameters $(\Sigma_{X^{*}},\,\sigma_{X^{*}}^{2})$ where the\textbf{
$j$}th column of $X^{*}$ is $\mathbf{x}_{j}^{*}$. \textbf{}\\
\textbf{}\\
\textbf{Assumption 3.4}: For every $j=1,...,p$, $\mathbf{w}_{j}:=Z_{j}v_{j}$
where $v_{j}\in\mathbb{K}(k_{1},\, d_{j}):=\mathbb{B}_{0}^{d_{j}}(k_{1})\cap\mathbb{B}_{2}^{d_{j}}(1)$.
The matrix $W\in\mathbb{R}^{n\times p}$ is sub-Gaussian with parameters
$(\Sigma_{W},\,\sigma_{W}^{2})$ where the\textbf{ $j$}th column
of $W$ is $\mathbf{w}_{j}$.\textbf{}\\
\textbf{}\\
\textbf{Assumption 3.5a}: The first-stage estimator $\hat{\pi}^{T}\in\mathbb{R}^{p\times d}$
satisfies the bound $\max_{j=1,...,p}|\hat{\pi}_{j}-\pi_{j}^{*}|_{1}\leq\frac{c\sigma_{\eta}}{\lambda_{\min}(\Sigma_{Z})}k_{1}\sqrt{\frac{\log\max(d,\, p)}{n}}$
with probability at least $1-c_{1}\exp(-c_{2}\log\max(d,\, p,\, n))$
for some universal constants $c_{1}$ and $c_{2}$, where\textbf{
}$d=\max_{j=1,...,p}d_{j}$ and $\lambda_{\min}(\Sigma_{Z})=\min_{j=1,...,p}\lambda_{\min}(\Sigma_{Z_{j}})$.
\textbf{}\\
\textbf{}\\
\textbf{Assumption 3.5b}: The first-stage estimator $\hat{\pi}^{T}\in\mathbb{R}^{p\times d}$
satisfies the bound \textbf{$\max_{j=1,...,p}|\hat{\pi}_{j}-\pi_{j}^{*}|_{2}\leq\frac{c\sigma_{\eta}}{\lambda_{\min}(\Sigma_{Z})}\sqrt{\frac{k_{1}\log\max(d,\, p)}{n}}$}
with probability at least $1-c_{1}\exp(-c_{2}\log\max(d,\, p,\, n))$
for some universal constants $c_{1}$ and $c_{2}$, where\textbf{
}$d=\max_{j=1,...,p}d_{j}$ and $\lambda_{\min}(\Sigma_{Z})=\min_{j=1,...,p}\lambda_{\min}(\Sigma_{Z_{j}})$.
\textbf{}\\
\textbf{}\\
\textbf{Assumption 3.6}: For every $j=1,...,p$, the first-stage estimator
$\hat{\pi}_{j}$ achieves the selection consistency (i.e., it recovers
the true support $J(\pi_{j}^{*})$) or has at most $k_{j}^{*}$ components
that are different from the components in $J(\pi_{j}^{*})$ where
$k_{j}^{*}\ll n$, with probability at least $1-c_{1}\exp(-c_{2}\log\max(d,\, p,\, n))$
for some universal constants $c_{1}$ and $c_{2}$, where\textbf{
}$d=\max_{j=1,...,p}d_{j}$. For simplicity, we consider the case
where the first-stage estimator recovers the true support $J(\pi_{j}^{*})$
for every $j=1,...,p$.\textbf{}\\
\textbf{}\\
\textbf{Remarks}\\
Assumption 3.2 is common in the literature (see, Loh and Wainwright,
2012; Negahban, et. al 2012; Rosenbaum and Tsybakov, 2013). The assumption
that \textbf{$Z_{j}\in\mathbb{R}^{n\times d_{j}}$ }is\textbf{ }sub-Gaussian
with parameters $(\Sigma_{Z_{j}},\,\sigma_{Z}^{2})$ for all $j$
provides a primitive condition which guarantees that the random matrix
formed by the instrumental variables satisfies the RE condition with
high probability. 

Based on the second part of Assumption 3.2 that \textbf{$Z_{j}\in\mathbb{R}^{n\times d_{j}}$
}is\textbf{ }sub-Gaussian with parameters $(\Sigma_{Z_{j}},\,\sigma_{Z}^{2})$
for all $j$, we have that $Z_{j}\pi_{j}^{*}:=\mathbf{x}_{j}^{*}$
and $Z_{j}v_{j}:=\mathbf{w}_{j}$ are sub-Gaussian vectors where $v_{j}\in\mathbb{K}(k_{1},\, d_{j}):=\mathbb{B}_{0}^{d_{j}}(k_{1})\cap\mathbb{B}_{2}^{d_{j}}(1)$.
Therefore, the conditions that $X^{*}\in\mathbb{R}^{n\times p}$ is
a sub-Gaussian matrix with parameters $(\Sigma_{X^{*}},\,\sigma_{X^{*}}^{2})$
where the\textbf{ $j$}th column of $X^{*}$ is $\mathbf{x}_{j}^{*}$(Assumption
3.3) and $W\in\mathbb{R}^{n\times p}$ is a sub-Gaussian matrix with
parameters $(\Sigma_{W},\,\sigma_{W}^{2})$ where the\textbf{ $j$}th
column of $W$ is $\mathbf{w}_{j}$ (Assumption 3.4) are mild extensions.
In terms of the instrumental variables and their linear combinations,
Assumptions 3.2-3.4 together with Assumption 3.5a (or 3.5b) on the
first-stage estimation error provide primitive conditions which guarantee
that the random matrix $\frac{1}{n}\hat{X}^{T}\hat{X}$ formed by
the fitted regressors $\hat{\mathbf{x}}_{j}(:=Z_{j}\hat{\pi}_{j})$
for $j=1,...,p$ satisfies the RE condition with high probability.

For Assumptions 3.5a(b), many existing high-dimensional estimation
procedures such as the Lasso or Dantzig selector (see, e.g., Cand�s
and Tao, 2007; Bickel, et. al, 2009; Negahban, et. al. 2012) simultaneously
satisfy the error bounds $\max_{j=1,...,p}|\hat{\pi}_{j}-\pi_{j}^{*}|_{1}\leq\frac{c\sigma_{\eta}}{\lambda_{\min}(\Sigma_{Z})}k_{1}\sqrt{\frac{\log\max(d,\, p)}{n}}$
(Assumption 3.5a) and $\max_{j=1,...,p}|\hat{\pi}_{j}-\pi_{j}^{*}|_{2}\leq\frac{c\sigma_{\eta}}{\lambda_{\min}(\Sigma_{Z})}\sqrt{\frac{k_{1}\log\max(d,\, p)}{n}}$
(Assumption 3.5b) with high probability. The reason I introduce Assumptions
3.5a and 3.5b separately will be explained shortly. It is worth noting
that while the $l_{1}-$error ($l_{2}-$error) from applying the Lasso
on a single first-stage equation should be of the order $O\left(k_{1}\sqrt{\frac{\log d}{n}}\right)$
(respectively, $O\left(\sqrt{\frac{k_{1}\log d}{n}}\right)$) with
probability at least $1-c_{1}\exp(-c_{2}\log\max(d,\, n))$, the extra
term $\log p$ in the errors $\max_{j=1,...,p}|\hat{\pi}_{j}-\pi_{j}^{*}|_{1}$
and $\max_{j=1,...,p}|\hat{\pi}_{j}-\pi_{j}^{*}|_{2}$ and the probability
guarantee $1-c_{1}\exp(-c_{2}\log\max(d,\, p,\, n))$ with which these
errors hold comes from the application of a union bound which takes
into account the fact that there are $p$ endogenous regressors in
the main equation and hence, $p$ equations to estimate in the first-stage.
As a result, it is not hard to see that the sample size required for
consistently estimating $p$ equations simultaneously when a Lasso-type
procedure is applied on each of the first-stage equations separately
should satisfy $\sqrt{\frac{k_{1}\log\max(d,\, p)}{n}}=o(1)$ as opposed
to the condition $\sqrt{\frac{k_{1}\log d}{n}}=o(1)$ for the case
where a single equation is estimated with a Lasso-type procedure. 

Assumption 3.6 says that the first-stage estimators correctly select
the non-zero coefficients with probability close to 1. In analogy
to the various sparsity assumptions on the true parameters in the
high-dimensional statistics literature (including the case of \textit{exact
sparsity} assumption meaning that the true parameter vector has only
a few non-zero components, or \textit{approximate sparsity} assumption
based on imposing a certain decay rate on the ordered entries of the
true parameter vector), Assumption 3.6 can be interpreted as an \textit{exact
sparsity} constraint on the first-stage estimate\textbf{ $\hat{\pi}_{j}$
}for\textbf{ $j=1,...,p$}, in terms of the $l_{0}$- ball, given
by 
\[
\mathbb{B}_{0}^{d_{j}}(k_{1}):=\left\{ \hat{\pi}_{j}\in\mathbb{R}^{d_{j}}\,\vert\,\sum_{l=1}^{d_{j}}1\{\hat{\pi}_{jl}\neq0\}\leq k_{1}\right\} \,\textrm{for }j=1,...,p.
\]
It is known that under some stringent conditions such as the {}``irrepresentable
condition'' (Zhao and Yu, 2006; B�hlmann and van de Geer, 2011) or
the {}``mutual incoherence condition'' (Wainwright, 2009) together
with the {}``beta-min condition'' (B�hlmann and van de Geer, 2011),
Lasso and Dantzig types of selectors can recover the support of the
true parameter vector with high probability. The {}``irrepresentable
condition'', as discussed in B�hlmann and van de Geer, 2011, is in
fact a sufficient and necessary condition to achieve variable-selection
consistency with the Lasso. Furthermore, they show that the {}``irrepresentable
condition'' implies the RE condition. Assumption 3.6 is the key condition
that differentiates the upper bounds in the two theorems to be presented
immediately. Similar to the problem of estimating $p$ equations as
in the discussion of Assumptions 3.5a(b), the sample size required
for consistently selecting the coefficients in each of the $p$ equations
simultaneously when a Lasso-type selector is applied on each of the
first-stage equations separately should satisfy $\sqrt{\frac{k_{1}\log\max(d,\, p)}{n}}=O(1)$
as opposed to the condition $\sqrt{\frac{k_{1}\log d}{n}}=O(1)$ for
the case where a single equation is estimated with a Lasso-type selector.
In addition, the {}``beta-min'' condition for consistent selection
in the $p$-equation problem needs to satisfy $\min_{j=1,...,p}\min_{l\in J(\pi_{j}^{*})}|\pi_{jl}^{*}|\geq O\left(\sqrt{\frac{\log\max(d,\, p)}{n}}\right)$
as opposed to $\min_{j=1,...,p}\min_{l\in J(\pi_{j}^{*})}|\pi_{jl}^{*}|\geq O\left(\sqrt{\frac{\log d}{n}}\right)$
for the consistent selection in a single equation problem. 

First, I present two results for the case where $p\gg n$ and $d_{j}\gg n$
for at least one $j$. As discussed earlier, the key difference between
the two theorems is that the bound in the second theorem hinges on
the additional assumption that the first-stage estimators correctly
select the non-zero coefficients with probability close to 1, i.e.,
Assumption 3.6. With this assumption, when the first-stage estimation
error dominates the second-stage error, the statistical error of the
parameters of interests in the main equation can be bounded by the
first-stage estimation error in $l_{2}-$norm. However, without Assumption
3.6, the statistical error of the parameters in the main equation
needs to be bounded by the first-stage estimation error in $l_{1}-$norm.
\textbf{}\\
\\
\textbf{Theorem 3.2} (Upper bound on the $l_{2}-$error and estimation
consistency): Suppose Assumptions 1.1, 3.1-3.3, and 3.5a hold. Then,
if %
\footnote{If the term $o(1)$ in {}``$\frac{k_{1}^{2}\log\max(d,\, p)}{n}=o(1)$''
(similarly, $o(1)$ in {}``$\frac{k_{1}\log\max(d,\, p)}{n}=o(1)$''
in Theorem 3.3, $o(1)$ in {}``$\frac{\log p}{n}=o(1)$'' in Corollary
3.4, $o(1)$ in {}``$\max\left\{ k_{1}M^{2}(d,\, p,\, k_{1},\, n),\:\frac{\log p}{n}\right\} =o(1)$''
in Theorem 3.5, and $o(1)$ in {}``$\max\left\{ M^{2}(d,\, p,\, k_{1},\, n),\:\frac{\log p}{n}\right\} =o(1)$''
in Theorem 3.6) is replaced by $O(1)$, the statistical error of the
parameters in the main equation will have the same scaling in terms
of $d$, $p$, $k_{1}$, $k_{2}$, and $n$ as before with the only
changes to the constants in $\varphi_{1}$ and $\varphi_{2}$.%
}
\begin{eqnarray*}
\frac{k_{1}^{2}k_{2}^{2}\log\max(d,\, p)}{n} & = & O(1),\\
\frac{k_{1}^{2}\log\max(d,\, p)}{n} & = & o(1),
\end{eqnarray*}
and the tuning parameter $\lambda_{n}$ satisfies 
\[
\lambda_{n}\asymp k_{1}k_{2}\sqrt{\frac{\log\max(d,\, p)}{n}},
\]
we have 
\[
|\hat{\beta}_{H2SLS}-\beta^{*}|_{2}\precsim\max\{\varphi_{1}\sqrt{k_{1}k_{2}}\sqrt{\frac{k_{1}\log\max(d,\, p)}{n}},\,\varphi_{2}\sqrt{\frac{k_{2}\log p}{n}}\},
\]
where 
\begin{eqnarray*}
\varphi_{1} & = & \frac{\sigma_{\eta}\max_{j,j^{'}}|\textrm{cov}(x_{1j^{'}}^{*},\,\mathbf{z}_{1j})|_{\infty}|\beta^{*}|_{1}}{\lambda_{\min}(\Sigma_{Z})\lambda_{\min}(\Sigma_{X^{*}})},\\
\varphi_{2} & = & \max\left\{ \frac{\sigma_{X^{*}}\sigma_{\eta}|\beta^{*}|_{1}}{\lambda_{\min}(\Sigma_{X^{*}})},\,\frac{\sigma_{X^{*}}\sigma_{\epsilon}}{\lambda_{\min}(\Sigma_{X^{*}})}\right\} ,
\end{eqnarray*}
with probability at least $1-c_{1}\exp(-c_{2}\log\max(\min(p,\, d),\, n))$
for some universal positive constants $c_{1}$ and $c_{2}$. If we
also have $k_{2}k_{1}\sqrt{\frac{k_{2}\log\max(d,\, p)}{n}}=o(1)$,
then%
\footnote{The extra factor $k_{2}$ in front of these scaling conditions for
consistency in Theorem 3.2 (as well as in the subsequent theorems
3.3, 3.5, 3.6, and Corollary 3.4) comes from the simple inequality
$|\beta^{*}|_{1}\leq k_{2}\max_{j}\beta_{j}^{*}$. %
} the two-stage estimator $\hat{\beta}_{H2SLS}$ is $l_{2}-$consistent
for \textbf{$\beta^{*}$}.\\
\\
\textbf{Theorem 3.3} (An improved upper bound on the $l_{2}-$error
and estimation consistency): Suppose Assumptions 1.1, 3.1-3.4, 3.5b,
and 3.6 hold. Then, if 
\begin{eqnarray*}
\frac{1}{n}\min\left\{ k_{1}^{2}k_{2}^{2}\log\max(d,\, p),\,\min_{r\in[0,\,1]}\max\left\{ k_{1}^{3-2r}\log d,\, k_{1}^{3-2r}\log p,\, k_{1}^{r}k_{2}\log d,\, k_{1}^{r}k_{2}\log p\right\} \right\}  & = & O(1)\\
\frac{k_{1}\log\max(d,\, p)}{n} & = & o(1),
\end{eqnarray*}
and the tuning parameter $\lambda_{n}$ satisfies 
\[
\lambda_{n}\asymp k_{2}\sqrt{\frac{k_{1}\log\max(d,\, p)}{n}},
\]
we have, 
\[
|\hat{\beta}_{H2SLS}-\beta^{*}|_{2}\precsim\max\{\varphi_{1}\sqrt{k_{2}}\sqrt{\frac{k_{1}\log\max(d,\, p)}{n}},\,\varphi_{2}\sqrt{\frac{k_{2}\log p}{n}}\},
\]
with probability at least $1-c_{1}\exp(-c_{2}\log\max(\min(p,\, d),\, n))$
for some universal positive constants $c_{1}$ and $c_{2}$, where
$\varphi_{1}$ and $\varphi_{2}$ are defined in Theorem 3.2. If we
also have $k_{2}\sqrt{\frac{k_{1}k_{2}\log\max(d,\, p)}{n}}=o(1)$,
then the two-stage estimator $\hat{\beta}_{H2SLS}$ is $l_{2}-$consistent
for \textbf{$\beta^{*}$}.\\

The proofs for Theorems 3.2 and 3.3 are provided in Sections 6.2 and
6.3, respectively.

The proofs for Theorems 3.2 and 3.3 each consist of two parts. The
first part is to show $\frac{1}{n}\hat{X}^{T}\hat{X}$ satisfies the
RE condition (3) and the second part is to bound the term $|\frac{1}{n}\hat{X}^{T}e|_{\infty}$
from above. Based on Lemma 3.1, the upper bound on $|\frac{1}{n}\hat{X}^{T}e|_{\infty}$
pins down the scaling requirement of $\lambda_{n}$, as mentioned
previously. The scaling conditions of $n$ and $\lambda_{n}$ depend
on the sparsity parameters $k_{1}$ and $k_{2}$, which are typically
unknown. Nevertheless, I will assume that upper bounds on $k_{1}$
and $k_{2}$ are available, i.e., we know that $k_{1}\leq\bar{k}_{1}$
and $k_{2}\leq\bar{k}_{2}$ for some integers $\bar{k}_{1}$ and $\bar{k}_{2}$
that grow with $n$ just like $k_{1}$ and $k_{2}$. Meaningful values
of $\bar{k}_{1}$ and $\bar{k}_{2}$ are small relative to $n$ presuming
that only a few regressors are relevant. This type of upper bound
assumption on the sparsity is called \textit{sparsity certificate}
in the literature (see, e.g., Gautier and Tsybakov, 2011). 

In Theorems 3.2 and 3.3, we see that the statistical errors of the
parameters of interests in the main equation depend on $\sigma_{\eta}$,
$\sigma_{\epsilon}$, $\sigma_{X^{*}}$, $\lambda_{\min}(\Sigma_{Z})$,
$\lambda_{\min}(\Sigma_{X^{*}})$, and $\max_{j,j^{'}}|\textrm{cov}(x_{1j^{'}}^{*},\,\mathbf{z}_{1j})|_{\infty}$.
In the simple case of $\sigma_{\eta}=0$ (for example, $\boldsymbol{\eta}=\mathbf{0}$
with probability 1 as in a high-dimensional linear regression model
without endogeneity), the $l_{2}-$errors in Theorems 3.2 and 3.3
reduce to $|\hat{\beta}_{H2SLS}-\beta^{*}|_{2}\precsim\frac{\sigma_{X^{*}}\sigma_{\epsilon}}{\lambda_{\min}(\Sigma_{X^{*}})}\sqrt{\frac{k_{2}\log p}{n}}$,
where the factor $\frac{\sigma_{X^{*}}\sigma_{\epsilon}}{\lambda_{\min}(\Sigma_{X^{*}})}$
has a natural interpretation of an inverse signal-to-noise ratio.
For instance, when $X^{*}$ is a zero-mean Gaussian matrix with covariance
$\Sigma_{X^{*}}=\sigma_{X^{*}}^{2}I$, one has $\lambda_{\min}(\Sigma_{X^{*}})=\sigma_{X^{*}}^{2}$,
so 
\[
\frac{\sigma_{X^{*}}\sigma_{\epsilon}}{\lambda_{\min}(\Sigma_{X^{*}})}=\frac{\sigma_{\epsilon}}{\sigma_{X^{*}}},
\]
which measures the inverse signal-to-noise ratio of the regressors
in a high-dimensional linear regression model without endogeneity.
Hence, the statistical error of the parameters of interests in the
main equation matches the scaling of the upper bound for the Lasso
in the context of the high-dimensional linear regression model without
endogeneity, i.e., $\sqrt{\frac{k_{2}\log p}{n}}$.

The terms $\max_{j,j^{'}}|\textrm{cov}(x_{1j^{'}}^{*},\,\mathbf{z}_{1j})|_{\infty}$
in Theorems 3.2 and 3.3 are related to the degree of dependency between
the columns of the design matrices formed by the instrumental variables
and their linear combinations. For instance, for any $l=1,...,d_{j}$
and $j=1,...,p$, notice that 
\[
\textrm{cov}(x_{1j}^{*},\, z_{1jl})=\textrm{cov}(\mathbf{z}_{1j}\pi_{j}^{*},\, z_{1jl}).
\]
The higher dependency between the columns of the design matrix $Z_{j}$
we have, the greater $\max_{l}\textrm{cov}(x_{1j}^{*},\, z_{1jl})$
is, and the harder the estimation problem becomes. In the special
case of $\max_{j,j^{'}}|\textrm{cov}(x_{1j^{'}}^{*},\,\mathbf{z}_{1j})|_{\infty}=\sigma_{Z}^{2}$,
$\lambda_{\min}(\Sigma_{Z})=\sigma_{Z}^{2}$, and $\lambda_{\min}(\Sigma_{X^{*}})=\sigma_{X^{*}}^{2}$,
$\varphi_{1}=\frac{\sigma_{\eta}}{\sigma_{X^{*}}^{2}}|\beta^{*}|_{1}=\frac{1}{\sigma_{X^{*}}}\left(\frac{\sigma_{\eta}}{\sigma_{X^{*}}}\right)|\beta^{*}|_{1}$,
where the multiplier $\frac{1}{\sigma_{X^{*}}}\left(\frac{\sigma_{\eta}}{\sigma_{X^{*}}}\right)$
in $\varphi_{1}$ is the inverse signal-to-noise ratio of $X^{*}$
scaled by $\frac{1}{\sigma_{X^{*}}}$. 

Under the assumption that the first-stage estimators correctly select
the non-zero coefficients with high probability (Assumption 3.6),
the scaling of the sample size required\textbf{ }in Theorem 3.3 is
guaranteed to be no greater (and in some cases strictly smaller) than
that in Theorem 3.2. For instance, if $p\leq d$ , then letting $r=1$
yields 
\[
\max\left\{ k_{1}\log d,\, k_{1}\log p,\, k_{1}k_{2}\log d,\, k_{1}k_{2}\log p\right\} =k_{1}k_{2}\log d\leq k_{1}^{2}k_{2}^{2}\log\max(d,\, p)=k_{1}^{2}k_{2}^{2}\log d.
\]
In this example, Theorem 3.2 suggests that the choice of sample size
needs to satisfy $\frac{k_{1}^{2}k_{2}^{2}\log d}{n}=O(1)$ and $\frac{k_{1}^{2}\log d}{n}=o(1)$
while Theorem 3.3 suggests that the choice of sample size only needs
to satisfy $\frac{k_{1}k_{2}\log d}{n}=O(1)$ and $\frac{k_{1}\log d}{n}=o(1)$.

From Theorem 3.2 (respectively, Theorem 3.3), we see that the estimation
error of the parameters of interests in the main equation is of the
order of the maximum of the first-stage estimation error in $l_{2}-$norm
multiplied by a factor of $\sqrt{k_{1}k_{2}}$ (respectively, $\sqrt{k_{2}}$)
and the second-stage estimation error. Upon the additional condition
that the first-stage estimators correctly select the non-zero coefficients
with probability close to 1, note that the bound on the $l_{2}-$error
of $\hat{\beta}_{H2SLS}$ in Theorem 3.3 is improved upon that in
Theorem 3.2 by a factor of $\sqrt{k_{1}}$\textbf{ }if the first term
in the braces dominates the second one. It is possible that the error
bound and scaling of the sample size required in Theorem 3.2 is suboptimal.
Section 5 provides a heuristic argument that may potentially improve
the bound on the $l_{2}-$error of $\hat{\beta}_{H2SLS}$ in Theorem
3.2 when the first-stage estimates fail to satisfy the exact sparsity
constraint specified by the $l_{0}$- ball discussed earlier. Intuitively,
the most direct effect on the $l_{2}-$error of the second-stage estimate
$\hat{\beta}_{H2SLS}$ should be attributed to the $l_{2}-$errors
(rather than the selection performance per se) of the first-stage
estimates. Imposing the exact sparsity constraint, namely, selection
consistency on the first-stage estimates such as Assumption 3.6 is
an example of showing how special structures that impose a certain
decay rate on the ordered entries of the first-stage estimates from
the $l_{1}-$regularized procedure can be utilized to tighten the
$l_{2}-$error bound. 

The estimation error of the parameters of interests in the main equation
can be bounded by the maximum of a term involving the first-stage
estimation error and a term involving the second-stage estimation
error, which partially confirms%
\footnote{To verify whether the rate $\sqrt{\frac{\log p}{n}}$ is achievable
for the triangular simultaneous linear equations models, a minimax
lower bound result needs to be established in future work. %
} the speculation in Gautier and Tsybakov (2011) (Section 7.2) that
the two-stage estimation procedure can achieve the estimation error
of an order $\sqrt{\frac{\log p}{n}}$. My results show that $\sqrt{\frac{\log p}{n}}$\textbf{
}is achieved either when the second-stage estimation error dominates
the first-stage estimation error, or when $p$ is large relative to
$d$.\textbf{ }In the case where the second-stage estimation error
dominates the first-stage estimation error, the statistical error
of the parameters of interests in the main equation matches (up to
a factor of $|\beta^{*}|_{1}$) the order of the upper bound for the
Lasso estimate in the context of the high-dimensional linear regression
model without endogeneity, i.e., $\sqrt{\frac{k_{2}\log p}{n}}$.
An example of the second case where $p$ is large relative to $d$
is when the first-stage estimation concerns regressions in low-dimensional
settings and the result for this specific example is formally stated
in Corollary 3.4 below. \\
\textbf{}\\
\textbf{Corollary 3.4} (First-stage estimation in low-dimensional
settings): Suppose Assumptions 1.1, 3.2, and 3.3 hold. Assume the
number of regressors $p(=p_{n})$ in (1) can grow with and exceed
the sample size $n$; the number of non-zero components in $\beta^{*}$
is at most $k_{2}$, which is allowed to increase to infinity with
$n$ but slowly compared to $n$; and $d=\max_{j=1,...,p}d_{j}\ll n$ and does \textit{not} grow with $n$.
Suppose that the first-stage estimator $\hat{\pi}$ satisfies the
bound \textbf{$\max_{j=1,...,p}|\hat{\pi}_{j}-\pi_{j}^{*}|_{2}\precsim\sqrt{\frac{\log p}{n}}$}
with probability at least $1-O(\frac{1}{\max(p,\, n)})$. Then, if
\begin{eqnarray*}
\frac{k_{2}\log p}{n} & = & O(1),\\
\frac{\log p}{n} & = & o(1),
\end{eqnarray*}
and the tuning parameter $\lambda_{n}$ satisfies 
\[
\lambda_{n}\asymp k_{2}\sqrt{\frac{\log p}{n}},
\]
we have 
\[
|\hat{\beta}_{H2SLS}-\beta^{*}|_{2}\precsim\max\{\varphi_{1}\sqrt{k_{2}}\sqrt{\frac{\log p}{n}},\,\varphi_{2}\sqrt{\frac{k_{2}\log p}{n}}\},
\]
with probability at least $1-O(\frac{1}{\max(p,\, n)})$, where $\varphi_{1}$
and $\varphi_{2}$ are defined in Theorem 3.2. If we also have $k_{2}\sqrt{\frac{k_{2}\log p}{n}}=o(1)$,
then the two-stage estimator $\hat{\beta}_{H2SLS}$ is $l_{2}-$consistent
for \textbf{$\beta^{*}$}.\\

Note that Corollary 3.4 is a special case of Theorem 3.3 and hence the result is obvious from Theorem 3.3.

Under the condition that the first-stage estimators correctly select
the non-zero coefficients with probability close to 1, we can also
compare the high-dimensional two-stage estimator $\hat{\beta}_{H2SLS}$\textbf{
}with another type of multi-stage procedure. These multi-stage procedures
include three steps. In the first step, one carries out the same first-stage
estimation as before such as applying the Lasso or Dantzig selector.
Under some stringent conditions that guarantee the selection-consistency
of these first-stage estimators (such as the {}``irrepresentable
condition'' or the {}``mutual incoherence condition'' described
earlier), we can recover the supports of the true parameter vectors
with high probability. In the second step, we apply OLS with the regressors
in the estimated support set to obtain $\hat{\pi}_{j}^{OLS}$ for
$j=1,...,p$. In the third step, we apply a Lasso technique to the
main equation with these fitted regressors based on the second-stage
OLS estimates. This type of procedure is in the similar spirit as
the literature on sparsity in high-dimensional linear models without
endogeneity (see, e.g., Cand�s and Tao, 2007; Belloni and Chernozhukov,
2013). 

Under this three-stage procedure, Corollary 3.4 above tells us that
the statistical error of the parameters of interests in the main equation\textbf{
}is of the order \textbf{$O\left(|\beta^{*}|_{1}\sqrt{\frac{k_{2}\log p}{n}}\right)$},
which is at least as good as $\hat{\beta}_{H2SLS}$. Nevertheless,
this improved statistical error is at the expense of imposing stringent
conditions that ensure the first-stage estimators to achieve selection
consistency. These assumptions only hold in a rather narrow range
of problems, excluding many cases where the design matrices exhibit
strong (empirical) correlations. If these stringent conditions in
fact do not hold, then the three-stage procedure may not work. On
the other hand, even in the absence of the selection-consistency in
the first-stage estimation, $\hat{\beta}_{H2SLS}$ is still a valid
procedure and the bound as well as the consistency result in Theorem
3.2 still hold. Therefore, \textbf{$\hat{\beta}_{H2SLS}$} may be
more appealing in the sense that it works for a broader range of problems
in which the first-stage design matrices (formed by the instruments)
$Z_{j}\in\mathbb{R}^{n\times d_{j}}$ for $j=1,...,p$ exhibit a high
amount of dependency among the covariates.

For Theorems 3.2 and 3.3, the results are derived for the case where
each of the first-stage equations is estimated separately with a Lasso-type
procedure. Depending on the specific structures of the first-stage
equations, other methods that take into account the interrelationships
between these equations might yield a smaller first-stage estimation
error and consequently a potential improvement on the $l_{2}-$error
of $\hat{\beta}_{H2SLS}$. This paper does not pursue these more efficient
first-stage estimators but rather considers the extensions of Theorems
3.2 and 3.3 in the following manner. Notice that for Theorem 3.2 (or
Theorem 3.3), we give an explicit form of the first-stage estimation
error in Assumptions 3.5a (respectively, 3.5b) and as discussed earlier,
Lasso type of techniques yield these estimation errors. However, the
estimation error of the parameters of interests in the main equation
can be bounded by the maximum of a term involving the first-stage
estimation error in $l_{2}-$norm multiplied by a factor of $\sqrt{k_{1}k_{2}}$
(or $\sqrt{k_{2}}$ if the first-stage estimators correctly select
the non-zero coefficients with probability close to 1) and a term
involving the second-stage estimation error, which holds for general
first-stage estimation errors as long as $|\hat{\pi}_{j}-\pi_{j}^{*}|_{1}\asymp\sqrt{k_{1}}|\hat{\pi}_{j}-\pi_{j}^{*}|_{2}$
for%
\footnote{Negahban, et. al (2012) discusses the type of penalized estimators
that satisfy such a relationship between the $l_{1}-$error and the
$l_{2}-$error.%
} $j=1,...,p$. This claim is formally stated in Theorems 3.5 and 3.6
below. \\
\\
\textbf{Theorem 3.5}: Suppose Assumptions 1.1 and 3.1-3.3 hold. Also,
assume the first-stage estimator $\hat{\pi}$ satisfies the bound
$\max_{j=1,...,p}|\hat{\pi}_{j}-\pi_{j}^{*}|_{1}\leq\sqrt{k_{1}}M(d,\, p,\, k_{1},\, n)$
with probability $1-\alpha$. Then, if 
\begin{eqnarray*}
\max\left\{ k_{2}^{2}k_{1}M^{2}(d,\, p,\, k_{1},\, n),\:\frac{k_{2}\log p}{n}\right\}  & = & O(1),\\
\max\left\{ k_{1}M^{2}(d,\, p,\, k_{1},\, n),\:\frac{\log p}{n}\right\}  & = & o(1),
\end{eqnarray*}
and the tuning parameter $\lambda_{n}$ satisfies
\[
\lambda_{n}\asymp k_{2}\max\left\{ \sqrt{k_{1}}M(d,\, p,\, k_{1},\, n),\:\sqrt{\frac{\log p}{n}}\right\} ,
\]
we have 
\[
|\hat{\beta}_{H2SLS}-\beta^{*}|_{2}\precsim\max\{\varphi_{1}\sqrt{k_{1}k_{2}}M(d,\, p,\, k_{1},\, n),\,\varphi_{2}\sqrt{\frac{k_{2}\log p}{n}}\},
\]
where 
\begin{eqnarray*}
\varphi_{1} & = & \frac{\max_{j,j^{'}}|\textrm{cov}(x_{1j^{'}}^{*},\,\mathbf{z}_{1j})|_{\infty}|\beta^{*}|_{1}}{\lambda_{\min}(\Sigma_{X^{*}})},\\
\varphi_{2} & = & \max\left\{ \frac{\sigma_{X^{*}}\sigma_{\eta}|\beta^{*}|_{1}}{\lambda_{\min}(\Sigma_{X^{*}})},\,\frac{\sigma_{X^{*}}\sigma_{\epsilon}}{\lambda_{\min}(\Sigma_{X^{*}})}\right\} ,
\end{eqnarray*}
with probability at least $1-\alpha-c_{1}\exp(-c_{2}\log\max(p,\, n))$
for some universal positive constants $c_{1}$ and $c_{2}$. If we
also have $k_{2}\max\left\{ \sqrt{k_{1}k_{2}}M(d,\, p,\, k_{1},\, n),\,\sqrt{\frac{k_{2}\log p}{n}}\right\} =o(1)$,
then the two-stage estimator $\hat{\beta}_{H2SLS}$ is $l_{2}-$consistent
for \textbf{$\beta^{*}$}.\textbf{}\\
\textbf{}\\
\textbf{Theorem 3.6}: Suppose Assumptions 1.1, 3.1-3.4, and 3.6 hold.
Also, assume the first stage estimator $\hat{\pi}$ satisfies the
bound \textbf{$\max_{j=1,...,p}|\hat{\pi}_{j}-\pi_{j}^{*}|_{2}\leq M(d,\, p,\, k_{1},\, n)$}
with probability $1-\alpha$. Then, if\\
\\
$\min\left\{ \max\left\{ k_{2}^{2}k_{1}M^{2}(d,\, p,\, k_{1},\, n),\:\frac{k_{2}\log p}{n}\right\} ,\,\min_{r\in[0,\,1]}\max\left\{ k_{1}^{2-2r}M^{2}(d,\, p,\, k_{1},\, n),\,\frac{k_{1}^{r}k_{2}\log d}{n},\,\frac{k_{1}^{r}k_{2}\log p}{n}\right\} \right\} $
\[
=O(1),
\]
\[
\max\left\{ M^{2}(d,\, p,\, k_{1},\, n),\:\frac{\log p}{n}\right\} =o(1),
\]
and the tuning parameter $\lambda_{n}$ satisfies

\[
\lambda_{n}\asymp k_{2}\max\left\{ M(d,\, p,\, k_{1},\, n),\:\sqrt{\frac{\log p}{n}}\right\} ,
\]
we have 
\[
|\hat{\beta}_{H2SLS}-\beta^{*}|_{2}\precsim\max\{\varphi_{1}\sqrt{k_{2}}M(d,\, p,\, k_{1},\, n),\,\varphi_{2}\sqrt{\frac{k_{2}\log p}{n}}\},
\]
with probability at least $1-\alpha-c_{1}\exp(-c_{2}\log\max(p,\, n))$
for some universal positive constants $c_{1}$ and $c_{2}$, where
$\varphi_{1}$ and $\varphi_{2}$ are defined in Theorem 3.5.\textbf{
}If we also have $k_{2}\max\left\{ \sqrt{k_{2}}M(d,\, p,\, k_{1},\, n),\,\sqrt{\frac{k_{2}\log p}{n}}\right\} =o(1)$,
then the two-stage estimator $\hat{\beta}_{H2SLS}$ is $l_{2}-$consistent
for \textbf{$\beta^{*}$}.\\

The proofs for Theorems 3.5 and 3.6 are provided in Section 6.4.

Upon an additional condition that the first-stage estimators correctly
select the non-zero coefficients with probability close to 1, note
that the bound on the $l_{2}-$error of $\hat{\beta}_{H2SLS}$ in
Theorem 3.6 is improved upon that in Theorem 3.5 by a factor of $\sqrt{k_{1}}$\textbf{
}if the first term in the braces dominates the second one. The scaling
of the sample size required in Theorem 3.6 is also improved upon that
in Theorem 3.5.

\subsection{Variable-selection consistency}

In this subsection, I address the following question:\textbf{ }given
an optimal two-stage Lasso solution $\hat{\beta}_{H2SLS}$, when do
we have $\mathbb{P}[J(\hat{\beta}_{H2SLS})=J(\beta^{*})]\rightarrow1$?
That is, when can we conclude $\hat{\beta}_{H2SLS}$ correctly selects
the non-zero coefficients in the main equation with high probability?
This property is referred to as \textit{variable-selection consistency}.\textbf{
}For consistent variable selection with the standard Lasso in the
context of linear models without endogeneity, it is known that the
so-called {}``neighborhood stability condition'' (Meinshausen and
B�hlmann, 2006) for the design matrix, re-formulated in a nicer form
as the {}``irrepresentable condition'' by Zhao and Yu, 2006, is
sufficient and necessary. A further refined analysis is given in Wainwright
(2009), which presents under a certain {}``incoherence condition''
the smallest sample size needed to recover a sparse signal. In this
paper, I adopt the analysis by Wainwright (2009), Ravikumar, Wainwright,
and Lafferty (2010), and Wainwright (2014)\textbf{ }to analyze the
selection consistency of $\hat{\beta}_{H2SLS}$. In particular, I
need the following assumptions.\\
\\
\textbf{Assumption 3.7}:\textbf{ $\left\Vert \mathbb{E}\left[X_{1,J(\beta^{*})^{c}}^{*T}X_{1,J(\beta^{*})}^{*}\right]\left[\mathbb{E}(X_{1,J(\beta^{*})}^{*T}X_{1,J(\beta^{*})}^{*})\right]^{-1}\right\Vert _{\infty}\leq1-\phi$
}for some $\phi\in(0,\,1]$.\textbf{}\\
\textbf{}\\
\textbf{Assumption 3.8}: The smallest eigenvalue of the submatrix
$\mathbb{E}\left[X_{1,J(\beta^{*})}^{*T}X_{1,J(\beta^{*})}^{*}\right]$
satisfies the bound 
\[
\lambda_{\min}\left(\mathbb{E}\left[X_{1,J(\beta^{*})}^{*T}X_{1,J(\beta^{*})}^{*}\right]\right)\geq C_{\min}>0.
\]
\textbf{}\\
\textbf{Remarks}\\
Assumption 3.7, the so-called {}``mutual incoherence condition''
originally formalized by Wainwright (2009), captures the intuition
that the large number of irrelevant covariates cannot exert an overly
strong effect on the subset of relevant covariates. In the most desirable
case, the columns indexed by $j\in J(\beta^{*})^{c}$ would all be
orthogonal to the columns indexed by $j\in J(\beta^{*})$\textbf{
}and then we would have $\phi=1$.\textbf{ }In the high-dimensional
setting, this perfect orthogonality is not possible, but one can still
hope for a type of {}``near orthogonality'' to hold. 

Notice that in order for the left-hand-side of the inequality in Assumption
3.7 to always fall in $[0,\,1)$, one needs some type of normalization
on the matrix $X_{j}^{*}=(X_{1j}^{*},\,...\,,X_{nj}^{*})^{T}$ for
all $j=1,...,p$. One possibility is to impose a column normalization
as follows: 
\[
\max_{j=1,...,p}\frac{|X_{j}^{*}|_{2}}{\sqrt{n}}\leq\kappa_{c},\quad0<\kappa_{c}<\infty.
\]
Under Assumptions 1.1 and 3.3, we know that each column $X_{j}^{*}$,
$j=1,...,p$ is consisted of \textit{i.i.d.} sub-Gaussian variables.
Without loss of generality, we can assume $\mathbb{E}(X_{1j}^{*})=0$
for all $j=1,...,p$. Consequently, the normalization above follows
from a standard bound for the norms of zero-mean sub-Gaussian vectors
and a union bound 
\[
\mathbb{P}\left[\max_{j=1,...,p}\frac{|X_{j}^{*}|_{2}}{\sqrt{n}}\leq\kappa_{c}\right]\geq1-2\exp(-cn+\log p)\geq1-2\exp(-c^{'}n),
\]
where the last inequality follows from $n\gg\log p$. For example,
if $X^{*}$ has a Gaussian design, then we have 
\[
\max_{j=1,...,p}\frac{|X_{j}^{*}|_{2}}{\sqrt{n}}\leq\max_{j=1,...,p}\Sigma_{jj}\left(1+\sqrt{\frac{32\log p}{n}}\right),
\]
where $\max_{j=1,..,p}\Sigma_{jj}$ corresponds to the maximal variance
of any element of $X^{*}$(see Raskutti, et. al, 2011).

Assumption 3.8 is required to ensure that the model is identifiable
even if the support set $J(\beta^{*})$ were known \textit{a} \textit{priori}.
Assumption 3.8 is relatively mild compared to Assumption 3.7.\\
\textbf{}\\
\textbf{Theorem 3.7} (Selection consistency):\textbf{ }Suppose Assumptions
1.1, 3.1-3.3, 3.5a, 3.7, and 3.8 hold. If 
\begin{eqnarray*}
\frac{1}{n}\max\left\{ k_{1}k_{2}^{3/2}\log p,\, k_{2}^{3}\log p\right\}  & = & O(1),\\
\frac{1}{n}k_{1}^{2}k_{2}^{2}\log\max(d,\, p) & = & o(1),
\end{eqnarray*}
and the tuning parameter $\lambda_{n}$ satisfies 
\[
\lambda_{n}\asymp k_{1}k_{2}\sqrt{\frac{\log\max(d,\, p)}{n}},
\]
then, we have: (a) The Lasso has a unique optimal solution $\hat{\beta}_{H2SLS}$,
(b) the support $J(\hat{\beta}_{H2SLS})\subseteq J(\beta^{*})$, 
\[
\textrm{(c)}\quad|\hat{\beta}_{H2SLS,\, J(\beta^{*})}-\beta_{H2SLS,\, J(\beta^{*})}^{*}|_{\infty}\leq c\max\left\{ \varphi_{1}k_{1}\sqrt{\frac{k_{2}\log\max(d,\, p)}{n}},\:\varphi_{2}\sqrt{\frac{k_{2}\log p}{n}}\right\} :=B_{1},
\]
where 
\begin{eqnarray*}
\varphi_{1} & = & \frac{\sigma_{\eta}\max_{j,j^{'}}|\textrm{cov}(x_{1j^{'}}^{*},\,\mathbf{z}_{1j})|_{\infty}|\beta^{*}|_{1}}{\lambda_{\min}(\Sigma_{Z})C_{\min}},\\
\varphi_{2} & = & \max\left\{ \frac{\sigma_{X^{*}}\sigma_{\eta}|\beta^{*}|_{1}}{C_{\min}},\;\frac{\sigma_{X^{*}}\sigma_{\epsilon}}{C_{\min}}\right\} ,
\end{eqnarray*}
with probability at least $1-c_{1}\exp(-c_{2}\log\max(\min(p,\, d),\, n))$,
(d) if $\min_{j\in J(\beta^{*})}|\beta_{j}^{*}|>B_{1}$, then \textbf{$J(\hat{\beta}_{H2SLS})\supseteq J(\beta^{*})$
}and hence $\hat{\beta}_{H2SLS}$ is variable-selection consistent,
i.e., \textbf{$J(\hat{\beta}_{H2SLS})=J(\beta^{*})$}.\textbf{}\\
\textbf{}\\
\textbf{Theorem 3.8} (Selection consistency): Suppose Assumptions
1.1, 3.1-3.4, 3.5b, 3.6-3.8 hold. If 
\begin{eqnarray*}
\frac{1}{n}\max\left\{ k_{1}^{\nicefrac{1}{2}}k_{2}^{\nicefrac{3}{2}}\log p\; k_{2}^{3}\log p\right\}  & = & O(1),\\
\frac{1}{n}\min\left\{ k_{1}^{2}k_{2}^{2}\log\max(d,\, p),\,\min_{r\in[0,\,1]}\max\left\{ k_{1}^{3-2r}\log d,\, k_{1}^{3-2r}\log p,\, k_{1}^{r}k_{2}\log d,\, k_{1}^{r}k_{2}\log p\right\} \right\}  & = & o(1),\\
\frac{1}{n}k_{1}k_{2}^{2}\log\max(d,\, p) & = & o(1),
\end{eqnarray*}
and the tuning parameter $\lambda_{n}$ satisfies 
\[
\lambda_{n}\asymp k_{2}\sqrt{\frac{k_{1}\log\max(d,\, p)}{n}},
\]
then, we have: (a) The Lasso has a unique optimal solution $\hat{\beta}_{H2SLS}$,
(b) the support $J(\hat{\beta}_{H2SLS})\subseteq J(\beta^{*})$, and
\[
\textrm{(c)}\quad|\hat{\beta}_{H2SLS,\, J(\beta^{*})}-\beta_{H2SLS,\, J(\beta^{*})}^{*}|_{\infty}\leq c^{'}\max\left\{ \varphi_{1}\sqrt{\frac{k_{1}k_{2}\log\max(d,\, p)}{n}},\:\varphi_{2}\sqrt{\frac{k_{2}\log p}{n}}\right\} :=B_{2}
\]
with probability at least $1-c_{1}\exp(-c_{2}\log\max(\min(p,\, d),\, n))$,
where $\varphi_{1}$ and $\varphi_{2}$ are defined in Theorem 3.7,
(d) if $\min_{j\in J(\beta^{*})}|\beta_{j}^{*}|>B_{2}$, then \textbf{$J(\hat{\beta}_{H2SLS})\supseteq J(\beta^{*})$
}and hence $\hat{\beta}_{H2SLS}$ is variable-selection consistent,
i.e., \textbf{$J(\hat{\beta}_{H2SLS})=J(\beta^{*})$}. \\

The proofs for Theorems 3.7 and 3.8 are provided in Section 6.6.

The proof for Theorems 3.7 and 3.8 hinges on an intermediate result
that shows the {}``mutual incoherence'' assumption on $\mathbb{E}[X_{1}^{*T}X_{1}^{*}]$
(the population version of $\frac{1}{n}X^{*T}X^{*}$) guarantees that,
with high probability, analogous conditions hold for the estimated
quantity $\frac{1}{n}\hat{X}^{T}\hat{X}$, formed by the fitted regressors
from the first-stage regression. This result is established in Lemma
6.5 in Section 6.5. 

The proofs for Theorems 3.7 and 3.8 are based on a construction called
Primal-Dual Witness (PDW) method developed by Wainwright (2009) (also
see Wainwright, 2014). This method constructs a pair $(\hat{\beta},\,\hat{\mu})$.
When this procedure succeeds, the constructed pair is primal-dual
optimal, and acts as a witness for the fact that the Lasso has a unique
optimal solution with the correct signed support. The procedure is
described in the following. 
\begin{enumerate}
\item Set $\hat{\beta}_{J(\beta^{*})^{c}}=0$.
\item Obtain $(\hat{\beta}_{J(\beta^{*})},\,\hat{\mu}_{J(\beta^{*})})$
by solving the oracle subproblem 
\[
\hat{\beta}_{J(\beta^{*})}\in\arg\min_{\beta_{J(\beta^{*})}\in\mathbb{R}^{k_{2}}}\{\frac{1}{2n}|y-\hat{X}_{J(\beta^{*})}\beta_{J(\beta^{*})}|_{2}^{2}+\lambda_{n}|\beta_{J(\beta^{*})}|_{1}\},
\]
and choose $\hat{\mu}_{J(\beta^{*})}\in\partial|\hat{\beta}_{J(\beta^{*})}|{}_{1}$,
where $\partial|\hat{\beta}_{J(\beta^{*})}|_{1}$ denotes the set
of subgradients at $\hat{\beta}_{J(\beta^{*})}$ for the function
$|\cdot|_{1}:\,\mathbb{R}^{k_{2}}\rightarrow\mathbb{R}$.
\item Solve for $\hat{\mu}_{J(\beta^{*})^{c}}$ via the zero-subgradient
equation 
\[
\frac{1}{n}\hat{X}^{T}(y-\hat{X}\hat{\beta})+\lambda_{n}\hat{\mu}=0,
\]
and check whether or not the \textit{strict dual feasibility} condition
$|\hat{\mu}_{J(\beta^{*})^{c}}|_{\infty}<1$ holds.
\end{enumerate}
Theorems 3.7 and 3.8 include four parts. Part (a) guarantees the uniqueness
of the optimal solution of the two-stage Lasso procedure, $\hat{\beta}_{H2SLS}$
(from the proofs for Theorems 3.7 and 3.8, we have that $\hat{\beta}_{H2SLS}=(\hat{\beta}_{J(\beta^{*})},\,\mathbf{0})$
where $\hat{\beta}_{J(\beta^{*})}$ is the solution obtained in step
2 of the PDW construction above). Based on this uniqueness claim,
one can then talk unambiguously about the support of the two-stage
Lasso estimate. Part (b) guarantees that the Lasso does not falsely
include elements that are not in the support of $\beta^{*}$. 

Part (c) ensures that $\hat{\beta}_{H2SLS,\, J(\beta^{*})}$ is uniformly
close to $\beta_{J(\beta^{*})}^{*}$ in the $l_{\infty}-$norm%
\footnote{The factor $\sqrt{k_{2}}$ in the $l_{\infty}-$bound in Theorems
3.7 and 3.8 seems to be extra and removing it may require a more involved
analysis in future work.%
}. Notice that the $l_{\infty}-$bound in Part (c) of Theorem 3.8 is
improved by a factor of $\sqrt{k_{1}}$ upon that in Part (c) of Theorem
3.7 if the first term in the braces dominates the second one. Also,
the scaling of the sample size required in Theorem 3.8 is improved
upon that in Theorem 3.7. Similar observations were made earlier when
we compared the bound in Theorem 3.2 with the bound in Theorem 3.3
(or, the bound in Theorem 3.5 with the bound in Theorem 3.6). Again,
these observations are attributed to that the additional assumption
of the first-stage estimators correctly selecting the non-zero coefficients
(Assumption 3.6) is imposed in Theorem 3.8 but not in Theorem 3.7.
Recall earlier comparison between Theorem 3.2 and Theorem 3.3 in the
scaling of the required sample size. A similar comparison can be made
between Theorem 3.7 and Theorem 3.8. Under the assumption that the
first-stage estimators correctly select the non-zero coefficients
with high probability (Assumption 3.6), the scaling of the sample
size required\textbf{ }in Theorem 3.8 is guaranteed to be no greater
(and in some cases strictly smaller) than that in Theorem 3.7. For
instance, if $p\leq d$ , then by letting $r=1$, 
\[
\min\left\{ k_{1}^{2}k_{2}^{2}\log\max(d,\, p),\,\max\left\{ k_{1}\log d,\, k_{1}\log p,\, k_{1}k_{2}\log d,\, k_{1}k_{2}\log p\right\} \right\} =k_{1}k_{2}\log d.
\]
In this example, Theorem 3.7 suggests that the choice of sample size
needs to satisfy $\frac{k_{1}^{2}k_{2}^{2}\log d}{n}=o(1)$ and $\frac{\max\left\{ k_{1}k_{2}^{3/2}\log p,\, k_{2}^{3}\log p\right\} }{n}=O(1)$
while Theorem 3.8 suggests that the choice of sample size only needs
to satisfy $\frac{k_{1}k_{2}^{2}\log d}{n}=o(1)$ and $\frac{\max\left\{ k_{1}^{\nicefrac{1}{2}}k_{2}^{\nicefrac{3}{2}}\log p,\; k_{2}^{3}\log p\right\} }{n}=O(1)$.
However, as discussed previously, it is possible that the error bound
and scaling of the sample size required in Theorem 3.7 is suboptimal.
Section 5 provides a heuristic argument that may potentially improve
the bound in Theorem 3.7 when the first-stage estimates fail to satisfy
the exact sparsity constraint specified by the $l_{0}$- ball. 

The last claim is a consequence of this uniform norm bound: as long
as the minimum value of $|\beta_{j}^{*}|$ over $j\in J(\beta^{*})$
is not too small, then the two-stage Lasso does not falsely exclude
elements that are in the support of $\beta^{*}$ with high probability.
The minimum value requirement of $|\beta_{j}^{*}|$ over $j\in J(\beta^{*})$
is comparable to the so-called {}``beta-min'' condition in B�hlmann
and van de Geer (2011). Combining the claims from (b) and (d), the
two-stage Lasso is variable-selection consistent with high probability.

\section{Simulations }

In this section, simulations are conducted to gain insight on the
finite sample performance of the regularized two-stage estimators.
I consider the triangular simultaneous equations model (1) and (2)
from Section 1 where $d_{j}=d$ for all $j=1,...,p$, $(y_{i},\,\mathbf{x}_{i}^{T},\,\mathbf{z}_{i}^{T},\,\epsilon_{i},\,\boldsymbol{\eta}_{i})$
are $i.i.d.$, and $(\epsilon_{i},\,\boldsymbol{\eta}_{i})$ have
the following joint normal distribution 
\[
(\epsilon_{i},\,\boldsymbol{\eta}_{i})\thicksim\mathcal{N}\left(\left(\begin{array}{c}
0\\
0\\
\vdots\\
0
\end{array}\right),\;\left(\begin{array}{ccccc}
\sigma_{\epsilon}^{2} & \rho\sigma_{\epsilon}\sigma_{\eta} & \cdots & \cdots & \rho\sigma_{\epsilon}\sigma_{\eta}\\
\rho\sigma_{\epsilon}\sigma_{\eta} & \sigma_{\eta}^{2} & 0 & \cdots & 0\\
\vdots & 0 & \sigma_{\eta}^{2} & \cdots & \vdots\\
\vdots & \vdots & \vdots & \ddots & 0\\
\rho\sigma_{\epsilon}\sigma_{\eta} & 0 & \cdots & 0 & \sigma_{\eta}^{2}
\end{array}\right)\right).
\]
The matrix $\mathbf{z}_{i}^{T}$ is a $p\times d$ matrix of normal
random variables with identical variances $\sigma_{z}$, and $\mathbf{z}_{ij}^{T}$
is independent of $(\epsilon_{i},\,\eta_{i1},...,\eta_{ip})$ for
all $j=1,...,p$. With this setup, I simulate 1000 sets of $(y_{i},\,\mathbf{x}_{i}^{T},\,\mathbf{z}_{i}^{T},\,\epsilon_{i},\,\boldsymbol{\eta}_{i})_{i=1}^{n}$
where $n$ is the sample size (i.e., the number of data points) in
each set, and perform 14 Monte Carlo simulation experiments constructed
from various combinations of model parameters ($d$, $k_{1}$, $p$,
$k_{2}$, $\beta^{*}$, $\sigma_{\epsilon}$, and $\sigma_{\eta}$),
the design of $\mathbf{z}_{i}$, the random matrix formed by the instrumental
variables, as well as the types of first-stage and second-stage estimators
employed (Lasso vs. OLS). For each replication $t=1,...,1000$, I
compute the estimates $\hat{\beta}^{t}$ of the main-equation parameters
$\beta^{*}$, $l_{2}-$errors of these estimates, $|\hat{\beta}^{t}-\beta^{*}|_{2}$,
and selection percentages of $\hat{\beta}^{t}$ (computed by the number
of the elements in $\hat{\beta}^{t}$ sharing the same sign as their
corresponding elements in $\beta^{*}$, divided by the total number
of elements in $\beta^{*}$). Table 4.1 displays the designs of the
14 experiments. For Experiment 1 and Experiments 3-14, I set the number
of parameters in each first-stage equation $d=100$, the number of
parameters in the main equation $p=50$, the number of non-zero parameters
in each first-stage equation $k_{1}=4$, the number of non-zero parameters
in the main equation $k_{2}=5$. Also, choose $(\pi_{j1}^{*},...,\pi_{j4}^{*})=\mathbf{1}$,
$(\pi_{j5}^{*},...,\pi_{j100}^{*})=\mathbf{0}$ for all $j=1,...,50$;
and $(\beta_{1}^{*},...,\beta_{5}^{*})=\mathbf{1}$, $(\beta_{6}^{*},...,\beta_{50}^{*})=\mathbf{0}$.
For convenience, in the following discussion, I will refer to those
non-zero parameters as {}``relevant'' parameters and those zero
parameters as {}``irrelevant'' parameters. Experiment 2 sets $d=4$,
$p=5$, $(\pi_{j1}^{*},...,\pi_{j4}^{*})=\mathbf{1}$, and $(\beta_{1}^{*},...,\beta_{5}^{*})=\mathbf{1}$.
The motivations of these experiments are explained in the following
discussion.\\
\\
\\

~~~~~~%
\begin{sideways}
\begin{tabular}{|c||c|c|c|c|c|c|c|c|c|c|c|c|c|c|}
\multicolumn{15}{c}{}\tabularnewline
\multicolumn{1}{c}{} & \multicolumn{1}{c}{} & \multicolumn{1}{c}{} & \multicolumn{1}{c}{} & \multicolumn{1}{c}{} & \multicolumn{1}{c}{} & \multicolumn{1}{c}{} & \multicolumn{1}{c}{} & \multicolumn{1}{c}{} & \multicolumn{1}{c}{} & \multicolumn{1}{c}{} & \multicolumn{1}{c}{} & \multicolumn{1}{c}{} & \multicolumn{1}{c}{} & \multicolumn{1}{c}{}\tabularnewline
\multicolumn{1}{c}{} & \multicolumn{1}{c}{} & \multicolumn{1}{c}{} & \multicolumn{1}{c}{} & \multicolumn{1}{c}{} & \multicolumn{1}{c}{} & \multicolumn{1}{c}{} & \multicolumn{1}{c}{} & \multicolumn{1}{c}{} & \multicolumn{1}{c}{} & \multicolumn{1}{c}{} & \multicolumn{1}{c}{} & \multicolumn{1}{c}{} & \multicolumn{1}{c}{} & \multicolumn{1}{c}{}\tabularnewline
\multicolumn{1}{c}{} & \multicolumn{1}{c}{} & \multicolumn{1}{c}{} & \multicolumn{1}{c}{} & \multicolumn{1}{c}{} & \multicolumn{1}{c}{} & \multicolumn{1}{c}{} & \multicolumn{1}{c}{} & \multicolumn{1}{c}{} & \multicolumn{1}{c}{} & \multicolumn{1}{c}{} & \multicolumn{1}{c}{} & \multicolumn{1}{c}{} & \multicolumn{1}{c}{} & \multicolumn{1}{c}{}\tabularnewline
\multicolumn{15}{c}{{\scriptsize Table 4.1: Designs of the Monte-Carlo simulation experiments,
1000 replications }}\tabularnewline
\hline 
{\tiny Experiment \#} & {\tiny 1} & {\tiny 2} & {\tiny 3} & {\tiny 4} & {\tiny 5} & {\tiny 6} & {\tiny 7} & {\tiny 8} & {\tiny 9} & {\tiny 10} & {\tiny 11} & {\tiny 12} & {\tiny 13} & {\tiny 14}\tabularnewline
\hline 
{\tiny $d$} & {\tiny 100} & {\tiny 4} & {\tiny NA} & {\tiny 100} & {\tiny 100} & {\tiny 100} & {\tiny 100} & {\tiny 100} & {\tiny 100} & {\tiny 100} & {\tiny 100} & {\tiny 100} & {\tiny 100} & {\tiny 100}\tabularnewline
\hline 
{\tiny $k_{1}$} & {\tiny 4} & {\tiny 4} & {\tiny NA} & {\tiny 4} & {\tiny 4} & {\tiny 4} & {\tiny 4} & {\tiny 4} & {\tiny 4} & {\tiny 4} & {\tiny 4} & {\tiny 4} & {\tiny 4} & {\tiny 4}\tabularnewline
\hline 
{\tiny $p$} & {\tiny 50} & {\tiny 5} & {\tiny 50} & {\tiny 50} & {\tiny 50} & {\tiny 50} & {\tiny 50} & {\tiny 50} & {\tiny 50} & {\tiny 50} & {\tiny 50} & {\tiny 50} & {\tiny 50} & {\tiny 50}\tabularnewline
\hline 
{\tiny $k_{2}$} & {\tiny 5} & {\tiny 5} & {\tiny 5} & {\tiny 5} & {\tiny 5} & {\tiny 5} & {\tiny 5} & {\tiny 5} & {\tiny 5} & {\tiny 5} & {\tiny 5} & {\tiny 5} & {\tiny 5} & {\tiny 5}\tabularnewline
\hline 
{\tiny $(\beta_{1},...,\beta_{5})$} & \textbf{\tiny 1} & \textbf{\tiny 1} & \textbf{\tiny 1} & \textbf{\tiny 1} & \textbf{\tiny 1} & \textbf{\tiny 1} & \textbf{\tiny 1} & \textbf{\tiny 1} & \textbf{\tiny 1} & \textbf{\tiny 1} & \textbf{\tiny 1} & \textbf{\tiny 1} & \textbf{\tiny 1} & \textbf{\tiny 0.01}\tabularnewline
\hline 
{\tiny $(\beta_{6},...,\beta_{50})$} & \textbf{\tiny 0} & {\tiny NA} & \textbf{\tiny 0} & \textbf{\tiny 0} & \textbf{\tiny 0} & \textbf{\tiny 0} & \textbf{\tiny 0} & \textbf{\tiny 0} & \textbf{\tiny 0} & \textbf{\tiny 0} & \textbf{\tiny 0} & \textbf{\tiny 0} & \textbf{\tiny 0} & \textbf{\tiny 0}\tabularnewline
\hline 
{\tiny $(\pi_{j1},...,\pi_{j4})$ for all $j$} & \textbf{\tiny 1} & \textbf{\tiny 1} & {\tiny NA} & \textbf{\tiny 1} & \textbf{\tiny 1} & \textbf{\tiny 1} & \textbf{\tiny 1} & \textbf{\tiny 1} & \textbf{\tiny 1} & \textbf{\tiny 1} & \textbf{\tiny 1} & \textbf{\tiny 1} & \textbf{\tiny 1} & \textbf{\tiny 1}\tabularnewline
\hline 
{\tiny $(\pi_{j5},...,\pi_{j100})$ for all $j$} & \textbf{\tiny 0} & {\tiny NA} & {\tiny NA} & \textbf{\tiny 0} & \textbf{\tiny 0} & \textbf{\tiny 0} & \textbf{\tiny 0} & \textbf{\tiny 0} & \textbf{\tiny 0} & \textbf{\tiny 0} & \textbf{\tiny 0} & \textbf{\tiny 0} & \textbf{\tiny 0} & \textbf{\tiny 0}\tabularnewline
\hline 
{\tiny $\sigma_{\epsilon}$} & {\tiny 0.4} & {\tiny 0.4} & {\tiny 0.4} & {\tiny 0.4} & {\tiny 0.4} & {\tiny 0.4} & {\tiny 1} & {\tiny 0.4} & {\tiny 0.4} & {\tiny 0.4} & {\tiny 1} & {\tiny 0.4} & {\tiny 0.4} & {\tiny 0.4}\tabularnewline
\hline 
{\tiny $\sigma_{\eta}$} & {\tiny 0.4} & {\tiny 0.4} & {\tiny NA} & {\tiny 0.4} & {\tiny 0.4} & {\tiny 0.4} & {\tiny 0.4} & {\tiny 1} & {\tiny 0.4} & {\tiny 0.4} & {\tiny 0.4} & {\tiny 1} & {\tiny 0.4} & {\tiny 0.4}\tabularnewline
\hline 
{\tiny $\sigma_{z}$} & {\tiny 1} & {\tiny 1} & {\tiny NA} & {\tiny 1} & {\tiny 1} & {\tiny 1} & {\tiny 1} & {\tiny 1} & {\tiny 0.4} & {\tiny 1} & {\tiny 1} & {\tiny 1} & {\tiny 0.4} & {\tiny 1}\tabularnewline
\hline 
{\tiny Row corr. in $\mathbf{z}_{i}^{T}$ for $i=1,...,n$} & {\tiny No} & {\tiny No} & {\tiny NA} & {\tiny No} & {\tiny No} & {\tiny No} & {\tiny No} & {\tiny No} & {\tiny No} & {\tiny Yes} & {\tiny Yes} & {\tiny Yes} & {\tiny Yes} & {\tiny No}\tabularnewline
\hline 
{\tiny 1st-stage estimation} & {\tiny Lasso} & {\tiny OLS} & {\tiny NA} & {\tiny OLS} & {\tiny Lasso} & {\tiny OLS} & {\tiny Lasso} & {\tiny Lasso} & {\tiny Lasso} & {\tiny Lasso} & {\tiny Lasso} & {\tiny Lasso} & {\tiny Lasso} & {\tiny Lasso}\tabularnewline
\hline 
{\tiny 2nd-stage estimation} & {\tiny Lasso} & {\tiny OLS} & {\tiny Lasso} & {\tiny Lasso} & {\tiny OLS} & {\tiny OLS} & {\tiny Lasso} & {\tiny Lasso} & {\tiny Lasso} & {\tiny Lasso} & {\tiny Lasso} & {\tiny Lasso} & {\tiny Lasso} & {\tiny Lasso}\tabularnewline
\hline 
\end{tabular}
\end{sideways}\\
\newpage{}The baseline experiment (Experiment 1) applies the two-stage
Lasso procedure to the endogenous sparse linear model with a triangular
simultaneous equations structure (1) and (2). For each data point
$i=1,...,n$, the instruments $\mathbf{z}_{i}^{T}$ is a $p\times d$
matrix of independent standard normal random variables. As a benchmark
for Experiment 1, Experiment 2 concerns the classical 2SLS procedure
when both stage equations are in the low-dimensional setting and the
supports of the true parameters in both stages are known \textit{a
priori}. As another benchmark for Experiment 1, Experiment 3 applies
a one-step Lasso procedure (without instrumenting the endogenous regressors)
to the same main equation model (1) as in Experiment 1. 

Experiments 4-6 concern, in a relatively large sample size setting
with sparsity, the performance of alternative {}``partially'' regularized
or non-regularized estimators: first-stage-OLS-second-stage-Lasso
(Experiment 4), first-stage-Lasso-second-stage-OLS (Experiment 5),
and first-stage-OLS-second-stage-OLS (Experiment 6). Experiments 7-14
return to the two-stage Lasso procedure with changes applied to the
model parameters that generate the data. Experiment 7 (Experiment
8) increases the standard deviation of the {}``noise'' in the main
equation, $\sigma_{\epsilon}$ (respectively, the standard deviation
of the {}``noise'' in the first-stage equations, $\sigma_{\eta}$);
Experiment 9 reduces $\sigma_{z}$, the standard deviation of the
{}``signal'', i.e., the instrumental variables; Experiment 10 introduces
correlations between the rows of the design matrix $\mathbf{z}_{i}^{T}$.
Notice that each row of $\mathbf{z}_{i}^{T}\in\mathbb{R}^{p\times d}$
is associated with each of the endogenous regressors and the row-wise
correlation in $\mathbf{z}_{i}^{T}$ hence introduces correlations
between the {}``purged'' regressors $X_{j}^{*}$ and $X_{j^{'}}^{*}$
for all $j\neq j^{'}$. The level of the correlation is set to $0.5$,
i.e., $\textrm{corr}(z_{ijl},\, z_{ij^{'}l})=0.5$ for $j\neq j^{'}$
and $l=1,...,d$ (notice that we still have $\textrm{corr}(z_{ijl},\, z_{ijl^{'}})=0$
for $l\neq l^{'}$ and $j=1,...,p$; i.e., there is no column-wise
correlation in $\mathbf{z}_{i}^{T}$). Experiment 11 (Experiment 12)
increases the {}``noise'' level in the main equation (respectively,
the {}``noise'' level in the first-stage equations) and introduces
the correlations between the {}``purged'' regressors $X_{j}^{*}$
and $X_{j^{'}}^{*}$ for all $j\neq j^{'}$ simultaneously. Experiment
13 reduces the {}``signal'' level of the instrumental variables
and introduces the correlations between the {}``purged'' regressors
$X_{j}^{*}$ and $X_{j^{'}}^{*}$ for all $j\neq j^{'}$ simultaneously.
Experiment 14 reduces the magnitude of $(\beta_{1}^{*},\,...\,,\beta_{5}^{*})$
from $(1,\,...\,,1)$ to $(0.01,\,...\,,0.01)$. 

The tuning parameters $\lambda_{1n}$ in the first-stage Lasso estimation
(in Experiments 1, 5, 7-14) are chosen according to the standard Lasso
theory of high-dimensional estimation techniques (e.g., Bickel, 2009);
in particular, $\lambda_{1n}=0.4\sqrt{\frac{\log d}{n}}$. The tuning
parameters $\lambda_{2n}$ in the second-stage Lasso estimation (in
Experiments 1, 3, 4, 7-14) are chosen according to the scaling condition
in Theorem 3.3; in particular, $\lambda_{2n}=0.1\cdot k_{2}\max\left\{ \sqrt{\frac{k_{1}\log d}{n}},\:\sqrt{\frac{\log p}{n}}\right\} $
in Experiments 1, 3, 4, 7-13 and $\lambda_{2n}=0.001\cdot k_{2}\max\left\{ \sqrt{\frac{k_{1}\log d}{n}},\:\sqrt{\frac{\log p}{n}}\right\} $
in Experiment 14. The value of $\lambda_{2n}$ in Experiments 1, 3,
4, 7-13 exceeds the value of $\lambda_{2n}$ in Experiments 14 by
a factor of $0.01$. This adjustment reflects the fact that the non-zero
parameters $(\beta_{1},...,\beta_{5})=(1,\,...\,,1)$ in Experiments
1, 3, 4, 7-13 exceed the non-zero parameters $(\beta_{1},...,\beta_{5})=(0.01,\,...\,,0.01)$
in Experiment 14 by a factor of $0.01$.

Figure 4.1a plots (in ascending values) the 1000 estimates of $\beta_{5}^{*}$
when the sample size $n=47$. The estimates of other {}``relevant''
main-equation parameters behave similarly as the estimates of $\beta_{5}^{*}$.
Figure 4.1b plots (in ascending values) the 1000 estimates of $\beta_{6}^{*}$
when the sample size $n=47$. The estimates of other {}``irrelevant''
main-equation parameters behave similarly as the estimates of $\beta_{6}^{*}$.
The sample size 47 satisfies the scaling condition in Theorem 3.3.
With the choice of $d=100$, $k_{1}=4$, $p=50$, $k_{2}=5$ in Experiments
1 and 3, the sample size $n=47$ represents a high-dimensional setting
with sparsity. Figure 4.1c (Figure 4.1d) is similar to Figure 4.1a
(Figure 4.1b) except that the sample size $n=4700$. 

With the 1000 estimates of the main-equation parameters from Experiments
1-3, Table 4.2 shows the mean of the $l_{2}-$errors of these estimates
(computed as $\frac{1}{1000}\sum_{t=1}^{1000}|\hat{\beta}^{t}-\beta^{*}|_{2}$),
the mean of the selection percentages (computed in a similar fashion
as the mean of the $l_{2}-$errors of the estimates of $\beta^{*}$),
the mean of the squared $l_{2}-$errors (i.e., the \textit{sample}
\textit{mean squared error}, SMSE, computed as $\frac{1}{1000}\sum_{t=1}^{1000}|\hat{\beta}^{t}-\beta^{*}|_{2}^{2}$),
and the sample squared bias $\sum_{j=1}^{50}(\bar{\hat{\beta}}_{j}-\beta_{j}^{*})^{2}$
(where $\bar{\hat{\beta}}_{j}=\frac{1}{1000}\sum_{t=1}^{1000}\hat{\beta}_{j}^{t}$
for $j=1,...,50$). To provide a sense of how well the first-stage
estimates behave, Table 4.2 also displays the {}``averaged'' mean
of the $l_{2}-$errors of the first-stage estimates (computed as $\frac{1}{50}\sum_{j=1}^{50}\frac{1}{1000}\sum_{t=1}^{1000}|\hat{\pi}_{j}^{t}-\pi_{j}^{*}|_{2}$),
the {}``averaged'' mean of the selection percentages of the first-stage
estimates (computed in a similar fashion as the {}``averaged'' mean
of the $l_{2}-$errors of the first-stage estimates), the {}``averaged''
mean of the squared $l_{2}-$errors (i.e., the {}``averaged'' SMSE,
computed as $\frac{1}{50}\sum_{j=1}^{50}\frac{1}{1000}\sum_{t=1}^{1000}|\hat{\pi}_{j}^{t}-\pi_{j}^{*}|_{2}^{2}$),
and the {}``averaged'' sample squared bias $\frac{1}{50}\sum_{j=1}^{50}\sum_{l=1}^{100}(\bar{\hat{\pi}}_{jl}-\pi_{jl}^{*})^{2}$
(where $\bar{\hat{\pi}}_{jl}=\frac{1}{1000}\sum_{t=1}^{1000}\hat{\pi}_{jl}^{t}$
for $j=1,...,50$ and $l=1,...,100$). 

Compared to the two-stage Lasso procedure, in estimating the {}``relevant''
main-equation parameters with both sample sizes $n=47$ and $n=4700$,
Figures 4.1a and 4.1c show that the classical 2SLS procedure where
the supports of the true parameters in both stages are known \textit{a
priori }produces larger estimates while the one-step Lasso procedure
(without instrumenting the endogenous regressors) produces smaller
estimates. The two-stage Lasso outperforms the classical 2SLS above
the $60^{th}$ percentile of the estimates while underestimates the
{}``relevant'' main-equation parameters below the $60^{th}$ percentile
relative to the classical 2SLS procedure. The one-step Lasso procedure
(without instrumenting the endogenous regressors) produces the poorest
estimates of the {}``relevant'' main-equation parameters. The mean
$\bar{\hat{\beta}}_{5}$ of the 1000 estimates $\hat{\beta}_{5}$
from the two-stage Lasso is $0.931$ (respectively, $1.000$ from
the classical 2SLS and $0.826$ from the one-step Lasso) when $n=47$
and $0.997$ (respectively, $1.000$ from the classical 2SLS and $0.988$
from the one-step Lasso) when $n=4700$. The fact that the two-stage
Lasso yields smaller estimates of the {}``relevant'' main-equation
parameters relative to the classical 2SLS for both sample sizes is
most likely due to the shrinkage effect from the $l_{1}-$penalization
in the second-stage estimation of the two-stage Lasso procedure.

In estimating the {}``irrelevant'' main-equation parameters, the
estimates of $(\beta_{6}^{*},...,\beta_{50}^{*})$ from both the two-stage
Lasso and the one-step Lasso are exactly $0$ at the $5^{th}$ percentile,
the median, and the $95^{th}$percentile when $n=47$ and $n=4700$.
The mean statistics of the estimates of $(\beta_{6}^{*},...,\beta_{50}^{*})$
range from $-0.001$ ($-2.938\times10^{-4}$) to $0.001$ ($3.403\times10^{-4}$)
when $n=47$, and $-6.926\times10^{-5}$ ($0$) to $5.593\times10^{-5}$
($3.076\times10^{-6}$) when $n=4700$ for the two-stage Lasso (respectively,
the one-step Lasso). Table 4.2 shows that the selection percentages
of the main-equation estimates from the two-stage Lasso and the one-step
Lasso are high for the designs considered. Figures 4.1b and 4.1d show
that, in estimating the {}``irrelevant'' main-equation parameters,
the one-step Lasso performs slightly better relative to the two-stage
Lasso procedure below the $2^{nd}$ percentile and above the $98^{th}$
percentile. 

In terms of estimation errors and sample bias, from Table 4.2 we see
that the mean of the $l_{2}-$errors of the estimates $\hat{\beta}_{H2SLS}$
of $\beta^{*}$ (or the {}``averaged'' mean of the $l_{2}-$errors
of the first-stage estimates) from the two-stage Lasso are greater
than those of $\hat{\beta}_{2SLS}$ (respectively, of the first-stage
estimates) from the classical 2SLS procedure for both $n=47$ and
$n=4700$. As $n$ increases, the mean of the $l_{2}-$errors of $\hat{\beta}_{H2SLS}$
and the mean of the $l_{2}-$errors of $\hat{\beta}_{2SLS}$ become
very close to each other as in the case when $n=4700$. Also, the
sample bias of $\hat{\beta}_{H2SLS}$ (or, the {}``averaged'' sample
bias of the first-stage estimates) from the two-stage Lasso are greater
by a magnitude of $100\sim1000$ (respectively, $1000\sim10^{4}$)
than those of $\hat{\beta}_{2SLS}$ (respectively, of the first-stage
estimates) from the classical 2SLS procedure for both sample sizes. 

For more investigation on how the $l_{2}-$error and sample bias of
$\hat{\beta}_{H2SLS}$ compare to those of $\hat{\beta}_{2SLS}$,
I have also considered designs where $\sigma_{\epsilon}$ and/or $\sigma_{\eta}$
are increased or decreased while everything else in Experiments 1
and 2 remains the same. In these modified designs except for those
with very large values of $\sigma_{\epsilon}$ under $n=47$, the
mean of the $l_{2}-$errors of $\hat{\beta}_{H2SLS}$ are generally
greater than those of $\hat{\beta}_{2SLS}$. The sample bias of $\hat{\beta}_{H2SLS}$
are consistently greater by a magnitude of $10\sim10^{5}$ than those
of $\hat{\beta}_{2SLS}$ for both sample sizes. This suggests that
the shrinkage effect from the $l_{1}-$penalization in both the first
and second stage estimations of the two-stage Lasso procedure might
have made its bias term converge to zero at a slower rate relative
to the classical 2SLS for the designs considered here. Whether this
conjecture holds true for general designs is an interesting question
for further research. Compared to the two-stage Lasso and the classical
2SLS, the one-step Lasso procedure without instrumenting the endogenous
regressors yields the largest $l_{2}-$errors as well as sample bias
of the main-equation estimates for both $n=47$ and $n=4700$, which
is expected. Finally notice that the $l_{2}-$errors (and the sample
bias) shrink as the sample size increases. $|\hat{\beta}_{2SLS}-\beta^{*}|_{2}$
being proportional to $\frac{1}{\sqrt{n}}$ is a known fact in low-dimensional
settings. From Section 3.1, we also have that the upper bounds for
$|\hat{\beta}_{H2SLS}-\beta^{*}|_{2}$ are proportional to $\frac{1}{\sqrt{n}}$
up to factors involving $\log d$, $\log p$, $k_{1}$, and $k_{2}$.
\\
\\
\\
\begin{tabular}{c}
\includegraphics[width=16cm,height=8cm,keepaspectratio]{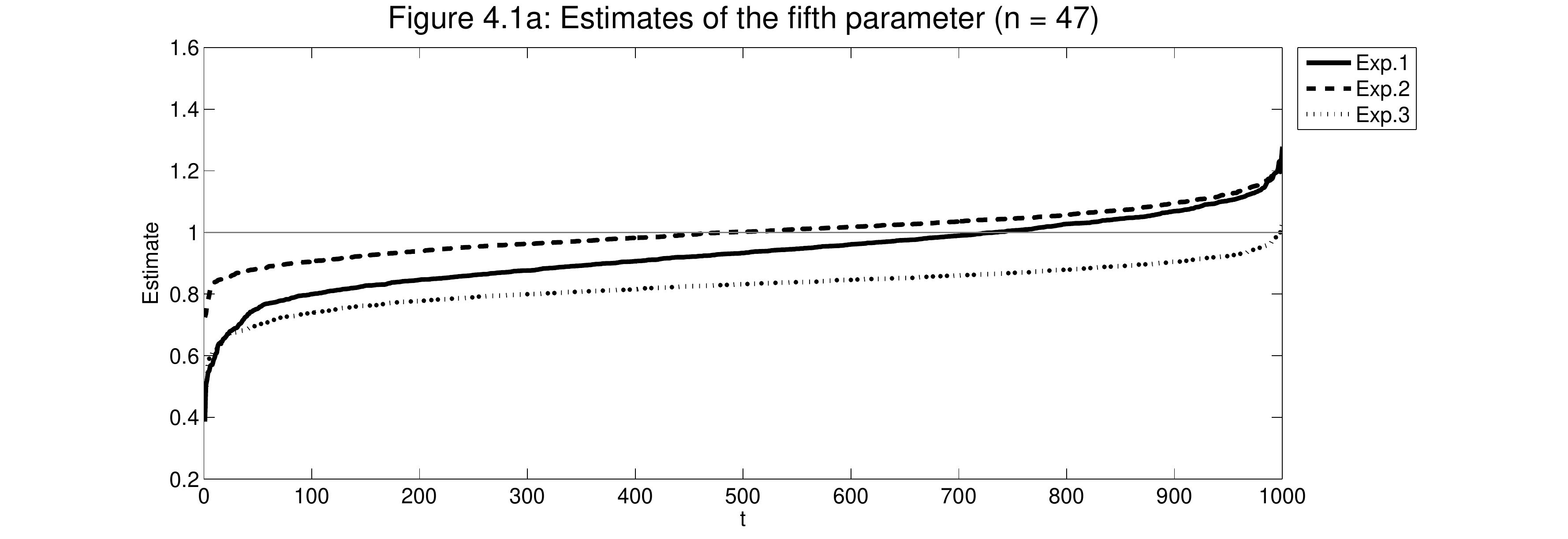}\tabularnewline
\end{tabular}\\
\begin{tabular}{c}
\includegraphics[width=16cm,height=8cm,keepaspectratio]{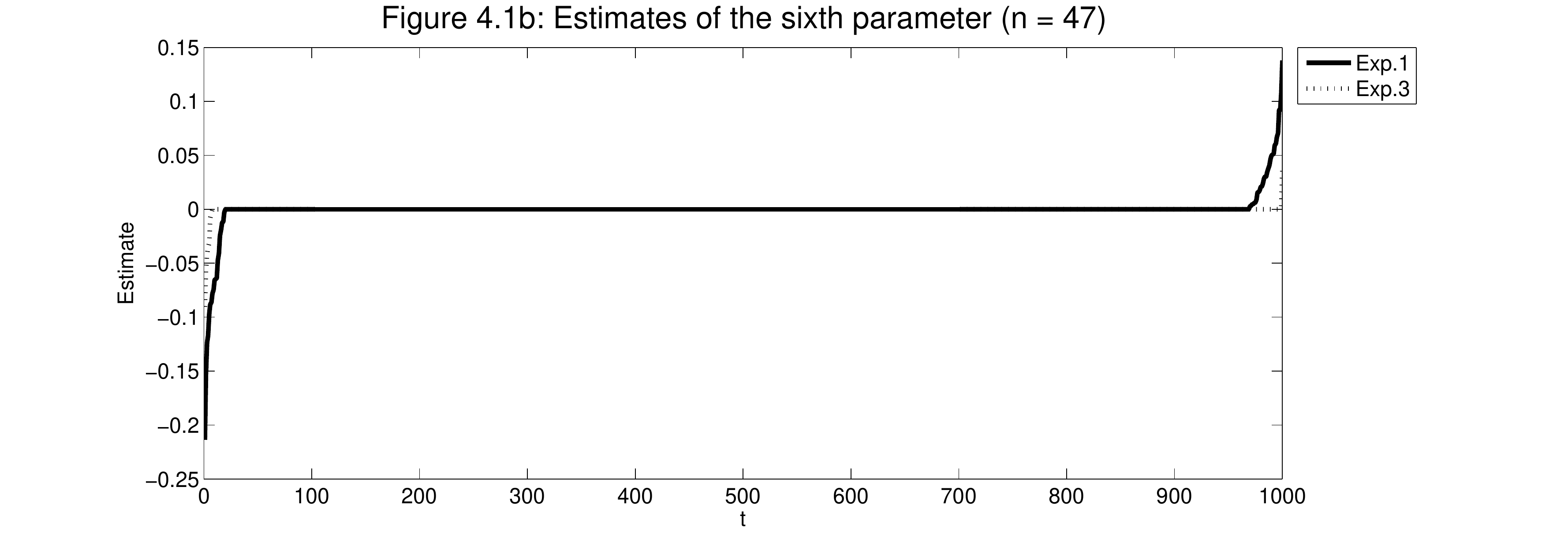}\tabularnewline
\end{tabular}\\
\begin{tabular}{c}
\includegraphics[width=16cm,height=8cm,keepaspectratio]{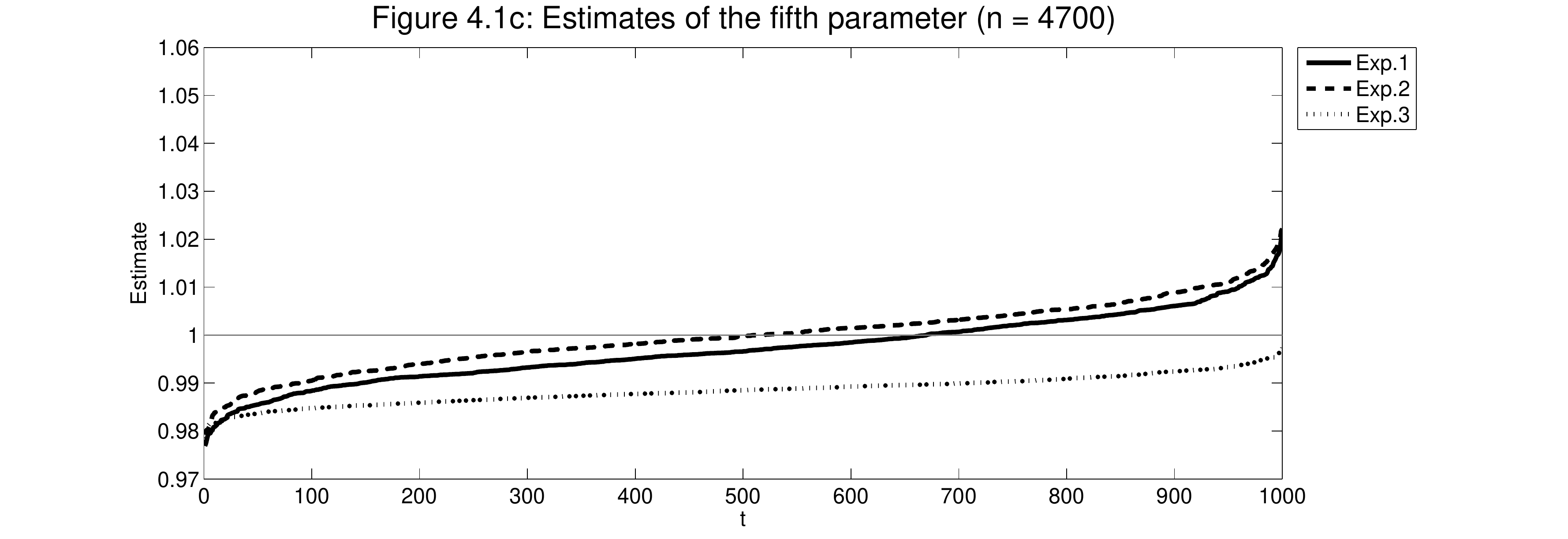}\tabularnewline
\end{tabular}\\
\begin{tabular}{c}
\includegraphics[width=16cm,height=8cm,keepaspectratio]{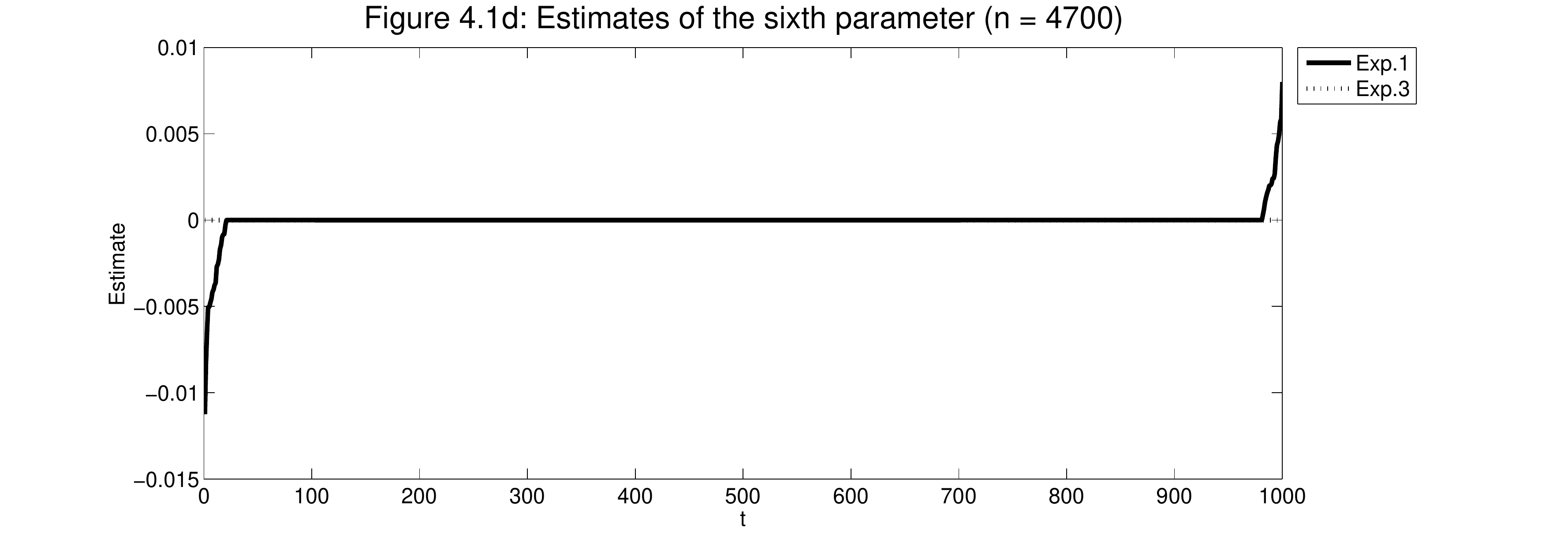}\tabularnewline
\end{tabular}\\
\\
\begin{tabular}{lll|c|c|c|c|c|c|c|}
 &  & \multicolumn{1}{l}{} & \multicolumn{7}{c}{{\scriptsize Table 4.2: $l_{2}-$errors, SMSE, bias, and selection
(Exp. 1-3)}}\tabularnewline
\cline{4-10} 
 &  &  & {\scriptsize Mean} & \multicolumn{3}{c|}{{\scriptsize $n=47$}} & \multicolumn{3}{c|}{{\scriptsize $n=4700$}}\tabularnewline
\cline{4-10} 
 &  &  & {\scriptsize Exp. \# } & {\scriptsize 1} & {\scriptsize 2} & {\scriptsize 3} & {\scriptsize 1} & {\scriptsize 2} & {\scriptsize 3}\tabularnewline
\cline{4-10} 
 &  &  & {\scriptsize $2^{nd}-$stage select \%} & {\scriptsize 97.3} & {\scriptsize NA} & {\scriptsize 99.5} & {\scriptsize 98.1} & {\scriptsize NA} & {\scriptsize 100}\tabularnewline
\cline{4-10} 
 &  &  & {\scriptsize $2^{nd}-$stage $l_{2}-$error} & {\scriptsize 0.288} & {\scriptsize 0.156} & {\scriptsize 0.412} & {\scriptsize 0.018} & {\scriptsize 0.015} & {\scriptsize 0.026}\tabularnewline
\cline{4-10} 
 &  &  & {\scriptsize $2^{nd}-$stage SMSE} & {\scriptsize 0.099} & {\scriptsize 0.028} & {\scriptsize 0.179} & {\scriptsize $3.38\times10^{-4}$} & {\scriptsize $2.43\times10^{-4}$} & {\scriptsize $7.07\times10^{-4}$}\tabularnewline
\cline{4-10} 
 &  &  & {\scriptsize $2^{nd}-$stage squared bias} & {\scriptsize 0.024} & {\scriptsize $1.49\times10^{-5}$} & {\scriptsize 0.154} & {\scriptsize $5.04\times10^{-5}$} & {\scriptsize $1.84\times10^{-7}$} & {\scriptsize $6.66\times10^{-4}$}\tabularnewline
\cline{4-10} 
 &  &  & {\scriptsize $1^{st}$-stage select \%} & {\scriptsize 0.977} & {\scriptsize NA} & {\scriptsize NA} & {\scriptsize 0.985} & {\scriptsize NA} & {\scriptsize NA}\tabularnewline
\cline{4-10} 
 &  &  & {\scriptsize $1^{st}-$stage $l_{2}-$error} & {\scriptsize 0.349} & {\scriptsize 0.115} & {\scriptsize NA} & {\scriptsize 0.028} & {\scriptsize 0.011} & {\scriptsize NA}\tabularnewline
\cline{4-10} 
 &  &  & {\scriptsize $1^{st}-$stage SMSE} & {\scriptsize 0.132} & {\scriptsize 0.003} & {\scriptsize NA} & {\scriptsize $7.92\times10^{-4}$} & {\scriptsize $2.78\times10^{-5}$} & {\scriptsize NA}\tabularnewline
\cline{4-10} 
 &  &  & {\scriptsize $1^{st}-$stage squared bias} & {\scriptsize 0.097} & {\scriptsize $1.07\times10^{-5}$} & {\scriptsize NA} & {\scriptsize $6.31\times10^{-4}$} & {\scriptsize $1.23\times10^{-7}$} & {\scriptsize NA}\tabularnewline
\cline{4-10} 
 &  & \multicolumn{1}{l}{} & \multicolumn{7}{c}{}\tabularnewline
\end{tabular}\\
\\
\\
In the following relatively large sample size setting (i.e., $n=4700$)
with sparsity, I compare the performance of the two-stage Lasso estimator
with the performances of the alternative {}``partially'' regularized
or non-regularized estimators as mentioned earlier. Figure 4.2a plots
(in ascending values) the 1000 estimates of $\beta_{5}^{*}$ when
the sample size $n=4700$. Figure 4.2b plots (in ascending values)
the 1000 estimates of $\beta_{6}^{*}$ when the sample size $n=4700$.
With the 1000 estimates of the {}``relevant'' ({}``irrelevant'')
main-equation parameters from Experiment 1 and Experiments 4-6, Figures
4.2c-4.2f (respectively, Figures 4.2g-4.2j) display the $5^{th}$
percentile, the median, the $95^{th}$ percentile, and the mean of
these estimates. The mean of the $l_{2}-$errors and the mean of the
selection percentages of the main-equation estimates together with
the {}``averaged'' mean of the $l_{2}-$errors and the {}``averaged''
mean of the selection percentages of the first-stage estimates from
these {}``partially'' regularized or non-regularized estimators
are displayed in Table 4.3. 

Figure 4.2a and Figures 4.2c-4.2f show that, compared to the two-stage
Lasso procedure, in estimating the {}``relevant'' main-equation
parameters when $n=4700$, the first-stage-Lasso-second-stage-OLS
estimator and the first-stage-OLS-second-stage-OLS estimator\textit{
}produce larger estimates while the first-stage-OLS-second-stage-Lasso
estimator produces smaller estimates. In estimating the {}``irrelevant''
main-equation parameters when $n=4700$, Figure 4.2b and Figures 4.2g-4.2j
show that the two-stage Lasso and the first-stage-OLS-second-stage-Lasso
estimator perform well while the first-stage-Lasso-second-stage-OLS
and the first-stage-OLS-second-stage-OLS do poorly (also see Table
4.3 for a comparison between the selection percentages of these estimators).
This suggests that employing regularization in the second-stage estimation
helps selecting the {}``relevant'' main-equation parameters. 

Turning to the comparison with {}``partially'' regularized or non-regularized
estimators in terms of $l_{2}-$errors, from Table 4.3 we see that
the two-stage Lasso estimator achieves the smallest $l_{2}-$error
of the main-equation estimates among all the estimators considered
here. The fact that the $l_{2}-$error (of the main-equation estimates)
of the two-stage Lasso estimator is smaller than the $l_{2}-$errors
of the first-stage-OLS-second-stage-Lasso estimator and the first-stage-OLS-second-stage-OLS
estimator could be attributed to the following. Based on the first-stage
estimation results from these experiments, the first-stage Lasso estimator
outperforms the first-stage OLS estimator in both estimation errors
and variable selections even in the relatively large sample size setting
with sparsity. Recall in Section 3, we have seen that, the estimation
error of the parameters of interests in the main equation can be bounded
by the maximum of a term involving the first-stage estimation error
and a term involving the second-stage estimation error. Given the
choices of $p$, $d$, $k_{1}$, and $k_{2}$ in Experiment 1 and
Experiments 4-6, these results agree with the theorems in Section
3.1. Additionally, compared to the first-stage-OLS-second-stage-OLS
estimator, the fact that the $l_{2}-$error (of the main-equation
estimates) of the two-stage Lasso estimator is smaller than the $l_{2}-$error
of the first-stage-OLS-second-stage-OLS estimator can also be explained
by the fact that the two-stage Lasso reduces the $l_{2}-$error of
the first-stage-OLS-second-stage-OLS estimates from $O\left(\sqrt{\frac{\max(p,\, d)}{n}}\right)$
to $O\left(\sqrt{\frac{\log\max(p,\, d)}{n}}\right)$, as we have
seen in Section 3. Similarly, the fact that the $l_{2}-$error (of
the main-equation estimates) of the two-stage Lasso estimator is smaller
than the $l_{2}-$error of the first-stage-Lasso-second-stage-OLS
estimator can be explained by the fact that the two-stage Lasso reduces
the $l_{2}-$error of the first-stage-Lasso-second-stage-OLS estimates
from $O\left(\sqrt{\frac{\max(p,\,\log d)}{n}}\right)$ to $O\left(\sqrt{\frac{\log\max(p,\, d)}{n}}\right)$.\\
\\
\begin{tabular}{c}
\includegraphics[width=16cm,height=10cm]{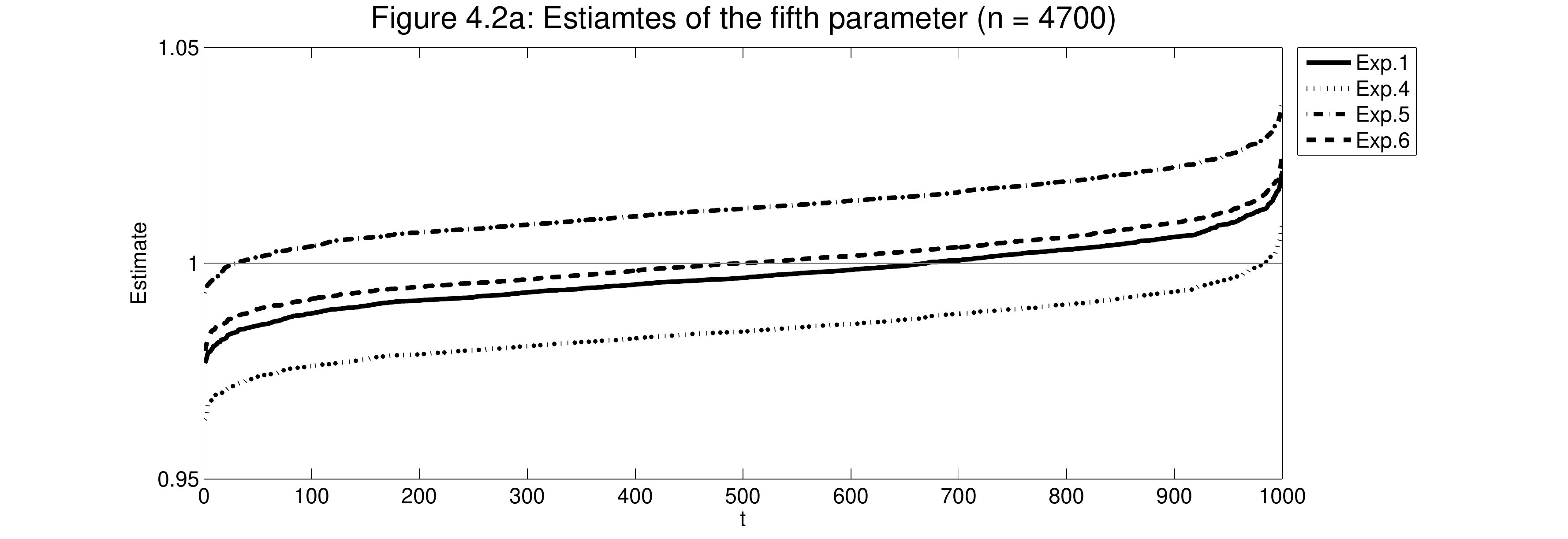}\tabularnewline
\includegraphics[width=16cm,height=10cm]{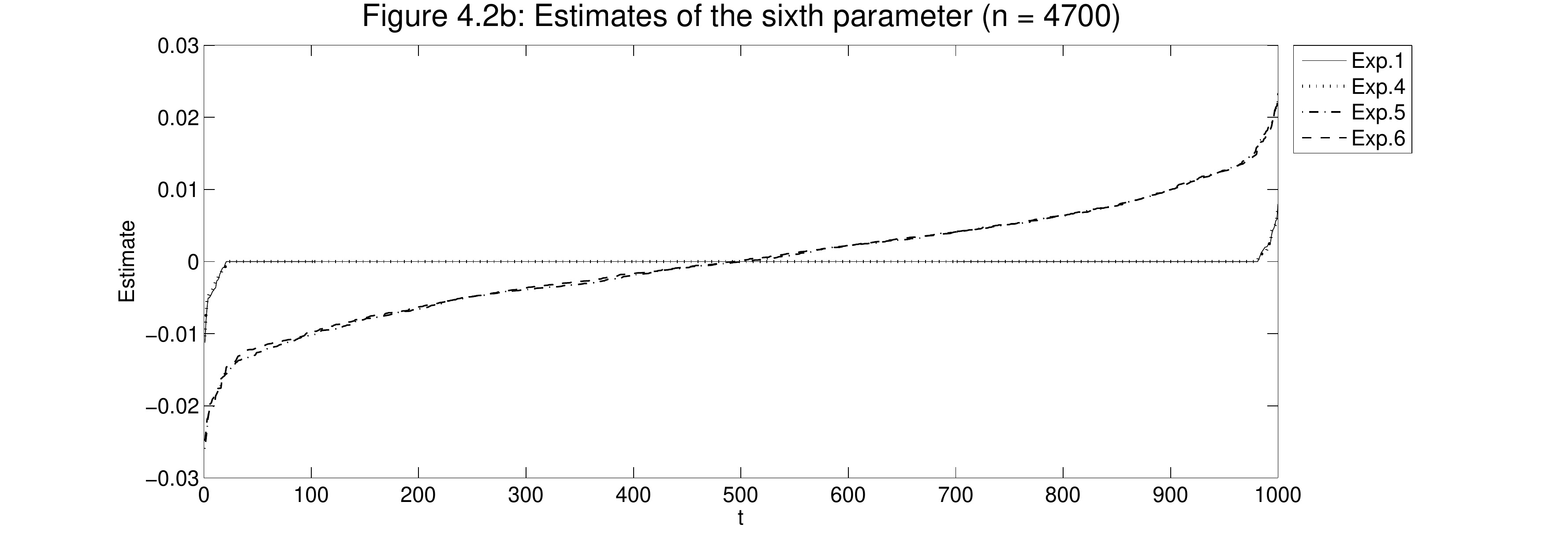}\tabularnewline
\end{tabular}\\
\\
\begin{tabular}{cc}
\includegraphics[width=8cm,height=8cm]{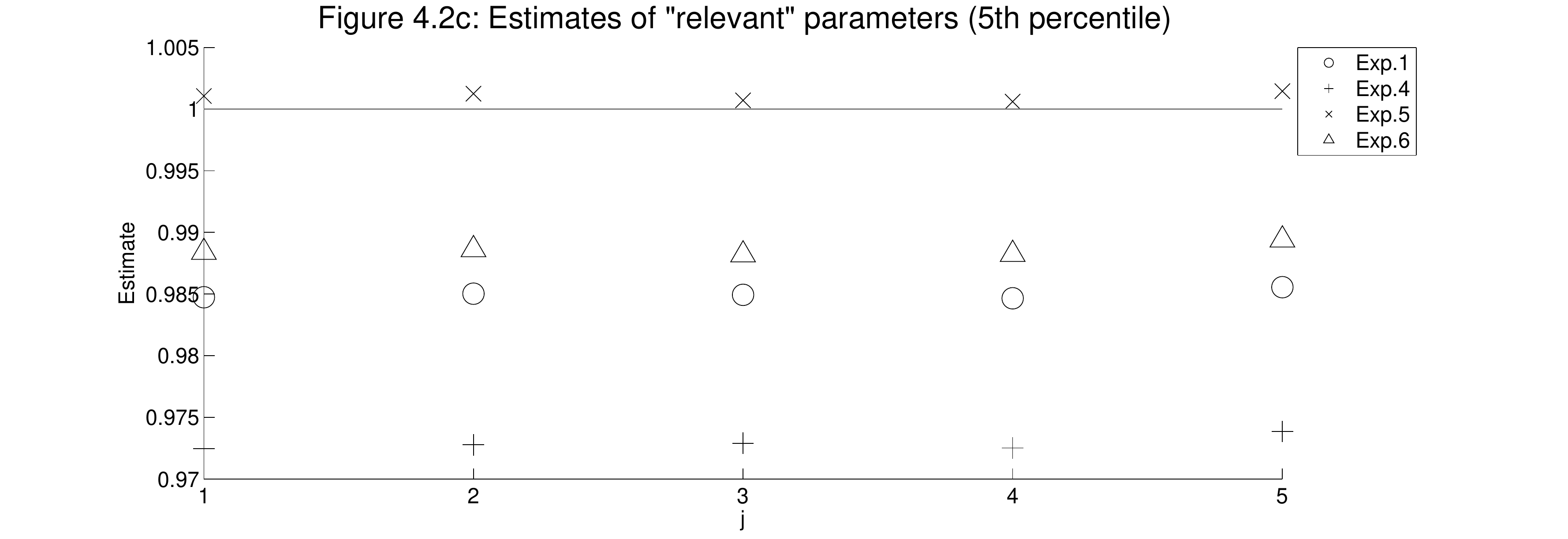} & \includegraphics[width=8cm,height=8cm]{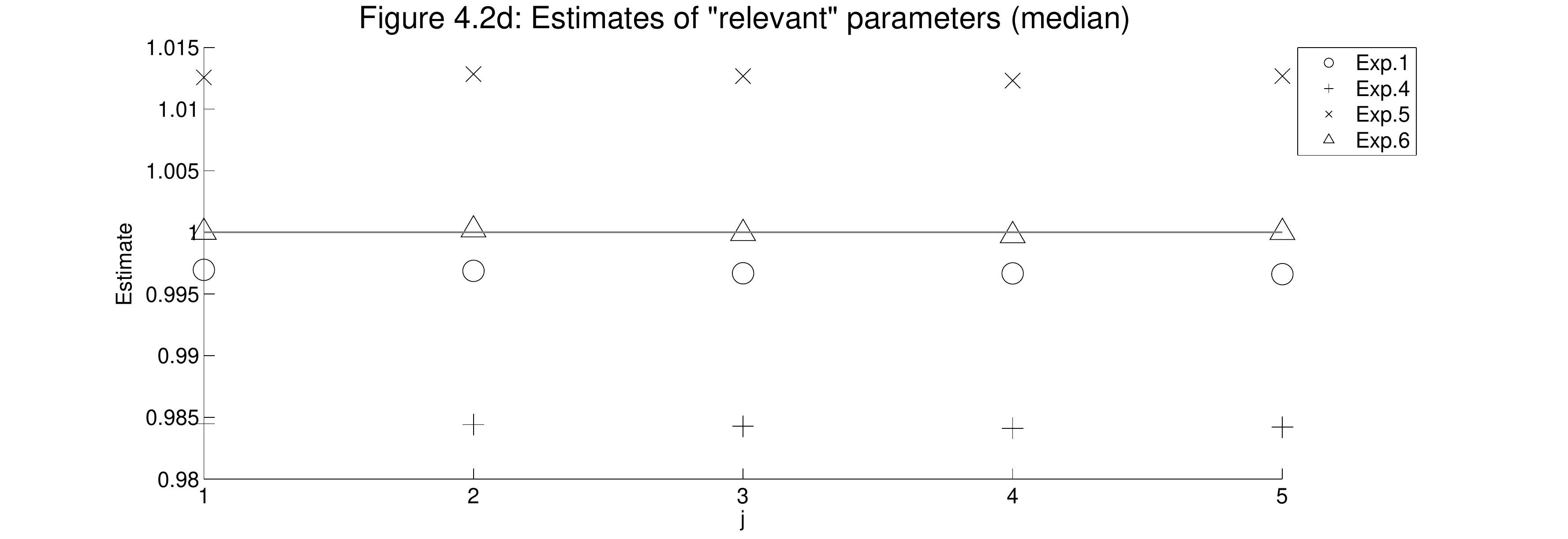}\tabularnewline
\end{tabular}\\
\begin{tabular}{cc}
\includegraphics[width=8cm,height=8cm]{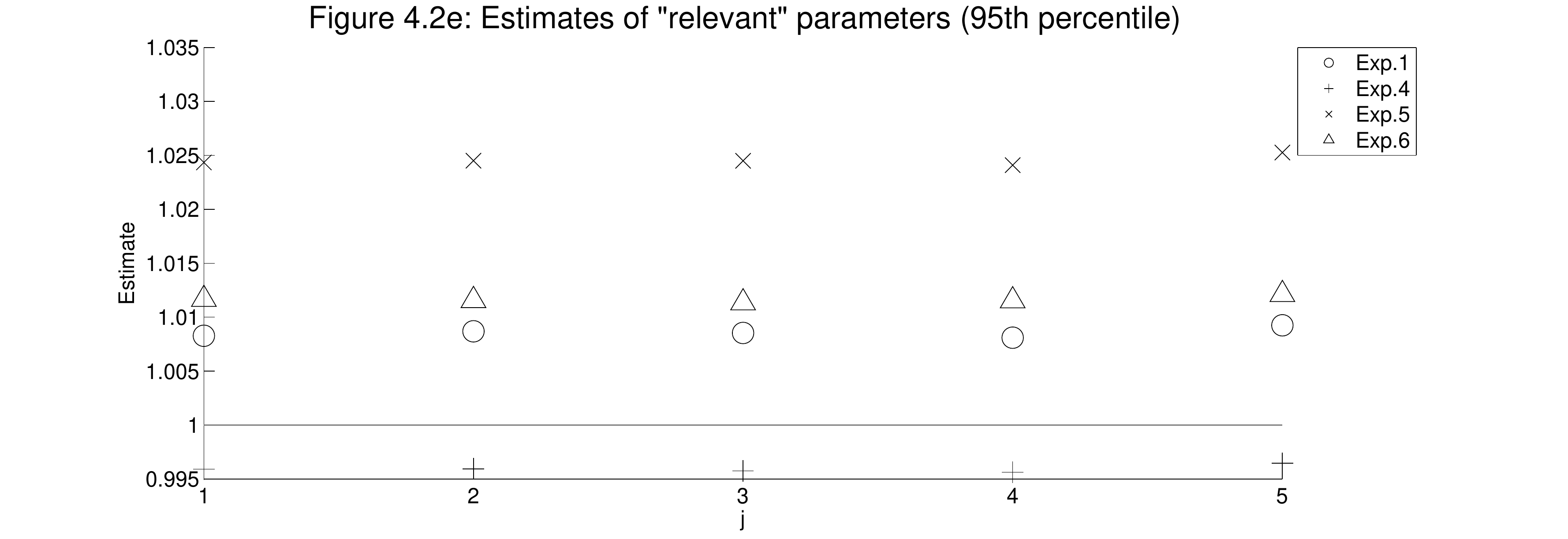} & \includegraphics[width=8cm,height=8cm]{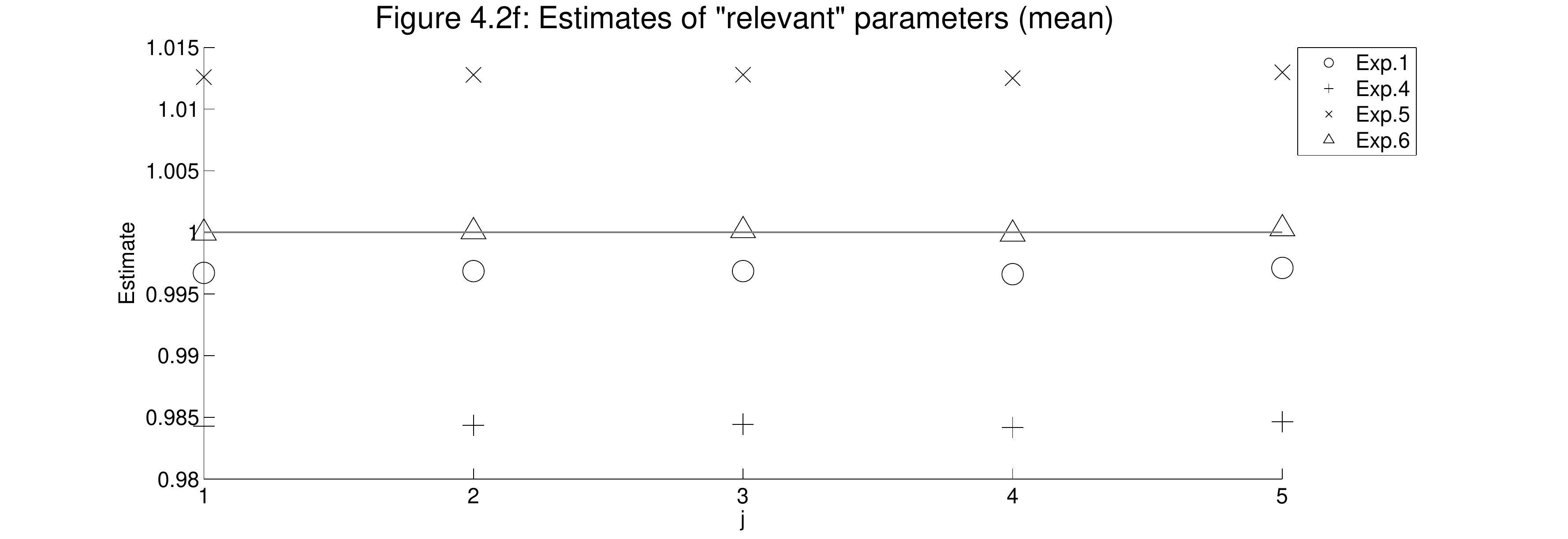}\tabularnewline
\end{tabular}\\
\begin{tabular}{cc}
\includegraphics[width=8cm,height=8cm]{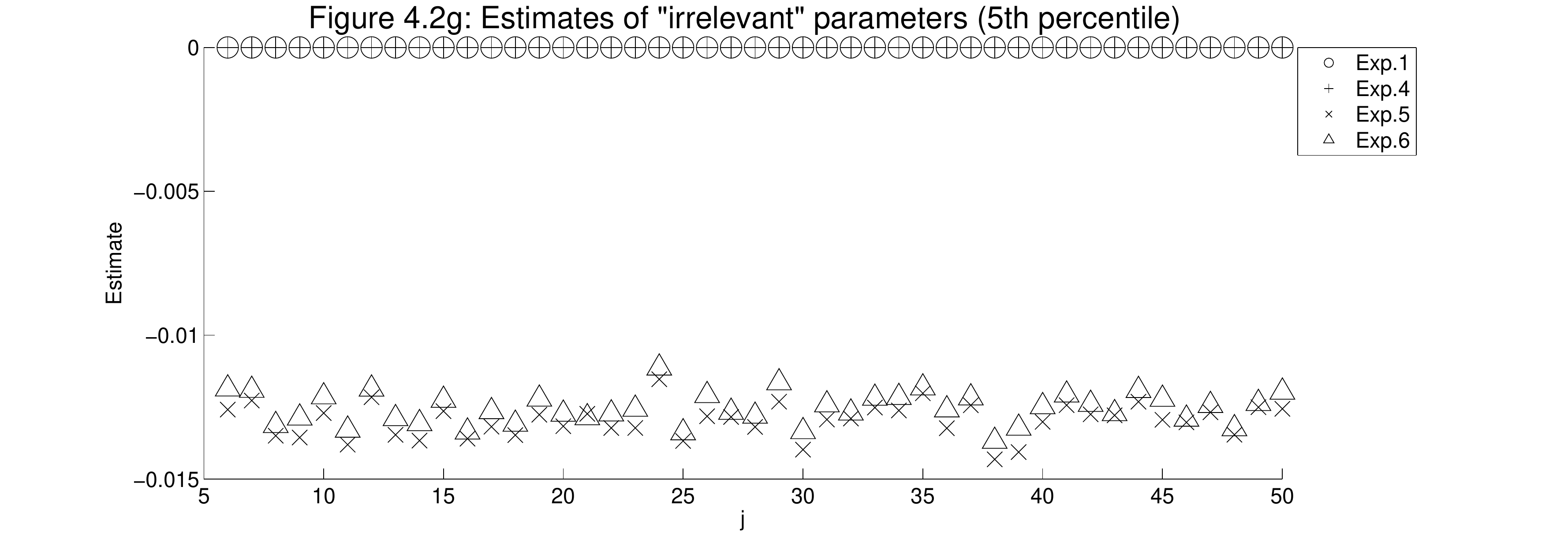} & \includegraphics[width=8cm,height=8cm]{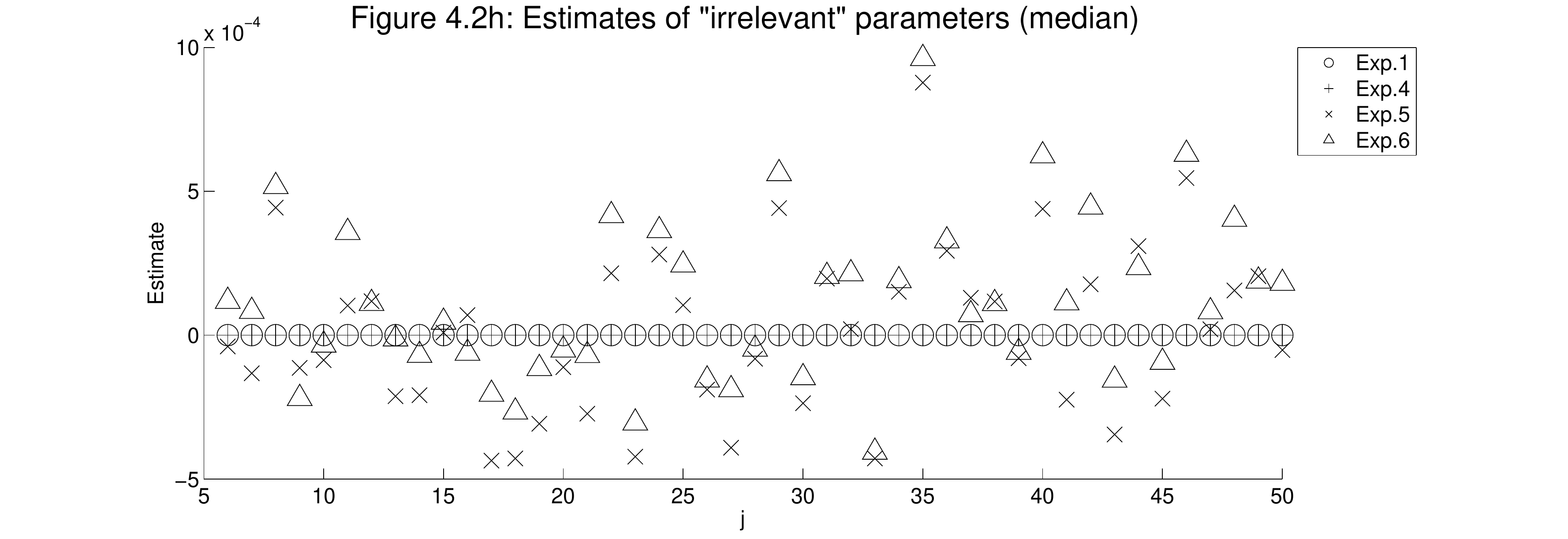}\tabularnewline
\end{tabular}\\
\begin{tabular}{cc}
\includegraphics[width=8cm,height=8cm]{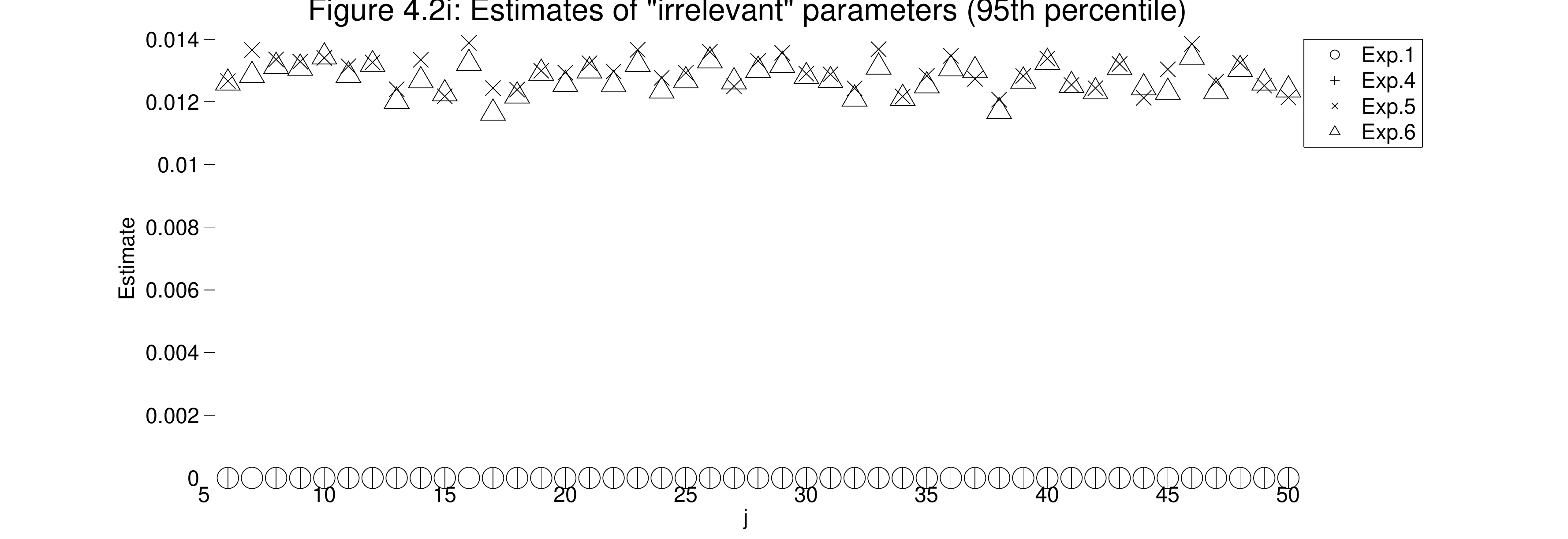} & \includegraphics[width=8cm,height=8cm]{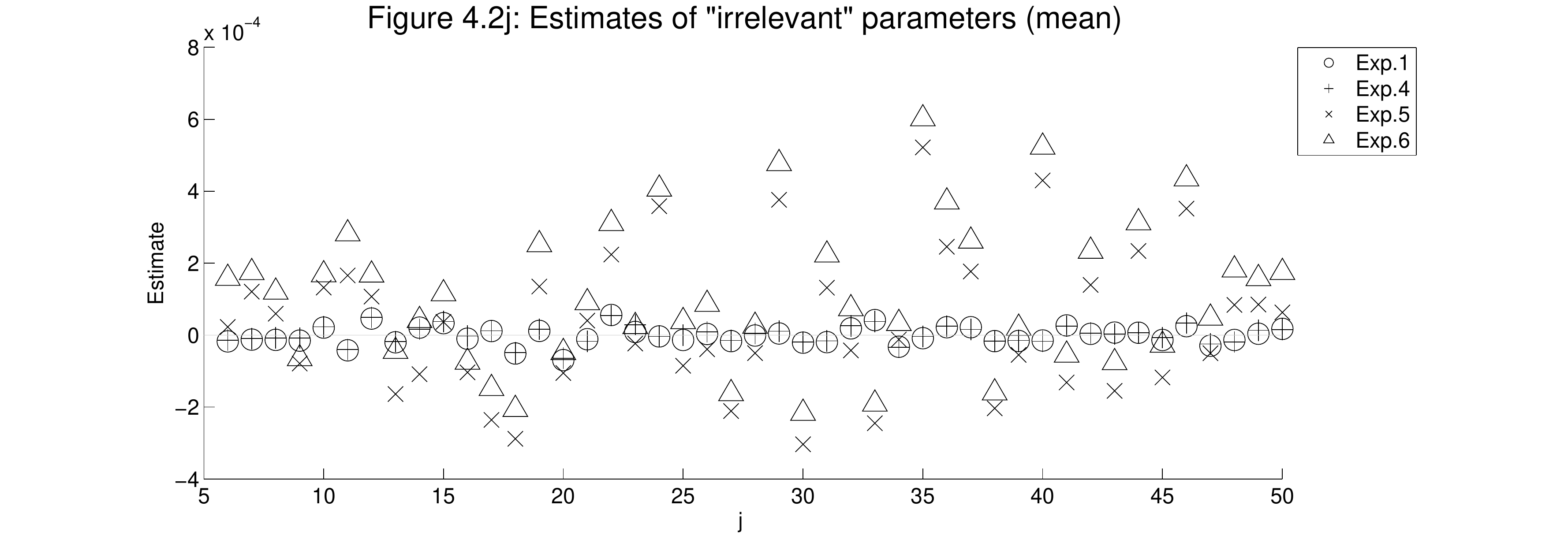}\tabularnewline
\end{tabular}\\
\\
{\scriptsize }%
\begin{tabular}{llllllllll|c|c|c|c|c|}
 &  &  &  &  &  &  &  &  & \multicolumn{1}{l}{} & \multicolumn{5}{c}{{\scriptsize Table 4.3: $l_{2}-$errors and selection (Exp. 1, 4-6)}}\tabularnewline
\cline{11-15} 
 &  &  &  &  &  &  &  &  &  & {\scriptsize Mean} & \multicolumn{4}{c|}{{\scriptsize $n=4700$}}\tabularnewline
\cline{11-15} 
 &  &  &  &  &  &  &  &  &  & {\scriptsize Exp. \# } & {\scriptsize 1} & {\scriptsize 4} & {\scriptsize 5} & {\scriptsize 6}\tabularnewline
\cline{11-15} 
 &  &  &  &  &  &  &  &  &  & {\scriptsize $2^{nd}-$stage select \%} & {\scriptsize 98.1} & {\scriptsize 98.1} & {\scriptsize 55.0} & {\scriptsize 54.5}\tabularnewline
\cline{11-15} 
 &  &  &  &  &  &  &  &  &  & {\scriptsize $2^{nd}-$stage $l_{2}-$err} & {\scriptsize 0.018} & {\scriptsize 0.038} & {\scriptsize 0.062} & {\scriptsize 0.054}\tabularnewline
\cline{11-15} 
 &  &  &  &  &  &  &  &  &  & {\scriptsize $1^{st}$-stage select \%} & {\scriptsize 98.5} & {\scriptsize 52.0} & {\scriptsize 98.5} & {\scriptsize 52.0}\tabularnewline
\cline{11-15} 
 &  &  &  &  &  &  &  &  &  & {\scriptsize $1^{st}-$stage $l_{2}-$err} & {\scriptsize 0.028} & {\scriptsize 0.059} & {\scriptsize 0.028} & {\scriptsize 0.059}\tabularnewline
\cline{11-15} 
\end{tabular}\\
\\
\\
In the next group of experiments which explore the sensitivity of
the results for the two-stage Lasso estimator to design parameters,
changes are applied to $\sigma_{\epsilon}$, $\sigma_{\eta}$, $\sigma_{z}$,
and the correlations between the rows of the design matrix $\mathbf{z}_{i}^{T}\in\mathbb{R}^{p\times d}$
for all $i=1,...,n$. Figure 4.3a plots (in ascending values) the
1000 estimates of $\beta_{5}^{*}$ when the sample size $n=47$. Figure
4.3b plots (in ascending values) the 1000 estimates of $\beta_{6}^{*}$
when the sample size $n=47$. Figures 4.3c-4.3f (Figures 4.3g-4.3j)
displays the $5^{th}$ percentile, the median, the $95^{th}$ percentile,
and the mean of the estimates of the {}``relevant'' ({}``irrelevant'')
main-equation parameters from Experiment 1 and Experiments 7-13. The
mean of the $l_{2}-$errors and the mean of the selection percentages
of the main-equation estimates together with the {}``averaged''
mean of the $l_{2}-$errors and the {}``averaged'' mean of the selection
percentages of the first-stage estimates from these experiments are
displayed in Table 4.4. 

Overall, we see from Table 4.4 that, relative to the baseline experiment
(Experiment 1), the mean of the $l_{2}-$errors of the estimates of
the main-equation parameters increase in Experiments 7-13; the mean
of the selection percentages of the estimates of the main-equation
parameters decrease in Experiments 7, 8, 10-13. The {}``averaged''
mean of the $l_{2}-$errors of the first-stage estimates increase
the most in Experiments 8, 9, 12, and 13 while those first-stage statistics
in Experiment 10 are comparable to those in Experiment 1. This makes
sense since Experiments 8, 9, 12, and 13 involve increasing the noise
level $\sigma_{\eta}$ or decreasing the signal level $\sigma_{z}$
of the instruments in the first-stage model while introducing correlations
between the rows of the design matrix $\mathbf{z}_{i}^{T}$ (for $i=1,...,n$)
(Experiment 10) should have little impact on the first-stage estimates,
which are obtained by performing the Lasso procedure on each of the
50 first-stage equations separately. Note that since Experiment 7
(Experiment 11) has exactly the same first-stage set up as Experiment
1 (respectively, Experiment 10), there is no need to look at the behavior
of their first-stage estimates separately. 

From Figure 4.3a, we see that, below the $10^{th}$ percentile, compared
to the baseline experiment (Experiment 1), introducing correlations
between the rows of the design matrix $\mathbf{z}_{i}^{T}$ (for $i=1,...,n$)
improves the estimates of the {}``relevant'' main-equation parameters
while the other changes to the data generating process yield worse
estimates; above the $80^{th}$ percentile, Figure 4.3a shows that,
any changes made to the data generating process yield worse estimates
of the {}``relevant'' main-equation parameters; at the median (or
mean), Figures 4.3a and 4.3d (respectively, Figure 4.3f) show that,
increasing the standard deviation of the {}``noise'' in the main
equation, $\sigma_{\epsilon}$, and reducing the standard deviation
of the {}``signal'', $\sigma_{z}$, yield worse estimates of the
{}``relevant'' main-equation parameters. Figures 4.3b and 4.3j show
that, in estimating the {}``irrelevant'' main-equation parameters,
any changes to the data generating process yield worse estimates;
in particular, those that involve introducing correlations between
the rows of the design matrix $\mathbf{z}_{i}^{T}$ (for $i=1,...,n$)
(Experiments 10-13) yield the poorest estimates of the {}``irrelevant''
main-equation parameters (also see Table 4.4 for a comparison between
the selection percentages of these experiments). Recall that each
row of $\mathbf{z}_{i}^{T}\in\mathbb{R}^{p\times d}$ is associated
with each of the endogenous regressors and row-wise correlations in
$\mathbf{z}_{i}^{T}$ hence introduce correlations between the {}``purged''
regressors $X_{j}^{*}$ and $X_{j^{'}}^{*}$ for all $j\neq j^{'}$.
As seen in Section 3.2, selection consistency hinges on Assumption
3.7 (the {}``mutual incoherence condition''), whose violation can
lead the Lasso to falsely include elements that are not in the support
of $\beta^{*}$, namely, the violation of Part (b) in Theorems 3.7
and 3.8. For the baseline experiment (Experiment 1), the quantity
$\left\Vert \mathbb{E}\left[X_{1,J(\beta^{*})^{c}}^{*T}X_{1,J(\beta^{*})}^{*}\right]\left[\mathbb{E}(X_{1,J(\beta^{*})}^{*T}X_{1,J(\beta^{*})}^{*})\right]^{-1}\right\Vert _{\infty}$
in Assumption 3.7 equals $0$ (because $\mathbf{z}_{i}^{T}$ is a
$p\times d$ matrix of independent standard normal random variables)
and Assumption 3.7 is easily satisfied. For Experiments 10-13, this
quantity increases, and therefore in these experiments, the estimates
of the {}``irrelevant'' main-equation parameters are worse relative
to the baseline experiment. \\
\\
\\
\begin{tabular}{c}
\includegraphics[width=20cm,height=16cm]{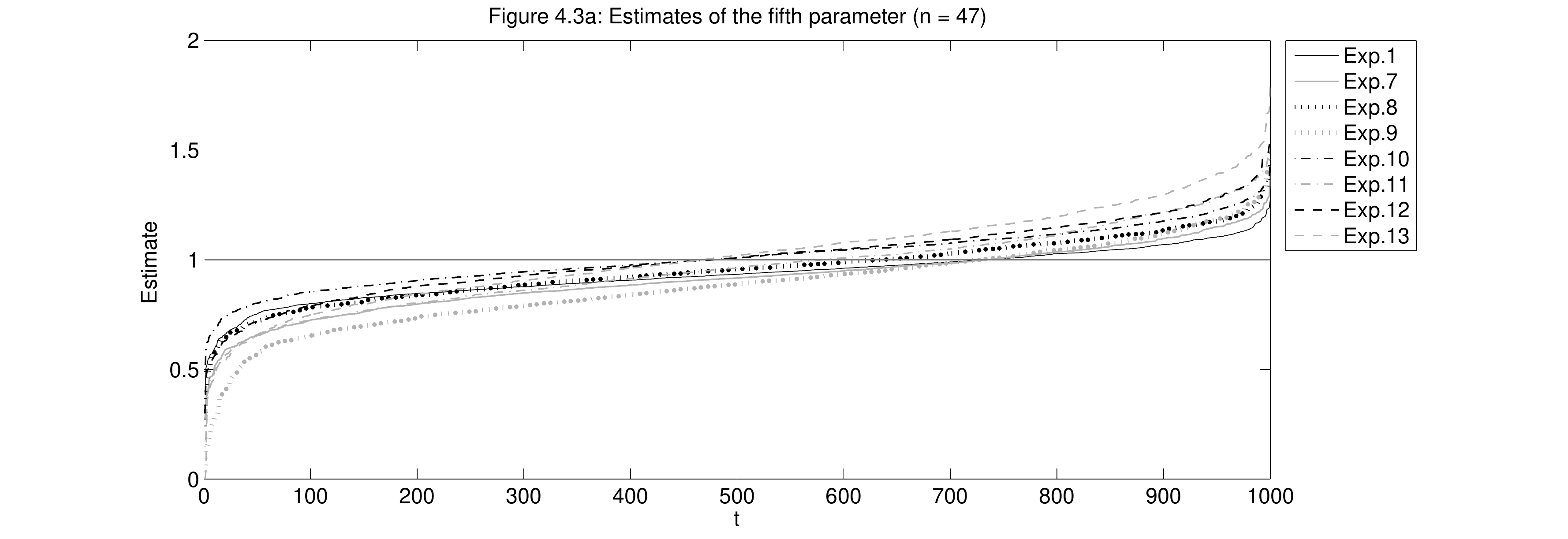}\tabularnewline
\end{tabular}\\
\begin{tabular}{c}
\includegraphics[width=20cm,height=16cm]{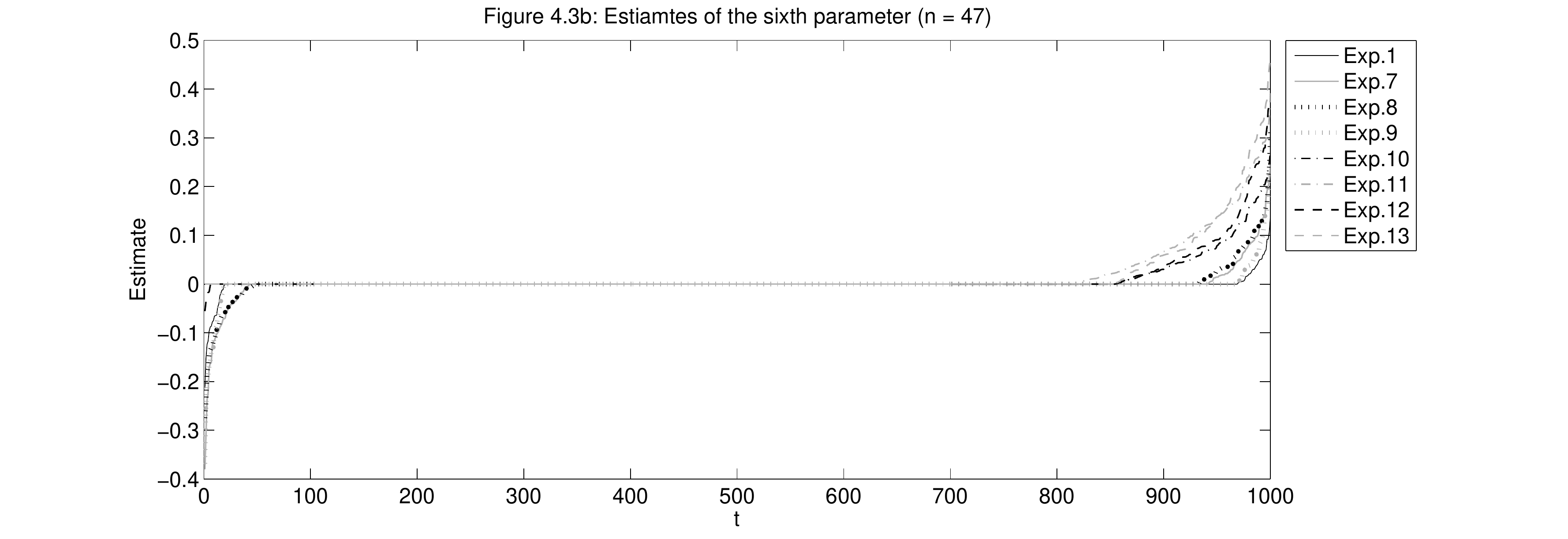}\tabularnewline
\end{tabular}\\
\begin{tabular}{cc}
\includegraphics[width=8cm,height=8cm]{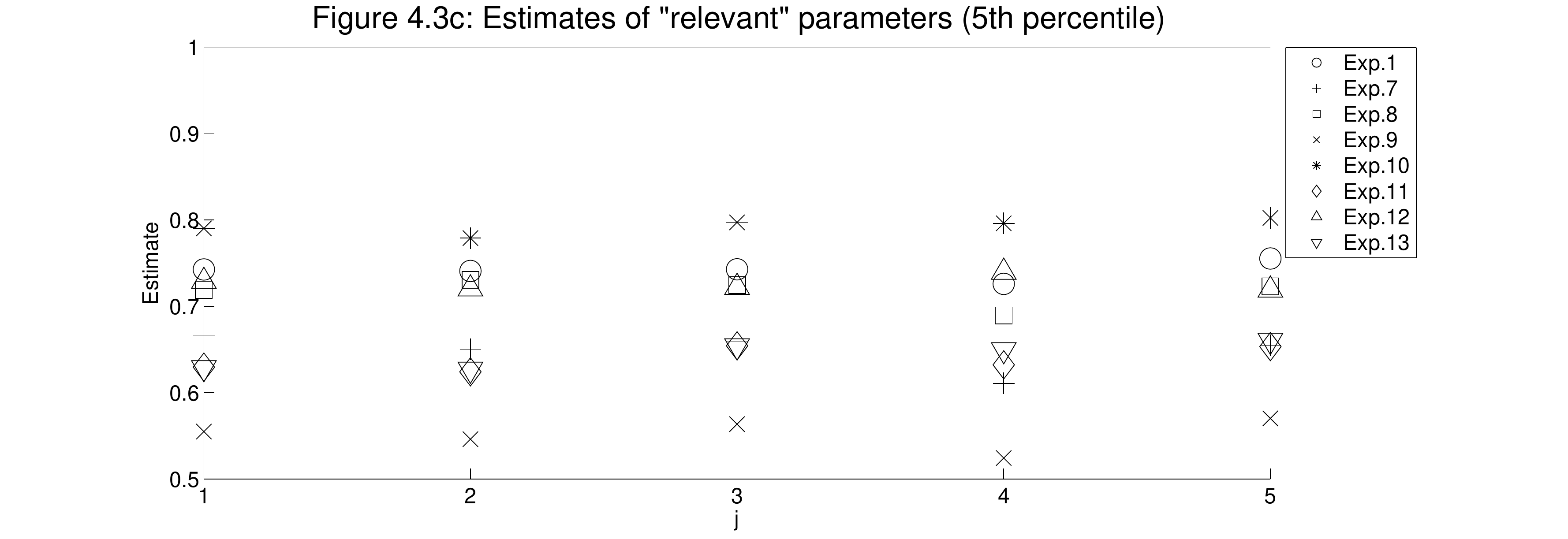} & \includegraphics[width=8cm,height=8cm]{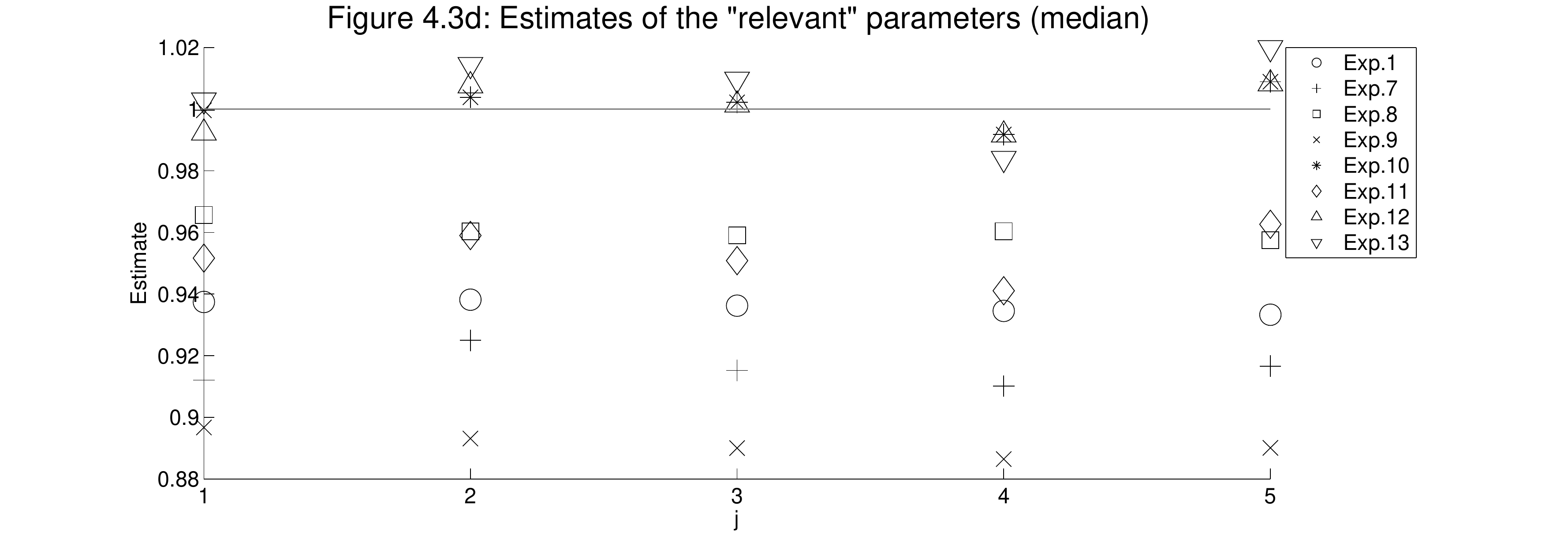}\tabularnewline
\end{tabular}\\
\begin{tabular}{cc}
\includegraphics[width=8cm,height=8cm]{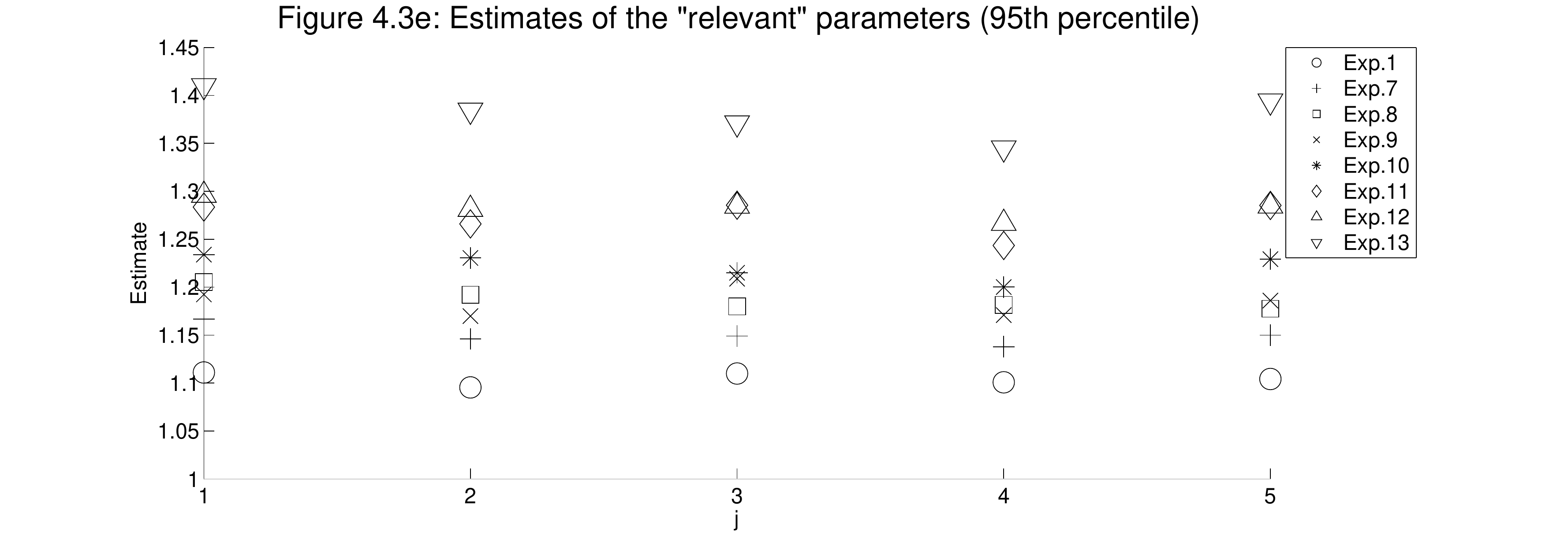} & \includegraphics[width=8cm,height=8cm]{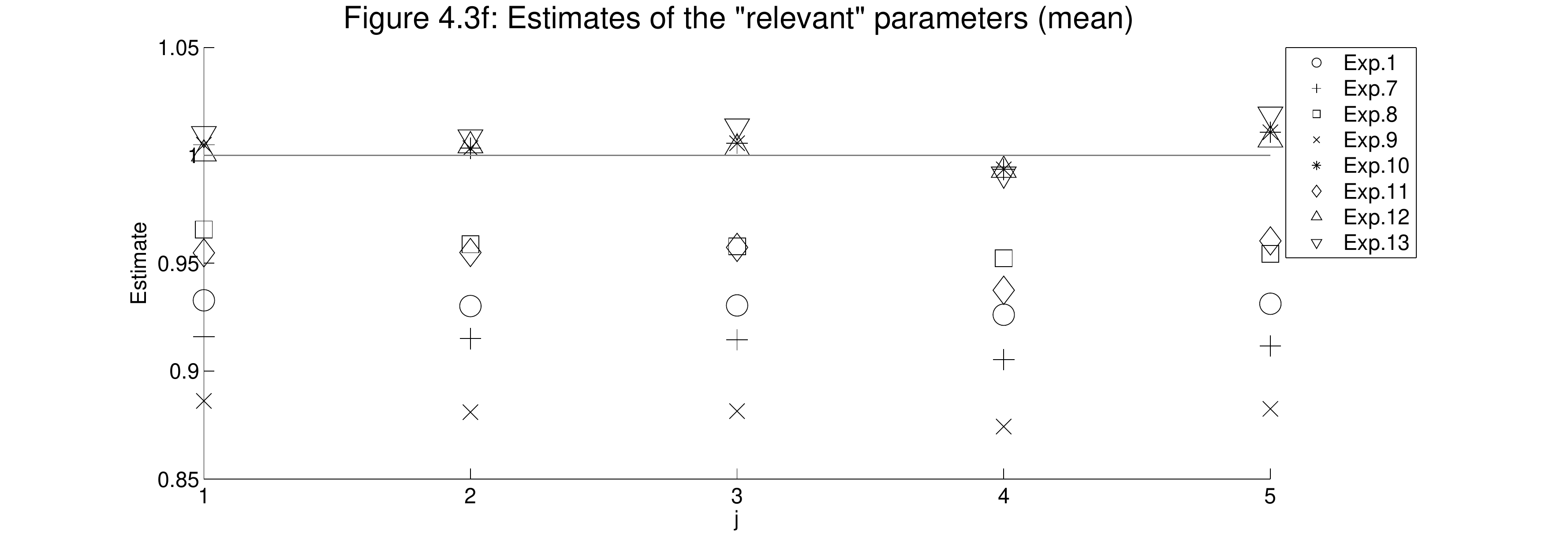}\tabularnewline
\end{tabular}\\
\begin{tabular}{cc}
\includegraphics[width=8cm,height=8cm]{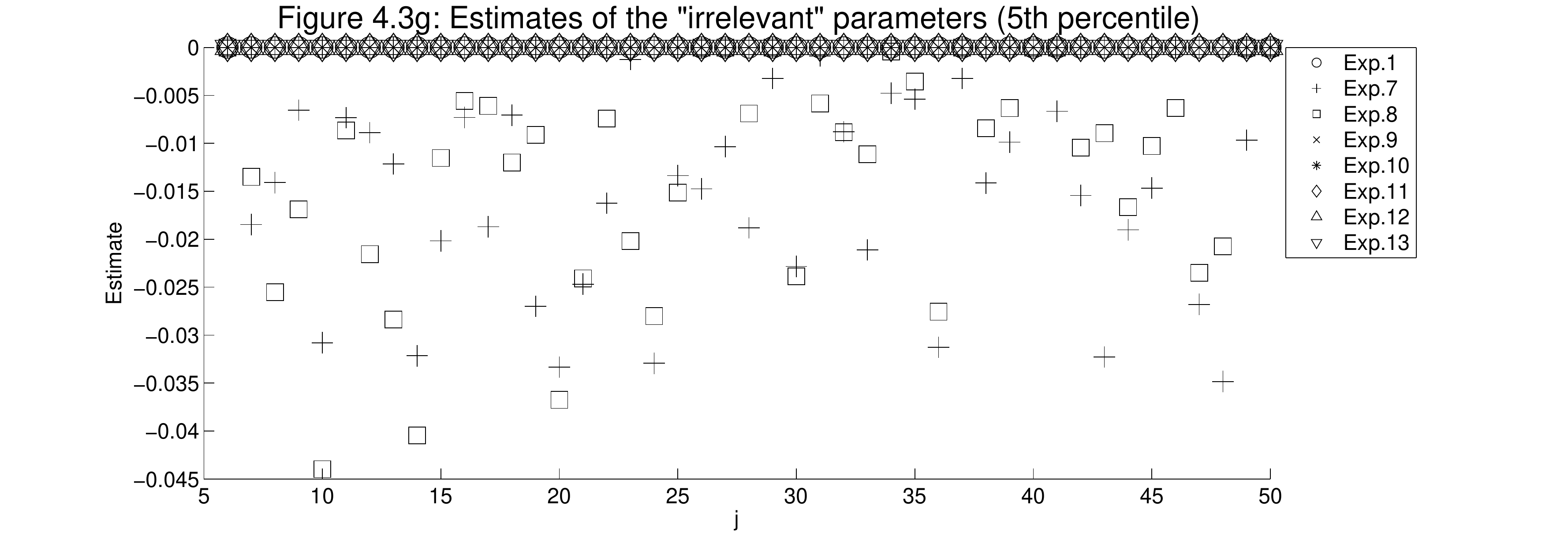} & \includegraphics[width=8cm,height=8cm]{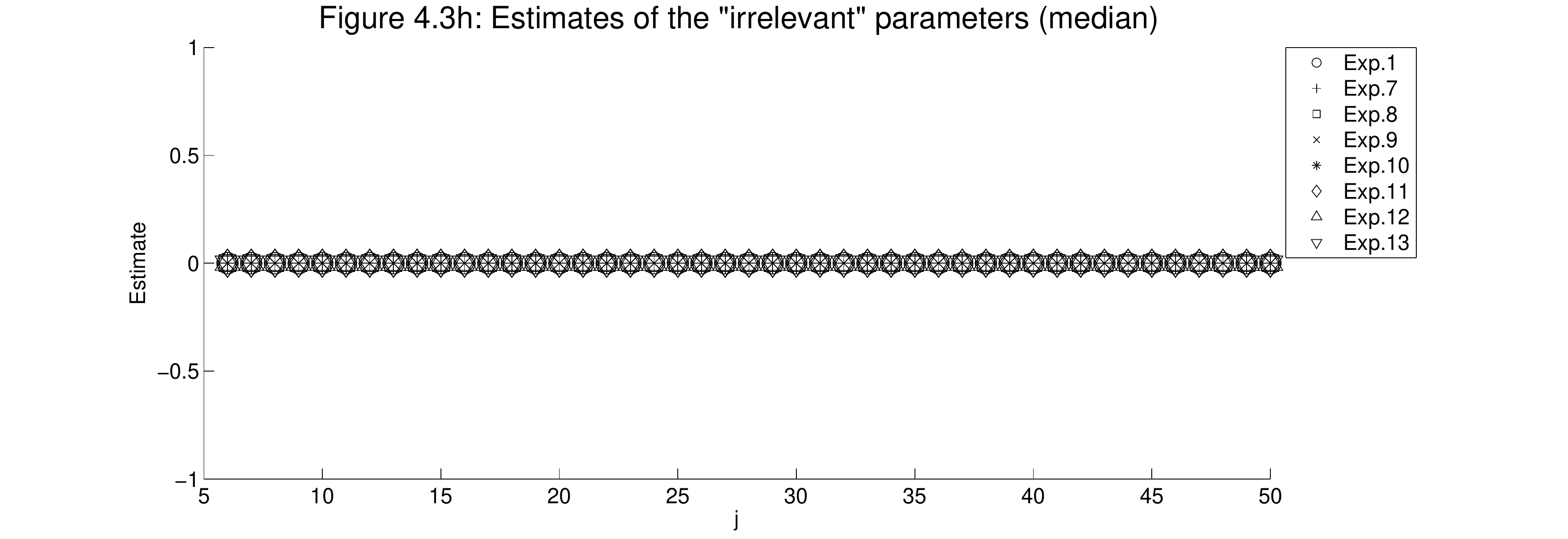}\tabularnewline
\end{tabular}\\
\begin{tabular}{cc}
\includegraphics[width=8cm,height=8cm]{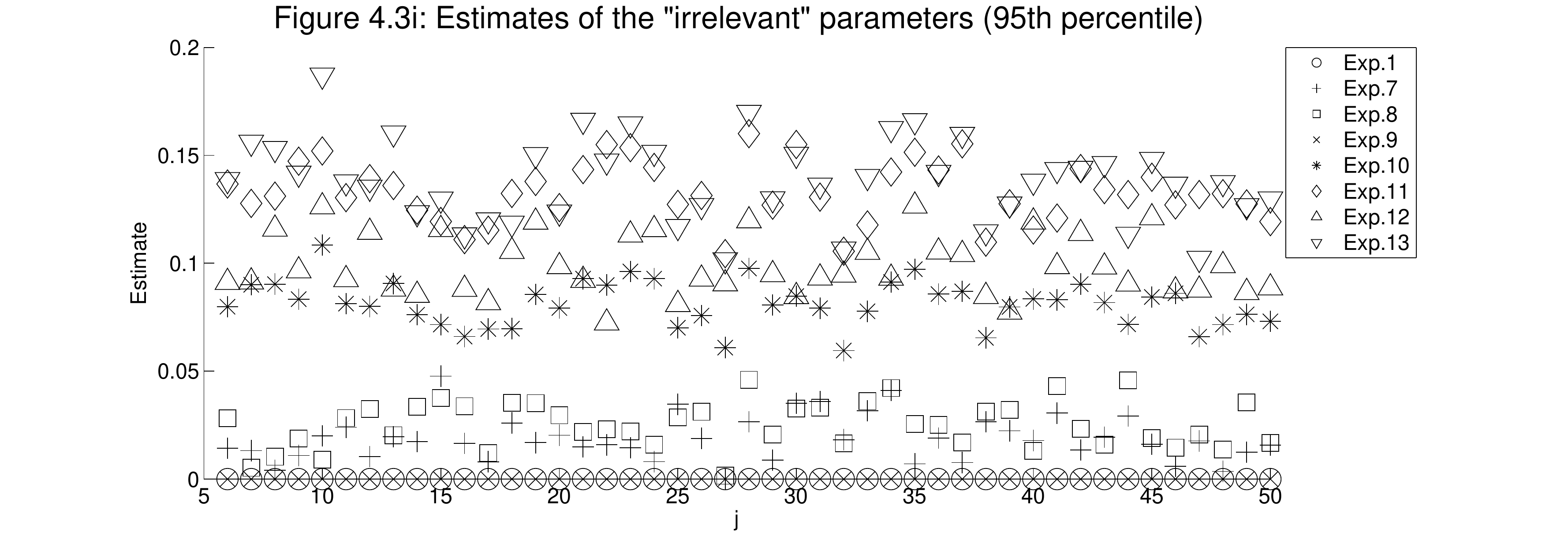} & \includegraphics[width=8cm,height=8cm]{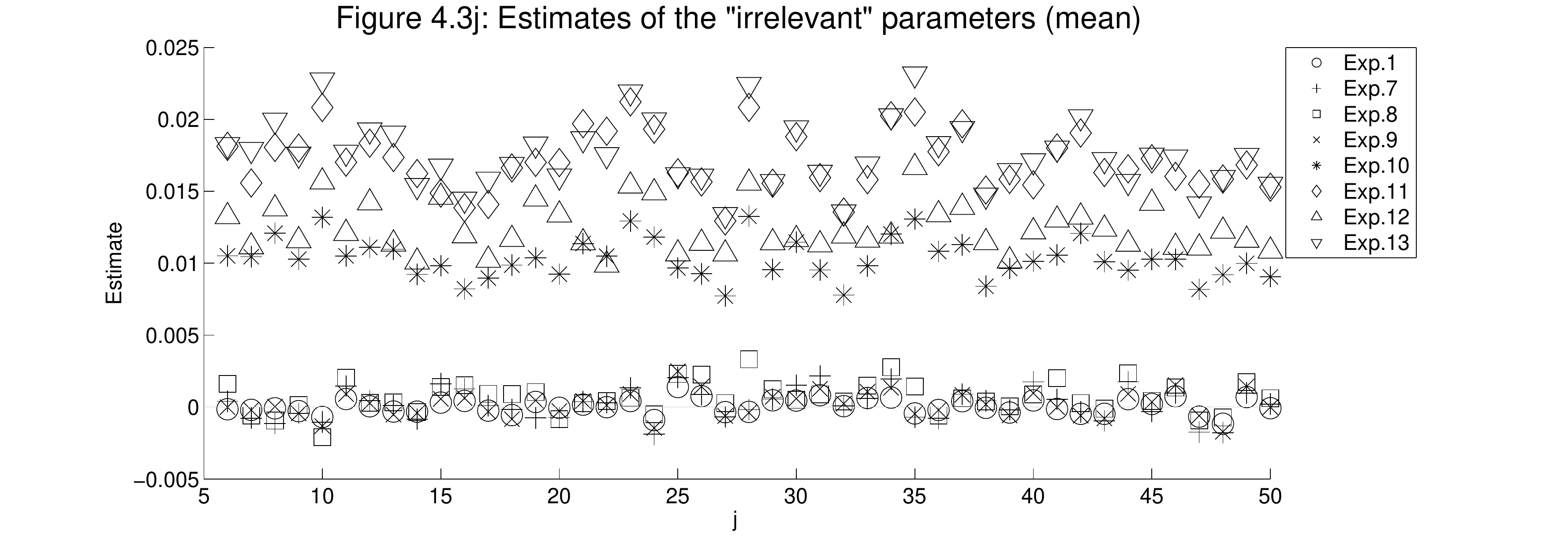}\tabularnewline
\end{tabular}\\
\\
{\scriptsize }%
\begin{tabular}{lllll|c|c|c|c|c|c|c|c|c|}
 &  &  &  & \multicolumn{1}{l}{} & \multicolumn{9}{c}{{\scriptsize Table 4.4: $l_{2}-$errors and selection (Exp. 1, 7-13)}}\tabularnewline
\cline{6-14} 
 &  &  &  &  & {\scriptsize Mean} & \multicolumn{8}{c|}{{\scriptsize $n=47$}}\tabularnewline
\cline{6-14} 
 &  &  &  &  & {\scriptsize Exp. \# } & {\scriptsize 1} & {\scriptsize 7} & {\scriptsize 8} & {\scriptsize 9} & {\scriptsize 10} & {\scriptsize 11} & {\scriptsize 12} & {\scriptsize 13}\tabularnewline
\cline{6-14} 
 &  &  &  &  & {\scriptsize $2^{nd}-$stage select \%} & {\scriptsize 97.3} & {\scriptsize 94.3} & {\scriptsize 93.8} & {\scriptsize 97.3} & {\scriptsize 87.0} & {\scriptsize 85.4} & {\scriptsize 87.8} & {\scriptsize 87.1}\tabularnewline
\cline{6-14} 
 &  &  &  &  & {\scriptsize $2^{nd}$-stage $l_{2}-$err} & {\scriptsize 0.288} & {\scriptsize 0.422} & {\scriptsize 0.376} & {\scriptsize 0.497} & {\scriptsize 0.365} & {\scriptsize 0.557} & {\scriptsize 0.471} & {\scriptsize 0.626}\tabularnewline
\cline{6-14} 
 &  &  &  &  & {\scriptsize $1^{st}-$stage select \%} & {\scriptsize 97.7} & {\scriptsize 97.7} & {\scriptsize 90.8} & {\scriptsize 97.7} & {\scriptsize 97.7} & {\scriptsize 97.7} & {\scriptsize 90.8} & {\scriptsize 97.7}\tabularnewline
\cline{6-14} 
 &  &  &  &  & {\scriptsize $1^{st}-$stage $l_{2}-$err} & {\scriptsize 0.349} & {\scriptsize 0.349} & {\scriptsize 0.789} & {\scriptsize 0.552} & {\scriptsize 0.352} & {\scriptsize 0.352} & {\scriptsize 0.793} & {\scriptsize 0.557}\tabularnewline
\cline{6-14} 
\end{tabular}\\
\\
\\
In the final experiment 14 where $(\beta_{1}^{*},\,...\,,\beta_{5}^{*})=(0.01,\,...\,,0.01)$
(as opposed to $(\beta_{1}^{*},\,...\,,\beta_{5}^{*})=(1,\,...\,,1)$
in the previous experiments), based on the estimates obtained from
the two-stage Lasso procedure, I count the number of occurrences that
each estimate $\hat{\beta}_{H2SLS,\,1},\,...\,,\hat{\beta}_{H2SLS,\,5}$
equals exactly $0$, respectively, over the 1000 replications (Table
4.5). Because the {}``relevant'' main-equation parameters are reduced
by a factor of 100, it is clearly more difficult for the two-stage
Lasso procedure to distinguish the {}``relevant'' coefficients from
the {}``irrelevant'' coefficients and Table 4.5 verifies this. Recall
in Experiments 10-13, by introducing correlations between the {}``purged''
regressors $X_{j}^{*}$ and $X_{j^{'}}^{*}$ for all $j\neq j^{'}$,
the estimates of the {}``irrelevant'' main-equation parameters become
worse. On the other hand, making the {}``relevant'' main-equation
parameters sufficiently smaller results in worse estimates of the
{}``relevant'' main-equation parameters. This observation confirms
Part (d) of Theorems 3.7 and 3.8; i.e., the violation of the {}``beta-min''
condition can lead the Lasso to mistake the {}``relevant'' coefficients
for the {}``irrelevant'' coefficients. In terms of the $l_{2}-$errors
and overall selection percentages, from Table 4.5 we see that poorer
estimation of the {}``relevant'' parameters also results in larger
$l_{2}-$errors%
\footnote{To compensate for the fact that $(\beta_{1},...,\beta_{5})=(1,\,...\,,1)$
in Experiment 1 exceeds the parameters $(\beta_{1},...,\beta_{5})=(0.01,\,...\,,0.01)$
in Experiments 14 by a factor of $0.01$, the $l_{2}-$error in Experiment
14 is adjusted as $\left[\frac{\sum_{j=1}^{5}(\hat{\beta}_{j}-\beta_{j}^{*})^{2}}{0.01^{2}}+\sum_{j=6}^{j=50}(\hat{\beta}_{j}-\beta_{j}^{*})^{2}\right]^{1/2}$.
The unadjusted $l_{2}-$error in Experiment 14 is 0.388.%
} and worse selection percentages, as expected. The significant drop
in the overall selection percentages suggests that not only the estimation
of the {}``relevant'' coefficients becomes less accurate in Experiment
14 but also the estimation of the {}``irrelevant'' coefficients.\\
\\
\begin{tabular}{ccccccc|c||c|c|c|c|c||c|c|}
 &  &  &  &  &  & \multicolumn{1}{c}{} & \multicolumn{8}{c}{{\scriptsize Table 4.5: Exp. 14}}\tabularnewline
\cline{8-15} 
 &  &  &  &  &  &  & \multirow{2}{*}{{\scriptsize $n=47$}} & \multicolumn{5}{c||}{{\scriptsize \# of zeros }} & \multirow{2}{*}{{\scriptsize $2^{nd}-$stg select \%}} & \multirow{2}{*}{{\scriptsize $2^{nd}-$stg $l_{2}-$err}}\tabularnewline
\cline{9-13} 
 &  &  &  &  &  &  &  & {\scriptsize $\beta_{1}$} & {\scriptsize $\beta_{2}$} & {\scriptsize $\beta_{3}$} & {\scriptsize $\beta_{4}$} & {\scriptsize $\beta_{5}$} &  & \tabularnewline
\cline{8-15} 
 &  &  &  &  &  &  & {\scriptsize Exp. 1} & {\scriptsize 0} & {\scriptsize 0} & {\scriptsize 0} & {\scriptsize 0} & {\scriptsize 0} & {\scriptsize 97.3} & {\scriptsize 0.288}\tabularnewline
\cline{8-15} 
 &  &  &  &  &  &  & {\scriptsize Exp. 14} & {\scriptsize 187} & {\scriptsize 187} & {\scriptsize 218} & {\scriptsize 194} & {\scriptsize 193} & {\scriptsize 57.7} & {\scriptsize 11.5}\tabularnewline
\cline{8-15} 
\end{tabular}\\

\section{Conclusion and extensions}

This paper has explored the validity of the two-stage estimation procedure
for sparse linear models in high-dimensional settings with possibly
many endogenous regressors. In particular, the number of endogenous
regressors in the main equation and the number of instruments in the
first-stage equations are permitted to grow with and exceed $n$.
Sufficient scaling conditions on the sample size for estimation consistency in $l_{2}-$ norm 
and variable-selection consistency of the high-dimensional two-stage estimators have been established. \textcolor{black}{I provide theoretical
justifications to a technical issue (regarding the RE condition and
the MI condition) that arises in the two-stage estimation procedure
from allowing the number of regressors in the main equation
to grow with and exceed $n$.} Depending on the underlying assumptions
that are imposed, the upper bounds on the $l_{2}-$error and the sample
size required to obtain these consistency results differ by factors
involving the sparsity parameters $k_{1}$ and/or $k_{2}$. Simulations
are conducted to gain insight on the finite sample performance of
the high-dimensional two-stage estimator.\textbf{ }

The approach and results of this paper suggest a number of possible
extensions including the ones listed in the following, which are left
to future research.\\
\textbf{}\\
\textbf{\textit{Revisiting the bound in Theorem 3.2}}. As discussed
earlier, Assumption 3.6 can be interpreted as a sparsity constraint
on the first-stage estimate\textbf{ $\hat{\pi}_{j}$ }for\textbf{
$j=1,...,p$}, in terms of the $l_{0}$- ball, given by 
\[
\mathbb{B}_{0}^{d_{j}}(k_{1}):=\left\{ \hat{\pi}_{j}\in\mathbb{R}^{d_{j}}\,\vert\,\sum_{l=1}^{d_{j}}1\{\hat{\pi}_{jl}\neq0\}\leq k_{1}\right\} \,\textrm{for }j=1,...,p.
\]
The sparsity constraint (namely, the selection consistency) regarding
these first-stage estimates is guaranteed under some conditions that
may be violated in many problems. It seems possible to extend Assumption
3.6 to the following \textit{approximate sparsity} constraint on the
first-stage estimates\textbf{ }in terms of $l_{1}$- balls, given
by
\[
\mathbb{B}_{1}^{d_{j}}(R_{j}):=\left\{ \hat{\pi}_{j}\in\mathbb{R}^{d_{j}}\,\vert\,|\hat{\pi}_{j}|_{1}=\sum_{l=1}^{d_{j}}|\hat{\pi}_{jl}|\leq R_{j}\right\} \,\textrm{for }j=1,...,p.
\]
If the first-stage estimation employs the Lasso or Dantzig selector
or some other procedures with the $l_{1}-$type of regularization,
then we are guaranteed to have $\hat{\pi}_{j}\in\mathbb{B}_{1}^{d_{j}}(R_{j})$
for every $j=1,...,p$. Depending on the type of sparsity assumptions
imposed on the first-stage estimates, the statistical error of the
high-dimensional two-stage estimator $\hat{\beta}_{H2SLS}$ in \textbf{$l_{2}-$}norm\textbf{
}and the required sample size differ. 

An inspection of the proof for Theorem 3.2 suggests that the error
bound and requirement of the sample size in Theorem 3.2 will hold
regardless of the sparsity assumption on the first-stage estimates.\textbf{
}However, under these special structures that impose a certain decay
rate on the ordered entries of the first-stage estimates, the bound
and scaling of the required sample size in Theorem 3.2 is likely to
be suboptimal.\textbf{ }To obtain sharper results, the proof technique
adopted for showing Theorem 3.3 seems more appropriate.\textbf{ }I
give a heuristic truncation argument to illustrate how the proof for
Theorem 3.3 might be extended to allow the weaker sparsity constraint
(in terms of $l_{1}-$balls) on the first-stage Lasso estimates. 

Suppose for every $j=1,...,p$, we choose the top $s^{j}$ coefficients
of $\hat{\pi}_{j}$ in absolute value, then the fast decay imposed
by the $l_{1}$- ball condition on $\hat{\pi}_{j}$ arising from the
Lasso procedure would mean that the remaining $d_{j}-s^{j}$ coefficients
would have relatively little impact. With this intuition, the proof
follows as if Assumption 3.6 were imposed with the only exception
that we also need to take into account the approximation error arising
from the the remaining $d_{j}-s^{j}$ coefficients of $\hat{\pi}_{j}$.\\
\textbf{}\\
\textbf{The }\textbf{\textit{approximate sparsity}}\textbf{ case}.
It is useful to extend the analysis for the high-dimensional 2SLS
estimator to the \textit{approximate sparsity} case, i.e., most of
the coefficients in the main equation and/or the first-stage equations
are too small to matter. One can have the approximate sparsity assumption
in the first-stage equations only (and assume the main equation parameters
are sparse), the main equation only (and assume the first-stage equations
parameters are sparse) or both-stage equations. When the first-stage
equations parameters are approximately sparse, the argument in the
proof for Theorem 3.2 can still be carried through while the proof
for Theorem 3.3 is no longer meaningful. \\
\textbf{}\\
\textbf{\textit{Control function approach in high-dimensional settings}}.
As an alternative to the {}``two-stage'' estimation proposed here,
it would be interesting to explore the validity of the high-dimensional
two-stage estimators based on the {}``control function'' approach
in the high-dimensional setting. When both the first and second-stage
equations are in low-dimensional settings (i.e., $p\ll n$ and $d_{j}\ll n$
for all $j=1,...,p$) and the supports of the true parameters in both
stages are known \textit{a priori}, the 2SLS procedure is algebraically
equivalent to a {}``control function'' estimator of $\beta^{*}$
that includes first-stage residuals $\hat{\eta}_{ij}=x_{ij}-\mathbf{z}_{ij}^{T}\hat{\pi}_{j}$
as {}``control variables'' in the regression of $y_{i}$ on $\mathbf{x}_{i}$
(e.g., Garen, 1984). Such algebraic equivalence no longer holds for
regularized estimators because the regularization employed destroys
the projection algebra. The extension for the 2SLS estimator from
low-dimensional settings to high-dimensional settings is somewhat
more natural than the extension for the two-stage estimator based
on the control function approach. One question to ask is: under what
conditions can we translate the sparsity or approximate sparsity assumption
on the coefficients $\beta^{*}$ in the triangular simultaneous equations
model (1) and (2) to the sparsity or approximate sparsity assumption
on the coefficients $\beta^{*}$ and $\alpha^{*}$ in the model $y_{i}=\mathbf{x}_{i}^{T}\beta^{*}+\mathbf{\boldsymbol{\eta}}_{i}\alpha^{*}+\xi_{i}$
where $\mathbb{E}(\boldsymbol{\eta}_{i}\xi_{i})=\mathbb{E}(\mathbf{x}_{i}\xi_{i})=\mathbf{0}$?
A simple sufficient condition for such a translation is to impose
the joint normality assumption of the error terms $\epsilon_{i}$
and $\boldsymbol{\eta}_{i}=(\eta_{i1},...,\eta_{ip})$. Then, by the
property of multivariate normal distributions, we have 
\[
\mathbb{E}(\epsilon_{i}\vert\boldsymbol{\eta}_{i})=\Sigma_{\epsilon\eta}\Sigma_{\eta\eta}^{-1}\boldsymbol{\eta}_{i}^{T}.
\]
If we further assume only a few of the correlation coefficients $(\rho_{\epsilon_{i}\eta_{i1}},...,\rho_{\epsilon_{i}\eta_{ip}})$
(associated with the covariance matrix $\Sigma_{\epsilon\eta}$) are
non-zero or most of these correlation coefficients are too small to
matter, the sparsity or approximate sparsity can be carried to the
model $y_{i}=\mathbf{x}_{i}^{T}\beta^{*}+\mathbf{\boldsymbol{\eta}}_{i}\alpha^{*}+\xi_{i}$.
Then, we can obtain consistent estimates of $\eta$, $\hat{\eta}$,
from the first-stage regression by either a standard least square
estimator when the first-stage regression concerns a small number
of regressors relative to $n$, or a least square estimator with $l_{1}$-
regularization (the Lasso or Dantzig selector) when the first-stage
regression concerns a large number of regressors relative to $n$,
and then apply a Lasso technique in the second stage as follows 
\[
\hat{\beta}_{HCF}\in\textrm{argmin}_{\beta,\alpha\in\mathbb{R}^{p}}:\:\frac{1}{2n}|y-X\beta-\hat{\eta}\alpha|_{2}^{2}+\lambda_{n}\left(|\beta|_{1}+|\alpha|_{1}\right).
\]
The statistical properties of $\hat{\beta}_{HCF}$ can be analyzed
in the same way as those of $\hat{\beta}_{H2SLS}$. How this argument
can be extended to non-Gaussian error settings is an interesting question
for future research. \\
 \\
\textbf{\textit{Minimax lower bounds for the high-dimensional linear
models with endogeneity}}. It would be worthwhile to establish the
minimax lower bounds on the parameters in the main equation for the
linear models in high-dimensional settings with endogeneity. In particular,
the goal is to derive lower bounds on the estimation error achievable
by any estimator, regardless of its computational complexity. Obtaining
lower bounds of this type is useful because on one hand, if the lower
bound matches the upper bound up to some constant factors, then there
is no need to search for estimators with a lower statistical error
(although it might still be useful to study estimators with lower
computational costs). On the other hand, if the lower bound does not
match the best known upper bounds, then it is worthwhile to search
for new estimators that potentially achieve the lower bound. To the
best of my knowledge, in econometric literature, there has been only
limited attention given to the minimax rates of linear models with
endogeneity in high-dimensional settings. 

\newpage{}

\section{Appendix: Proofs}

For technical simplifications, in the following proofs, I assume without
loss of generality that the first moment of $(y_{i},\,\mathbf{x}_{i},\,\mathbf{z}_{i})$
is zero for all $i$ (if it is not the case, we can simply subtract
their population mean). Also, for notational simplicity, assume $d_{j}=d$
for all $j=1,...,p$; additionally, as in most high-dimensional statistics
literature, I assume the regime of interest is $p\geq n$ and $d\geq n$
(except for Corollary 3.4 where $d\ll n$ is assumed). The modification
to allow $p<n$ or $d<n$ or $d_{j}\neq d_{j^{'}}$ for some $j$
and $j^{'}$ is trivial. Also, as a general rule for the proofs, $b$
constants denote positive constants that do not involve $n$, $p$,
$d$, $k_{1}$ and $k_{2}$ but possibly the sub-Gaussian parameters
defined in Assumptions 3.2-3.4; $c$ constants denote universal positive
constants that are independent of both $n$, $p$, $d$, $k_{1}$
and $k_{2}$ as well as the sub-Gaussian parameters. The specific
values of these constants may change from place to place.

\subsection{Lemma 3.1}

\textbf{Proof}. First, write 
\begin{eqnarray*}
y & = & X\beta^{*}+\epsilon=X^{*}\beta^{*}+(X\beta^{*}-X^{*}\beta^{*}+\epsilon)\\
 & = & X^{*}\beta^{*}+(\boldsymbol{\eta}\beta^{*}+\epsilon)\\
 & = & \hat{X}\beta^{*}+(X^{*}-\hat{X})\beta^{*}+\boldsymbol{\eta}\beta^{*}+\epsilon\\
 & = & \hat{X}\beta^{*}+e,
\end{eqnarray*}
where $e:=(X^{*}-\hat{X})\beta^{*}+\boldsymbol{\eta}\beta^{*}+\epsilon$.
Define $\hat{v}^{0}=\hat{\beta}_{H2SLS}-\beta^{*}$ and the Lagrangian
$L(\beta;\,\lambda_{n})=\frac{1}{2n}|y-\hat{X}\beta|_{2}^{2}+\lambda_{n}|\beta|_{1}$.
Since $\hat{\beta}_{H2SLS}$ is optimal, we have 
\[
L(\hat{\beta}_{H2SLS};\,\lambda_{n})\leq L(\beta^{*};\,\lambda_{n})=\frac{1}{2n}|e|_{2}^{2}+\lambda_{n}|\beta^{*}|_{1},
\]
Some algebraic manipulation of the \textit{basic inequality} above
yields 
\begin{eqnarray*}
0 & \leq & \frac{1}{2n}|\hat{X}\hat{v}^{0}|_{2}^{2}\leq\frac{1}{n}e^{T}\hat{X}\hat{v}^{0}+\lambda_{n}\left\{ |\beta_{J(\beta^{*})}^{*}|_{1}-|(\beta_{J(\beta^{*})}^{*}+\hat{v}_{J(\beta^{*})}^{0},\,\hat{v}_{J(\beta^{*})^{c}}^{0})|_{1}\right\} \\
 & \leq & |\hat{v}^{0}|_{1}|\frac{1}{n}\hat{X}^{T}e|_{\infty}+\lambda_{n}\left\{ |\hat{v}_{J(\beta^{*})}^{0}|_{1}-|\hat{v}_{J(\beta^{*})^{c}}^{0}|_{1}\right\} \\
 & \leq & \frac{\lambda_{n}}{2}\left\{ 3|\hat{v}_{J(\beta^{*})}^{0}|_{1}-|\hat{v}_{J(\beta^{*})^{c}}^{0}|_{1}\right\} ,
\end{eqnarray*}
where the last inequality holds as long as $\lambda_{n}\geq2|\frac{1}{n}\hat{X}^{T}e|_{\infty}>0$.
Consequently, $|\hat{v}^{0}|_{1}\leq4|\hat{v}_{J(\beta^{*})}^{0}|_{1}\leq4\sqrt{k_{2}}|\hat{v}_{J(\beta^{*})}^{0}|_{2}\leq4\sqrt{k_{2}}|\hat{v}^{0}|_{2}$.
Note that we also have 
\begin{eqnarray*}
\frac{1}{2n}|\hat{X}\hat{v}^{0}|_{2}^{2} & \leq & |\hat{v}^{0}|_{1}|\frac{1}{n}\hat{X}^{T}e|_{\infty}+\lambda_{n}\left\{ |\hat{v}_{J(\beta^{*})}^{0}|_{1}-|\hat{v}_{J(\beta^{*})^{c}}^{0}|_{1}\right\} \\
 & \leq & 4\sqrt{k_{2}}|\hat{v}^{0}|_{2}\lambda_{n}.
\end{eqnarray*}
Since we assume in Lemma 3.1 that the random matrix $\hat{\Gamma}=\hat{X}^{T}\hat{X}$
satisfies the RE condition (3) with $\gamma=3$, we have 
\[
|\hat{\beta}_{H2SLS}-\beta^{*}|_{2}\leq\frac{c^{'}}{\delta}\sqrt{k_{2}}\lambda_{n}.
\]

\subsection{Theorem 3.2}

As discussed in Section 3, the $l_{2}$-consistency of $\hat{\beta}_{H2SLS}$
requires verifications of two conditions: (i) $\hat{\Gamma}=\hat{X}^{T}\hat{X}$
satisfies the RE condition (3) with $\gamma=3$, and (ii) the term
$|\frac{1}{n}\hat{X}^{T}e|_{\infty}\precsim f(k_{1},\, k_{2},\, d,\, p,\, n)$
with high probability. This is done via Lemmas 6.1 and 6.2.\\
\textbf{}\\
\textbf{Lemma 6.1} (RE condition): Under Assumptions 1.1, 3.1, 3.3,
3.5a and the condition 
\[
n\succsim\max\{k_{1}^{2}\log d,\, k_{1}^{2}\log p\},
\]
we have, for some universal constants $c$, $c_{1}$, and $c_{2}$,
\[
\frac{|\hat{X}v^{0}|_{2}^{2}}{n}\geq\kappa_{1}|v^{0}|_{2}^{2}-c\kappa_{2}k_{1}\sqrt{\frac{\log\max(p,\, d)}{n}}|v^{0}|_{1}^{2},\qquad\textrm{for all }v^{0}\in\mathbb{R}^{p},
\]
with probability at least $1-c_{1}\exp(-c_{2}\log\max(p,\, d))$,
where 
\begin{eqnarray*}
\kappa_{1} & = & \frac{\lambda_{\min}(\Sigma_{X^{*}})}{2},\;\kappa_{2}=\max\left\{ b_{0},\: b_{1},\: b_{2}\right\} ,\\
b_{0} & = & \lambda_{\min}(\Sigma_{X^{*}})\max\left\{ \frac{\sigma_{X^{*}}^{4}}{\lambda_{\min}^{2}(\Sigma_{X^{*}})},\,1\right\} ,\\
b_{1} & = & \max\left\{ \frac{\sigma_{\eta}\max_{j^{'},\, j}|\mathbb{E}(x_{ij^{'}}^{*},\,\mathbf{z}_{ij})|_{\infty}}{\lambda_{\min}(\Sigma_{Z})},\;\frac{\sigma_{\eta}\sigma_{X^{*}}\sigma_{Z}}{\lambda_{\min}(\Sigma_{Z})}\right\} ,\\
b_{2} & = & \max\left\{ \frac{\sigma_{\eta}^{2}\max_{j^{'},\, j}|\mathbb{E}(\mathbf{z}_{ij^{'}},\,\mathbf{z}_{ij})|_{\infty}}{\lambda_{\min}^{2}(\Sigma_{Z})},\:\frac{\sigma_{\eta}^{2}\sigma_{Z}^{2}}{\lambda_{\min}^{2}(\Sigma_{Z})}\right\} .
\end{eqnarray*}
\textbf{}\\
\textbf{Proof. }We have 
\[
\left|v^{0T}\frac{\hat{X}^{T}\hat{X}}{n}v^{0}\right|+\left|v^{0T}\left(\frac{X^{*T}X^{*}-\hat{X}^{T}\hat{X}}{n}\right)v^{0}\right|\geq\left|v^{0T}\frac{X^{*T}X^{*}}{n}v^{0}\right|,
\]
which implies \textbf{
\begin{eqnarray*}
\left|v^{0T}\frac{\hat{X}^{T}\hat{X}}{n}v^{0}\right| & \geq & \left|v^{0T}\frac{X^{*T}X^{*}}{n}v^{0}\right|-\left|v^{0T}\left(\frac{X^{*T}X^{*}-\hat{X}^{T}\hat{X}}{n}\right)v^{0}\right|\\
 & \geq & \left|v^{0T}\frac{X^{*T}X^{*}}{n}v^{0}\right|-\left|\frac{X^{*T}X^{*}-\hat{X}^{T}\hat{X}}{n}\right|_{\infty}\left|v^{0}\right|_{1}^{2}\\
 & \geq & \left|v^{0T}\frac{X^{*T}X^{*}}{n}v^{0}\right|-\left(\left|\frac{X^{*T}(\hat{X}-X^{*})}{n}\right|_{\infty}+\left|\frac{(\hat{X}-X^{*})^{T}\hat{X}}{n}\right|_{\infty}\right)\left|v^{0}\right|_{1}^{2}\\
 & \geq & \left|v^{0T}\frac{X^{*T}X^{*}}{n}v^{0}\right|-\left|\frac{X^{*T}(\hat{X}-X^{*})}{n}\right|_{\infty}\left|v^{0}\right|_{1}^{2}\\
 &  & -\left|\frac{(\hat{X}-X^{*})^{T}X^{*}}{n}\right|_{\infty}\left|v^{0}\right|_{1}^{2}-\left|\frac{(\hat{X}-X^{*})^{T}(\hat{X}-X^{*})}{n}\right|_{\infty}\left|v^{0}\right|_{1}^{2}.
\end{eqnarray*}
}To bound the term $\left|\frac{X^{*T}(\hat{X}-X^{*})}{n}\right|_{\infty}$,
let us first fix $(j^{'},\, j)$\textbf{ }and bound the $(j^{'},\, j)$
element of the matrix\textbf{ }$\frac{X^{*T}(\hat{X}-X^{*})}{n}$.
Notice that 
\begin{eqnarray*}
\left|\frac{1}{n}\mathbf{x}_{j^{'}}^{*T}(\hat{\mathbf{x}}_{j}-\mathbf{x}_{j}^{*})\right| & = & \left|\left(\frac{1}{n}\sum_{i=1}^{n}x_{ij^{'}}^{*}\mathbf{z}_{ij}\right)(\hat{\pi}_{j}-\pi_{j}^{*})\right|\\
 & \leq & \left|\hat{\pi}_{j}-\pi_{j}^{*}\right|_{1}\left|\frac{1}{n}\sum_{i=1}^{n}x_{ij^{'}}^{*}\mathbf{z}_{ij}\right|_{\infty}.
\end{eqnarray*}
Under Assumptions 3.2 and 3.3, we have that the random matrix\textbf{
$Z_{j}\in\mathbb{R}^{n\times d_{j}}$ }is\textbf{ }a sub-Gaussian
with parameters at most $(\Sigma_{Z_{j}},\,\sigma_{Z}^{2})$ for all
$j=1,...,p$, and $x_{j^{'}}^{*}$ is a sub-Gaussian vector with a
parameter at most $\sigma_{X^{*}}$ for every $j^{'}=1,...,p$. Therefore,
by Lemma 6.8 and an application of union bound,\textbf{ }we have 
\[
\mathbb{P}\left[\max_{j^{'},\, j}|\frac{1}{n}\mathbf{x}_{j^{'}}^{*T}Z_{j}-\mathbb{E}(x_{ij^{'}}^{*},\,\mathbf{z}_{ij})|{}_{\infty}\geq t\right]\leq6p^{2}d\exp(-cn\min\{\frac{t^{2}}{\sigma_{X^{*}}^{2}\sigma_{Z}^{2}},\,\frac{t}{\sigma_{X^{*}}\sigma_{Z}}\}),
\]
so as long as $n\succsim\log\max(p,\, d)$, 

\[
\mathbb{P}\left[\max_{j^{'},\, j}|\frac{1}{n}\mathbf{x}_{j^{'}}^{*T}Z_{j}-\mathbb{E}(x_{ij^{'}}^{*},\,\mathbf{z}_{ij})|{}_{\infty}\geq c_{0}\sigma_{X^{*}}\sigma_{Z}\sqrt{\frac{\log\max(p,\, d)}{n}}\right]\leq c_{1}\exp(-c_{2}\log\max(p,\, d)),
\]
where\textbf{ }$c_{0}$, $c_{1}$, and $c_{2}$ are some universal
constants. Under Assumption 3.5a, if $n\succsim\log\max(p,\, d)$,
then, \textbf{
\begin{eqnarray*}
\left|\frac{X^{*T}(\hat{X}-X^{*})}{n}\right|_{\infty} & \leq & \frac{c\sigma_{\eta}}{\lambda_{\min}(\Sigma_{Z})}k_{1}\sqrt{\frac{\log\max(p,\, d)}{n}}\left(\max_{j^{'},\, j}|\mathbb{E}(x_{ij^{'}}^{*},\,\mathbf{z}_{ij})|_{\infty}+c_{0}\sigma_{X^{*}}\sigma_{Z}\sqrt{\frac{\log\max(p,\, d)}{n}}\right)\\
 & \leq & c_{3}\max\left\{ \frac{\sigma_{\eta}\max_{j^{'},\, j}|\mathbb{E}(x_{ij^{'}}^{*},\,\mathbf{z}_{ij})|_{\infty}}{\lambda_{\min}(\Sigma_{Z})},\;\frac{\sigma_{\eta}\sigma_{X^{*}}\sigma_{Z}}{\lambda_{\min}(\Sigma_{Z})}\right\} k_{1}\sqrt{\frac{\log\max(p,\, d)}{n}},
\end{eqnarray*}
}with probability at least $1-c_{1}\exp(-c_{2}\log\max(p,\, d))$. 

To bound the term $\left|\frac{(\hat{X}-X^{*})^{T}(\hat{X}-X^{*})}{n}\right|_{\infty}$,
again let us first fix $(j^{'},\, j)$\textbf{ }and bound the $(j^{'},\, j)$
element of the matrix\textbf{ }$\frac{(\hat{X}-X^{*})^{T}(\hat{X}-X^{*})}{n}$.
Using the similar argument as above, if $n\succsim\log\max(p,\, d)$,
we have, 
\begin{eqnarray*}
\left|\frac{(\hat{X}-X^{*})^{T}(\hat{X}-X^{*})}{n}\right|_{\infty} & = & \max_{j^{'},\, j}\left|(\hat{\pi}_{j^{'}}-\pi_{j^{'}}^{*})^{T}\left(\frac{1}{n}\sum_{i=1}^{n}\mathbf{z}_{ij^{'}}^{T}\mathbf{z}_{ij}\right)(\hat{\pi}_{j}-\pi_{j}^{*})\right|\\
 & \leq & \max_{j^{'},\, j}\left(\left|\hat{\pi}_{j^{'}}-\pi_{j^{'}}^{*}\right|_{1}\left|\hat{\pi}_{j}-\pi_{j}^{*}\right|_{1}\left|\frac{1}{n}\sum_{i=1}^{n}\mathbf{z}_{ij^{'}}^{T}\mathbf{z}_{ij}\right|_{\infty}\right)\\
 & \leq & \left(\frac{c\sigma_{\eta}}{\lambda_{\min}(\Sigma_{Z})}k_{1}\sqrt{\frac{\log\max(p,\, d)}{n}}\right)^{2}\left(\max_{j^{'},\, j}|\mathbb{E}(\mathbf{z}_{ij^{'}},\,\mathbf{z}_{ij})|_{\infty}+c_{0}\sigma_{Z}^{2}\sqrt{\frac{\log\max(p,\, d)}{n}}\right)\\
 & \leq & c_{3}\max\left\{ \frac{\sigma_{\eta}^{2}\max_{j^{'},\, j}|\mathbb{E}(\mathbf{z}_{ij^{'}},\,\mathbf{z}_{ij})|_{\infty}}{\lambda_{\min}^{2}(\Sigma_{Z})},\:\frac{\sigma_{\eta}^{2}\sigma_{Z}^{2}}{\lambda_{\min}^{2}(\Sigma_{Z})}\right\} k_{1}^{2}\frac{\log\max(p,\, d)}{n},
\end{eqnarray*}
with probability at least $1-c_{1}\exp(-c_{2}\log\max(p,\, d))$.

Putting everything together, under the condition $n\succsim\max\{k_{1}^{2}\log d,\, k_{1}^{2}\log p\}$
and applying Lemma 6.10 with $r=0$, we have 
\begin{eqnarray*}
\left|v^{0T}\frac{\hat{X}^{T}\hat{X}}{n}v^{0}\right| & \geq & \left|v^{0T}\frac{X^{*T}X^{*}}{n}v^{0}\right|\\
 & - & \left(2c_{3}b_{1}k_{1}\sqrt{\frac{\log\max(p,\, d)}{n}}+c_{4}b_{2}k_{1}^{2}\frac{\log\max(p,\, d)}{n}\right)\left|v^{0}\right|_{1}^{2}\\
 & \geq & \left|v^{0T}\frac{X^{*T}X^{*}}{n}v^{0}\right|-\left(c_{5}\max\{b_{1},\, b_{2}\}k_{1}\sqrt{\frac{\log\max(p,\, d)}{n}}\right)\left|v^{0}\right|_{1}^{2}\\
 & \geq & \frac{\lambda_{\min}(\Sigma_{X^{*}})}{2}\left|v^{0}\right|_{2}^{2}-c_{0}b_{0}\frac{\log\max(p,\, d)}{n}\left|v^{0}\right|_{1}^{2}\\
 & - & \left(c_{5}\max\{b_{1},\, b_{2}\}k_{1}\sqrt{\frac{\log\max(p,\, d)}{n}}\right)\left|v^{0}\right|_{1}^{2},
\end{eqnarray*}
with probability at least $1-c_{1}^{'}\exp(-c_{2}^{'}n)-c_{1}^{''}\exp(-c_{2}^{''}\log\max(p,\, d))=1-c_{1}\exp(-c_{2}\log\max(p,\, d))$
(given $d>n$ and $p>n$ is the regime of our interests), where $b_{0}=\lambda_{\min}(\Sigma_{X^{*}})\max\left\{ \frac{\sigma_{X^{*}}^{4}}{\lambda_{\min}^{2}(\Sigma_{X^{*}})},\,1\right\} $,
$b_{1}=\max\left\{ \frac{\sigma_{\eta}\max_{j^{'},\, j}|\mathbb{E}(x_{ij^{'}}^{*},\,\mathbf{z}_{ij})|_{\infty}}{\lambda_{\min}(\Sigma_{Z})},\;\frac{\sigma_{\eta}\sigma_{X^{*}}\sigma_{Z}}{\lambda_{\min}(\Sigma_{Z})}\right\} $,
and $b_{2}=\max\left\{ \frac{\sigma_{\eta}^{2}\max_{j^{'},\, j}|\mathbb{E}(\mathbf{z}_{ij^{'}},\,\mathbf{z}_{ij})|_{\infty}}{\lambda_{\min}^{2}(\Sigma_{Z})},\:\frac{\sigma_{\eta}^{2}\sigma_{Z}^{2}}{\lambda_{\min}^{2}(\Sigma_{Z})}\right\} $.\textbf{
}Notice the last inequality can be written in the form 
\begin{eqnarray*}
\left|v^{0T}\frac{\hat{X}^{T}\hat{X}}{n}v^{0}\right| & \geq & \kappa_{1}|v^{0}|_{2}^{2}-\kappa_{2}\max\left\{ k_{1}\sqrt{\frac{\log\max(p,\, d)}{n}},\:\frac{\log d}{n},\,\frac{\log p}{n}\right\} |v^{0}|_{1}^{2}\\
 & \geq & \kappa_{1}|v^{0}|_{2}^{2}-\kappa_{2}k_{1}\sqrt{\frac{\log\max(p,\, d)}{n}}|v^{0}|_{1}^{2}
\end{eqnarray*}
where $\kappa_{1}=\frac{\lambda_{\min}(\Sigma_{X^{*}})}{2},\;\kappa_{2}=\max\left\{ b_{0},\: b_{1},\: b_{2}\right\} $,
and the second inequality follows since $n\succsim\log\max(p,\, d)$.
$\square$\\
\\
In proving Lemma 3.1, upon our choice of $\lambda_{n}$, we have shown\textbf{
\[
\hat{v}=\hat{\beta}_{H2SLS}-\beta^{*}\in\mathbb{C}(J(\beta^{*}),\,3),
\]
}which implies $|\hat{v}^{0}|_{1}^{2}\leq16|\hat{v}_{J(\beta^{*})}^{0}|_{1}^{2}\leq16k_{2}|\hat{v}_{J(\beta^{*})}^{0}|_{2}^{2}$.
Therefore, if we have the scaling 
\[
\frac{1}{n}k_{1}^{2}k_{2}^{2}\log\max(p,\, d)=O(1),
\]
so that 
\[
\kappa_{2}k_{1}k_{2}\sqrt{\frac{\log\max(p,\, d)}{n}}<\kappa_{1},
\]
then, \textbf{
\[
\left|\hat{v}^{0T}\frac{\hat{X}^{T}\hat{X}}{n}\hat{v}^{0}\right|\geq c_{0}\lambda_{\min}(\Sigma_{X^{*}})\left|\hat{v}^{0}\right|_{2}^{2},
\]
}provided $\sigma_{\eta}$, $\sigma_{Z}$, $\sigma_{X^{*}}$, $\max_{j^{'},\, j}|\mathbb{E}(x_{ij^{'}}^{*},\,\mathbf{z}_{ij})|_{\infty}$,
and $\max_{j^{'},\, j}|\mathbb{E}(\mathbf{z}_{ij^{'}},\,\mathbf{z}_{ij})|_{\infty}$
are bounded from above while $\lambda_{\min}(\Sigma_{Z})$ and $\lambda_{\min}(\Sigma_{X^{*}})$
are bounded away from $0$. The above inequality implies RE (3).\textbf{}\\
\textbf{}\\
\textbf{Lemma 6.2} (Upper bound on $|\frac{1}{n}\hat{X}^{T}e|_{\infty}$):
Under Assumptions 1.1, 3.1-3.3, 3.5a, and the condition $\frac{\max\{k_{1}^{2}\log d,\, k_{1}^{2}\log p\}}{n}=o(1)$,
we have 
\[
|\frac{1}{n}\hat{X}^{T}e|_{\infty}\precsim\max\left\{ \psi_{1}k_{1}\sqrt{\frac{\log\max(p,\, d)}{n}},\,\psi_{2}\sqrt{\frac{\log p}{n}}\right\} ,
\]
where 
\begin{eqnarray*}
\psi_{1} & = & \frac{\sigma_{\eta}\max_{j,j^{'}}|\textrm{cov}(x_{1j^{'}}^{*},\,\mathbf{z}_{1j})|_{\infty}|\beta^{*}|_{1}}{\lambda_{\min}(\Sigma_{Z})},\\
\psi_{2} & = & \max\left\{ \sigma_{X^{*}}\sigma_{\eta}|\beta^{*}|_{1},\,\sigma_{X^{*}}\sigma_{\epsilon}\right\} ,
\end{eqnarray*}
with probability at least $1-c_{1}\exp(-c_{2}\log\min(p,\, d))$ for
some universal constants $c_{1}$ and $c_{2}$.\\
\\
\textbf{Proof}. We have 
\begin{eqnarray*}
\frac{1}{n}\hat{X}^{T}e & = & \frac{1}{n}\hat{X}^{T}\left[(X^{*}-\hat{X})\beta^{*}+\boldsymbol{\eta}\beta^{*}+\epsilon\right]\\
 & = & \frac{1}{n}X^{*T}\left[(X^{*}-\hat{X})\beta^{*}+\boldsymbol{\eta}\beta^{*}+\epsilon\right]+\frac{1}{n}(\hat{X}-X^{*})^{T}\left[(X^{*}-\hat{X})\beta^{*}+\boldsymbol{\eta}\beta^{*}+\epsilon\right].
\end{eqnarray*}
Hence, 
\begin{eqnarray}
|\frac{1}{n}\hat{X}^{T}e|_{\infty} & \leq & |\frac{1}{n}X^{*T}(\hat{X}-X^{*})\beta^{*}|_{\infty}+|\frac{1}{n}X^{*T}\boldsymbol{\eta}\beta^{*}|_{\infty}+|\frac{1}{n}X^{*T}\epsilon|_{\infty}\\
 & + & |\frac{1}{n}(\hat{X}-X^{*})^{T}(\hat{X}-X^{*})\beta^{*}|_{\infty}+|\frac{1}{n}(\hat{X}-X^{*})^{T}\boldsymbol{\eta}\beta^{*}|_{\infty}+|\frac{1}{n}(\hat{X}-X^{*})^{T}\epsilon|_{\infty}.\nonumber 
\end{eqnarray}
We need to bound each of the terms on the right-hand-side of the above
inequality. Let us first bound $|\frac{1}{n}X^{*T}(\hat{X}-X^{*})\beta^{*}|_{\infty}$.
We have 
\[
\frac{1}{n}X^{*T}(\hat{X}-X^{*})\beta^{*}=\left[\begin{array}{c}
\sum_{j=1}^{p}\beta_{j}^{*}\frac{1}{n}\sum_{i=1}^{n}x_{i1}^{*}(\hat{x}_{ij}-x_{ij}^{*})\\
\vdots\\
\sum_{j=1}^{p}\beta_{j}^{*}\frac{1}{n}\sum_{i=1}^{n}x_{ip}^{*}(\hat{x}_{ij}-x_{ij}^{*})
\end{array}\right].
\]
For any $j^{'}=1,...,p$, we have 
\begin{eqnarray*}
|\sum_{j=1}^{p}\beta_{j}^{*}\frac{1}{n}\sum_{i=1}^{n}x_{ij^{'}}^{*}(\hat{x}_{ij}-x_{ij}^{*})| & \leq & \max_{j^{'},\, j}|\frac{1}{n}\sum_{i=1}^{n}x_{ij^{'}}^{*}(\hat{x}_{ij}-x_{ij}^{*})||\beta^{*}|_{1}\\
 & = & \left|\frac{X^{*T}(\hat{X}-X^{*})}{n}\right|_{\infty}|\beta^{*}|_{1}.
\end{eqnarray*}
In proving Lemma 6.1, under the condition $\frac{\log\max(p,\, d)}{n}=o(1)$,
we have, 
\[
\left|\frac{X^{*T}(\hat{X}-X^{*})}{n}\right|_{\infty}\leq c\frac{\sigma_{\eta}\max_{j^{'},\, j}|\mathbb{E}(x_{ij^{'}}^{*},\,\mathbf{z}_{ij})|_{\infty}}{\lambda_{\min}(\Sigma_{Z})}k_{1}\sqrt{\frac{\log\max(p,\, d)}{n}},
\]
with probability at least $1-c_{1}\exp(-c_{2}\log\max(p,\, d))$.
Therefore, 
\[
|\frac{1}{n}X^{*T}(\hat{X}-X^{*})\beta^{*}|_{\infty}\leq c\frac{\sigma_{\eta}\max_{j^{'},\, j}|\mathbb{E}(x_{ij^{'}}^{*},\,\mathbf{z}_{ij})|_{\infty}}{\lambda_{\min}(\Sigma_{Z})}|\beta^{*}|_{1}k_{1}\sqrt{\frac{\log\max(p,\, d)}{n}}.
\]
The term $|\frac{1}{n}(\hat{X}-X^{*})^{T}(\hat{X}-X^{*})\beta^{*}|_{\infty}$
can be bounded using a similar argument and we have, 
\[
|\frac{1}{n}(\hat{X}-X^{*})^{T}(\hat{X}-X^{*})\beta^{*}|_{\infty}\leq c\frac{\sigma_{\eta}^{2}\max_{j^{'},\, j}|\mathbb{E}(\mathbf{z}_{ij^{'}},\,\mathbf{z}_{ij})|_{\infty}}{\lambda_{\min}^{2}(\Sigma_{Z})}|\beta^{*}|_{1}k_{1}^{2}\frac{\log\max(p,\, d)}{n},
\]
with probability at least $1-c_{1}\exp(-c_{2}\log\max(p,\, d))$.
For the term $|\frac{1}{n}X^{*T}\boldsymbol{\eta}\beta^{*}|_{\infty}$,
we have 
\begin{eqnarray*}
|\frac{1}{n}X^{*T}\boldsymbol{\eta}\beta^{*}|_{\infty} & \leq & \max_{j^{'},\, j}|\frac{1}{n}\sum_{i=1}^{n}x_{ij^{'}}^{*}\eta_{ij}||\beta^{*}|_{1}\\
 & \leq & c\sigma_{X^{*}}\sigma_{\eta}|\beta^{*}|_{1}\sqrt{\frac{\log p}{n}},
\end{eqnarray*}
with probability at least $1-c_{1}\exp(-c_{2}\log p)$. The last inequality
follows from Lemma 6.8 and Assumption 1.1 that $\mathbb{E}(\mathbf{z}_{ij^{'}}\eta_{ij})=\mathbf{0}$
for all $j^{'},\, j$ as well as Assumption 3.2 that $\eta_{j}$ is
an \textit{i.i.d.} zero-mean sub-Gaussian vector with parameter $\sigma_{\eta}^{2}$
for $j=1,...,p$, and the random matrix\textbf{ $Z_{j}\in\mathbb{R}^{n\times d_{j}}$
}is\textbf{ }sub-Gaussian with parameters $(\Sigma_{Z_{j}},\,\sigma_{Z}^{2})$
for $j=1,...,p$. For the term $|\frac{1}{n}(X^{*}-\hat{X})^{T}\boldsymbol{\eta}\beta^{*}|_{\infty}$,
we have, 
\begin{eqnarray*}
|\frac{1}{n}(X^{*}-\hat{X})^{T}\boldsymbol{\eta}\beta^{*}|_{\infty} & \leq & \max_{j^{'}}|\hat{\pi}_{j^{'}}-\pi_{j^{'}}^{*}|_{1}\max_{j^{'},\, j}|\frac{1}{n}\sum_{i=1}^{n}\mathbf{z}_{ij^{'}}^{T}\eta_{ij}|_{\infty}|\beta^{*}|_{1}\\
 & \leq & c\frac{\sigma_{Z}\sigma_{\eta}^{2}}{\lambda_{\min}(\Sigma_{Z})}|\beta^{*}|_{1}k_{1}\frac{\log\max(p,\, d)}{n},
\end{eqnarray*}
with probability at least $1-c_{1}\exp(-c_{2}\log\max(p,\, d))$.
Again, the last inequality follows from Lemma 6.8 and Assumption 1.1
that $\mathbb{E}(\mathbf{z}_{ij^{'}}\eta_{ij})=\mathbf{0}$ for all
$j^{'},\, j$ as well as Assumption 3.2. 

To bound the term $|\frac{1}{n}X^{*T}\epsilon|_{\infty}$, note under
Assumptions 3.2 and 3.3 as well as Assumption 1.1 that $\mathbb{E}(\mathbf{z}_{ij}\epsilon_{i})=\mathbf{0}$
for all $j=1,...,p$, again by Lemma 6.8, 
\[
|\frac{1}{n}X^{*T}\epsilon|_{\infty}\leq c\sigma_{X^{*}}\sigma_{\epsilon}\sqrt{\frac{\log p}{n}},
\]
with probability at least $1-c_{1}\exp(-c_{2}\log p)$. 

For the term $|\frac{1}{n}(X^{*}-\hat{X})^{T}\epsilon|_{\infty}$,
we have 
\begin{eqnarray*}
|\frac{1}{n}(X^{*}-\hat{X})^{T}\epsilon|_{\infty} & \leq & \max_{j}|\hat{\pi}_{j}-\pi_{j}^{*}|_{1}\max_{j}|\frac{1}{n}\sum_{i=1}^{n}\mathbf{z}_{ij}^{T}\epsilon_{i}|_{\infty}\\
 & \leq & c\frac{\sigma_{Z}\sigma_{\epsilon}\sigma_{\eta}}{\lambda_{\min}(\Sigma_{Z})}k_{1}\frac{\log\max(p,\, d)}{n},
\end{eqnarray*}
with probability at least $1-c_{1}\exp(-c_{2}\log\max(p,\, d))$. 

Putting everything together, under the condition $\frac{\max\{k_{1}^{2}\log d,\, k_{1}^{2}\log p\}}{n}=o(1)$,
the claim in Lemma 6.2 follows. $\square$

Under the conditions 
\begin{eqnarray*}
\frac{k_{1}^{2}k_{2}^{2}\log\max(p,\, d)}{n} & = & O(1),\\
\frac{k_{1}^{2}\log\max(p,\, d)}{n} & = & o(1),
\end{eqnarray*}
and $\lambda_{n}\asymp k_{2}k_{1}\sqrt{\frac{\log\max(p,\, d)}{n}}$
(the $k_{2}$ factor in the choice of $\lambda_{n}$ comes from the
simple inequality $|\beta^{*}|_{1}\leq k_{2}\max_{j=1,...,p}\beta^{*}$
by exploring the sparsity of $\beta^{*}$), combining Lemmas 3.1,
6.1, and 6.2, we have 
\[
|\hat{\beta}_{H2SLS}-\beta^{*}|_{2}\precsim\max\{\varphi_{1}\sqrt{k_{1}k_{2}}\sqrt{\frac{k_{1}\log\max(d,\, p)}{n}},\,\varphi_{2}\sqrt{\frac{k_{2}\log p}{n}}\},
\]
where 
\begin{eqnarray*}
\varphi_{1} & = & \frac{\sigma_{\eta}\max_{j,j^{'}}|\textrm{cov}(x_{1j^{'}}^{*},\,\mathbf{z}_{1j})|_{\infty}|\beta^{*}|_{1}}{\lambda_{\min}(\Sigma_{Z})\lambda_{\min}(\Sigma_{X^{*}})},\\
\varphi_{2} & = & \max\left\{ \frac{\sigma_{X^{*}}\sigma_{\eta}|\beta^{*}|_{1}}{\lambda_{\min}(\Sigma_{X^{*}})},\,\frac{\sigma_{X^{*}}\sigma_{\epsilon}}{\lambda_{\min}(\Sigma_{X^{*}})}\right\} ,
\end{eqnarray*}
with probability at least $1-c_{1}\exp(-c_{2}\log\min(p,\, d))$ for
some universal positive constants $c_{1}$ and $c_{2}$, which proves
Theorem 3.2. $\square$

\subsection{Theorem 3.3}

Again, we verify the conditions: i) $\hat{\Gamma}=\hat{X}^{T}\hat{X}$
satisfies the RE condition (3) with $\gamma=3$, and (ii) the term
$|\frac{1}{n}\hat{X}^{T}\hat{\epsilon}|_{\infty}\precsim f(k_{1},\, k_{2},\, d,\, p,\, n)$
with high probability. This is done via Lemmas 6.3 and 6.4.\\
\textbf{}\\
\textbf{Lemma 6.3} (RE condition): Let $r\in[0,\,1]$. Under Assumptions
1.1, 3.1, 3.3, 3.4, 3.5b, 3.6, and the condition $n\succsim k_{1}^{3-2r}\log\max(p,\, d)$,
we have, for some universal constants $c$, $c^{'}$, $c_{1}$, and
$c_{2}$, 
\[
\frac{|\hat{X}v^{0}|_{2}^{2}}{n}\geq\left(\kappa_{1}-c\kappa_{2}k_{1}^{3/2-r}\sqrt{\frac{\log\max(p,\, d)}{n}}\right)|v^{0}|_{2}^{2}-c^{'}\kappa_{3}\frac{k_{1}^{r}\log\max(p,\, d)}{n}|v^{0}|_{1}^{2},\quad\textrm{for all }v^{0}\in\mathbb{R}^{p},
\]
with probability at least $1-c_{1}\exp(-c_{2}\log\max(p,\, d))$,
where 
\begin{eqnarray*}
\kappa_{1} & = & \frac{\lambda_{\min}(\Sigma_{X^{*}})}{2},\quad\kappa_{2}=\max(b_{2}b_{1}^{-1},\, b_{3}b_{1}^{-2}),\quad\kappa_{3}=\max\left\{ b_{0},\, b_{2}b_{1}^{-1},\, b_{3}b_{1}^{-2}\right\} ,\\
b_{0} & = & \lambda_{\min}(\Sigma_{X^{*}})\max\left\{ \frac{\sigma_{X^{*}}^{4}}{\lambda_{\min}^{2}(\Sigma_{X^{*}})},\,1\right\} ,\\
b_{1} & = & \frac{\lambda_{\min}(\Sigma_{Z})}{\sigma_{\eta}},\quad b_{2}=\max\left\{ \sigma_{X^{*}}\sigma_{W},\;\sup_{v\in\mathbb{K}(2s,\, p)\times\mathbb{K}(k_{1},\, d_{1})\times...\times\mathbb{K}(k_{1},\, d_{p})}\left|v^{0T}\left[\mathbb{E}(x_{1j^{'}}^{*}\mathbf{z}_{1j}v^{j})\right]v^{0}\right|\right\} ,\\
b_{3} & = & \max\left\{ \sigma_{W}^{2},\;\sup_{v\in\mathbb{K}(2s,\, p)\times\mathbb{K}^{2}(k_{1},\, d_{1})\times...\times\mathbb{K}^{2}(k_{1},\, d_{p})}\left|v^{0T}\left[\mathbb{E}(v^{j^{'}}\mathbf{z}_{1j^{'}}^{T}\mathbf{z}_{1j}v^{j})\right]v^{0}\right|\right\} .
\end{eqnarray*}
\textbf{}\\
\textbf{Proof}. Again, 
\begin{eqnarray*}
\left|v^{0T}\frac{\hat{X}^{T}\hat{X}}{n}v^{0}\right| & \geq & \left|v^{0T}\frac{X^{*T}X^{*}}{n}v^{0}\right|-\left|v^{0T}\left(\frac{X^{*T}X^{*}-\hat{X}^{T}\hat{X}}{n}\right)v^{0}\right|\\
 & \geq & \left|v^{0T}\frac{X^{*T}X^{*}}{n}v^{0}\right|-\left(\left|v^{0T}\frac{X^{*T}(\hat{X}-X^{*})}{n}v^{0}\right|+\left|v^{0T}\frac{(\hat{X}-X^{*})^{T}\hat{X}}{n}v^{0}\right|\right)\\
 & \geq & \left|v^{0T}\frac{X^{*T}X^{*}}{n}v^{0}\right|-\left|v^{0T}\frac{X^{*T}(\hat{X}-X^{*})}{n}v^{0}\right|\\
 &  & -\left|v^{0T}\frac{(\hat{X}-X^{*})^{T}X^{*}}{n}v^{0}\right|-\left|v^{0T}\frac{(\hat{X}-X^{*})^{T}(\hat{X}-X^{*})}{n}v^{0}\right|.
\end{eqnarray*}
To bound the above terms, I apply a discretization argument motivated
by the idea in Loh and Wainwright (2012). This type of argument is
often used in statistical problems requiring manipulating and controlling
collections of random variables indexed by sets with an infinite number
of elements. For the particular problem in this paper, I work with
the product space $\mathbb{K}(2s,\, p)\times\mathbb{K}(k_{1},\, d_{1})\times...\times\mathbb{K}(k_{1},\, d_{p})$
and $\mathbb{K}(2s,\, p)\times\mathbb{K}^{2}(k_{1},\, d_{1})\times...\times\mathbb{K}^{2}(k_{1},\, d_{p})$.
For $s\geq1$ and $L\geq1$, recall the notation $\mathbb{K}(s,\, L):=\{v\in\mathbb{R}^{L}\,|\,|v|_{2}\leq1,\,|v|_{0}\leq s\}$.
Given $V^{j}\subseteq\{1,...,d_{j}\}$ and $V^{0}\subseteq\{1,...,p\}$,
define $S_{V^{j}}=\{v\in\mathbb{R}^{d_{j}}:\,|v|_{2}\leq1,\, J(v)\subseteq V^{j}\}$
and $S_{V^{0}}=\{v\in\mathbb{R}^{p}:\,|v|_{2}\leq1,\, J(v)\subseteq V^{0}\}$.
Note that $\mathbb{K}(k_{1},\, d_{j})=\cup_{|V^{j}|\leq k_{1}}S_{V^{j}}$
and $\mathbb{K}(2s,\, p)=\cup_{|V^{0}|\leq2s}S_{V^{0}}$ with $s=s(r):=\frac{1}{c}\frac{n}{k_{1}^{r}\log\max(p,\, d)}\min\left\{ \frac{\lambda_{\min}^{2}(\Sigma_{X^{*}})}{\sigma_{X^{*}}^{4}},\,1\right\} ,\; r\in[0,\,1]$.
The choice of $s$ is explained in the proof for Lemma 6.10. If $\mathcal{V}^{j}=\{t_{1}^{j},...,t_{m_{j}}^{j}\}$
is a $\frac{1}{9}$-cover of $S_{V^{j}}$ ($\mathcal{V}^{0}=\{t_{1}^{0},...,t_{m_{0}}^{0}\}$
is a $\frac{1}{9}$-cover of $S_{V^{0}}$), for every $v^{j}\in S_{V^{j}}$
($v^{0}\in S_{V^{0}}$), we can find some $t_{i}^{j}\in\mathcal{V}^{j}$
($t_{i^{'}}^{0}\in\mathcal{V}^{0}$) such that $|\triangle v^{j}|_{2}\leq\frac{1}{9}$
($|\triangle v^{0}|_{2}\leq\frac{1}{9}$), where $\triangle v^{j}=v^{j}-t_{i}^{j}$
(respectively, $\triangle v^{0}=v^{0}-t_{i^{'}}^{0}$). By Ledoux
and Talagrand (1991), we can construct $\mathcal{V}^{j}$ with $|\mathcal{V}^{j}|\leq81^{k_{1}}$
and $|\mathcal{V}^{0}|\leq81^{2s}$. Therefore, for $v^{0}\in\mathbb{K}(2s,\, p)$,
there is some $S_{V^{0}}$ and $t_{i^{'}}^{0}\in\mathcal{V}^{0}$
such that 
\begin{eqnarray*}
v^{0T}\frac{X^{*T}(\hat{X}-X^{*})}{n}v^{0} & = & (t_{i^{'}}^{0}+v^{0}-t_{i^{'}}^{0})^{T}\frac{X^{*T}(\hat{X}-X^{*})}{n}(t_{i^{'}}^{0}+v^{0}-t_{i^{'}}^{0})\\
 & = & t_{i^{'}}^{0T}\frac{X^{*T}(\hat{X}-X^{*})}{n}t_{i^{'}}^{0}+2\triangle v^{0T}\frac{X^{*T}(\hat{X}-X^{*})}{n}t_{i^{'}}^{0}+\triangle v^{0T}\frac{X^{*T}(\hat{X}-X^{*})}{n}\triangle v^{0}
\end{eqnarray*}
with $|\triangle v^{0}|_{2}\leq\frac{1}{9}$. 

Recall for the $(j^{'},\, j)$ element of the matrix\textbf{ }$\frac{X^{*T}(\hat{X}-X^{*})}{n}$,
we have 
\[
\frac{1}{n}\mathbf{x}_{j^{'}}^{*T}(\hat{\mathbf{x}}_{j}-\mathbf{x}_{j}^{*})=\left(\frac{1}{n}\sum_{i=1}^{n}x_{ij^{'}}^{*}\mathbf{z}_{ij}\right)(\hat{\pi}_{j}-\pi_{j}^{*}).
\]
Let $\frac{\lambda_{\min}(\Sigma_{Z})}{c\sigma_{\eta}}=b_{1}$. Notice
that, under Assumptions 3.5b and 3.6, $|\hat{\pi}_{j}-\pi_{j}^{*}|_{2}b_{1}\sqrt{\frac{n}{k_{1}\log\max(p,\, d)}}\leq1$
and $|\textrm{supp}(\hat{\pi}_{j}-\pi_{j}^{*})|\leq k_{1}$ for every
$j=1,...,p$. Define $\bar{\pi}_{j}=(\hat{\pi}_{j}-\pi_{j}^{*})b_{1}\sqrt{\frac{n}{k_{1}\log\max(p,\, d)}}$
and hence, $\bar{\pi}_{j}\in\mathbb{K}(k_{1},\, d_{j})=\cup_{|V^{j}|\leq k_{1}}S_{V^{j}}$.
Therefore, there is some $S_{V^{j}}$ with $|V^{j}|\leq k_{1}$ and
$t_{i}^{j}\in\mathcal{V}^{j}$ (where $\mathcal{V}^{j}=\{t_{1}^{j},...,t_{m_{j}}^{j}\}$
is a $\frac{1}{9}$-cover of $S_{V^{j}}$) such that 
\begin{eqnarray*}
\frac{1}{n}x_{j^{'}}^{*T}\mathbf{z}_{j}(\hat{\pi}_{j}-\pi_{j}^{*}) & = & \frac{1}{n}x_{j^{'}}^{*T}\mathbf{z}_{j}(t_{i}^{j}+\bar{\pi}_{j}-t_{i}^{j})b_{1}^{-1}\sqrt{\frac{k_{1}\log\max(p,\, d)}{n}}\\
 & = & b_{1}^{-1}\sqrt{\frac{k_{1}\log\max(p,\, d)}{n}}\left(\frac{1}{n}x_{j^{'}}^{*T}\mathbf{z}_{j}t_{i}^{j}+\frac{1}{n}x_{j^{'}}^{*T}\mathbf{z}_{j}\triangle v^{j}\right)
\end{eqnarray*}
with $|\triangle v^{j}|_{2}\leq\frac{1}{9}$. 

Denote a matrix $A$ by $\left[A_{j^{'}j}\right]$, where the $(j^{'},\, j)$
element of $A$ is $A_{j^{'}j}$. Define $v=(v^{0},\, v^{1},\,...,\, v^{p})\in S_{V}:=S_{V^{0}}\times S_{V^{1}}\times...\times S_{V^{p}}$.
Hence, 
\[
\left|v^{0T}\frac{X^{*T}(\hat{X}-X^{*})}{n}v^{0}-\mathbb{E}(v^{0T}\frac{X^{*T}(\hat{X}-X^{*})}{n}v^{0})\right|
\]
\begin{eqnarray*}
 & \leq & \sup_{v\in S_{V}}b_{1}^{-1}\sqrt{\frac{k_{1}\log\max(p,\, d)}{n}}\left|v^{0T}\left[\frac{1}{n}x_{j^{'}}^{*T}\mathbf{z}_{j}v^{j}-\mathbb{E}(x_{1j^{'}}^{*}\mathbf{z}_{1j}v^{j})\right]v^{0}\right|\\
 & \leq & b_{1}^{-1}\sqrt{\frac{k_{1}\log\max(p,\, d)}{n}}\{\max_{i^{'},\, i}\left|t_{i^{'}}^{0T}\left[\frac{1}{n}x_{j^{'}}^{*T}\mathbf{z}_{j}t_{i}^{j}-\mathbb{E}(x_{1j^{'}}^{*}\mathbf{z}_{1j}t_{i}^{j})\right]t_{i^{'}}^{0}\right|\\
 &  & +\sup_{v\in S_{V}}\left|t_{i^{'}}^{0T}\left[\frac{1}{n}x_{j^{'}}^{*T}\mathbf{z}_{j}\triangle v^{j}-\mathbb{E}(x_{1j^{'}}^{*}\mathbf{z}_{1j}\triangle v^{j})\right]t_{i^{'}}^{0}\right|+\sup_{v\in S_{V}}2\left|\triangle v^{0T}\left[\frac{1}{n}x_{j^{'}}^{*T}\mathbf{z}_{j}t_{i}^{j}-\mathbb{E}(x_{1j^{'}}^{*}\mathbf{z}_{1j}t_{i}^{j})\right]t_{i^{'}}^{0}\right|\\
 &  & +\sup_{v\in S_{V}}2\left|\triangle v^{0T}\left[\frac{1}{n}x_{j^{'}}^{*T}\mathbf{z}_{j}\triangle v^{j}-\mathbb{E}(x_{1j^{'}}^{*}\mathbf{z}_{1j}\triangle v^{j})\right]t_{i^{'}}^{0}\right|+\sup_{v\in S_{V}}\left|\triangle v^{0T}\left[\frac{1}{n}x_{j^{'}}^{*T}\mathbf{z}_{j}t_{i}^{j}-\mathbb{E}(x_{1j^{'}}^{*}\mathbf{z}_{1j}t_{i}^{j})\right]\triangle v^{0}\right|\\
 &  & +\sup_{v\in S_{V}}\left|\triangle v^{0T}\left[\frac{1}{n}x_{j^{'}}^{*T}\mathbf{z}_{j}\triangle v^{j}-\mathbb{E}(x_{1j^{'}}^{*}\mathbf{z}_{1j}\triangle v^{j})\right]\triangle v^{0}\right|\}\\
 & \leq & b_{1}^{-1}\sqrt{\frac{k_{1}\log\max(p,\, d)}{n}}\{\max_{i^{'},\, i}\left|t_{i^{'}}^{0T}\left[\frac{1}{n}x_{j^{'}}^{*T}\mathbf{z}_{j}t_{i}^{j}-\mathbb{E}(x_{1j^{'}}^{*}\mathbf{z}_{1j}t_{i}^{j})\right]t_{i^{'}}^{0}\right|\\
 &  & +\sup_{v\in S_{V}}\frac{1}{9}\left|v^{0T}\left[\frac{1}{n}x_{j^{'}}^{*T}\mathbf{z}_{j}v^{j}-\mathbb{E}(x_{1j^{'}}^{*}\mathbf{z}_{1j}v^{j})\right]v^{0}\right|+\sup_{v\in S_{V}}\frac{2}{9}\left|v^{0T}\left[\frac{1}{n}x_{j^{'}}^{*T}\mathbf{z}_{j}v^{j}-\mathbb{E}(x_{1j^{'}}^{*}\mathbf{z}_{1j}v^{j})\right]v^{0}\right|\\
 &  & +\sup_{v\in S_{V}}\frac{2}{81}\left|v^{0T}\left[\frac{1}{n}x_{j^{'}}^{*T}\mathbf{z}_{j}v^{j}-\mathbb{E}(x_{1j^{'}}^{*}\mathbf{z}_{1j}v^{j})\right]v^{0}\right|+\sup_{v\in S_{V}}\frac{1}{81}\left|v^{0T}\left[\frac{1}{n}x_{j^{'}}^{*T}\mathbf{z}_{j}v^{j}-\mathbb{E}(x_{1j^{'}}^{*}\mathbf{z}_{1j}v^{j})\right]v^{0}\right|\\
 &  & +\sup_{v\in S_{V}}\frac{1}{729}\left|v^{0T}\left[\frac{1}{n}x_{j^{'}}^{*T}\mathbf{z}_{j}v^{j}-\mathbb{E}(x_{1j^{'}}^{*}\mathbf{z}_{1j}v^{j})\right]v^{0}\right|\},
\end{eqnarray*}
where the last inequality uses the fact that $9\triangle v^{j}\in S_{V^{j}}$
and $9\triangle v^{0}\in S_{V^{0}}$. Therefore, 
\[
\sup_{v\in S_{V}}b_{1}^{-1}\sqrt{\frac{k_{1}\log\max(p,\, d)}{n}}\left|v^{0T}\left[\frac{1}{n}x_{j^{'}}^{*T}\mathbf{z}_{j}v^{j}-\mathbb{E}(x_{1j^{'}}^{*}\mathbf{z}_{1j}v^{j})\right]v^{0}\right|
\]
\begin{eqnarray*}
 & \leq & \frac{729}{458}b_{1}^{-1}\sqrt{\frac{k_{1}\log\max(p,\, d)}{n}}\max_{i^{'},\, i}t_{i^{'}}^{0T}\left[\frac{1}{n}x_{j^{'}}^{*T}\mathbf{z}_{j}t_{i}^{j}-\mathbb{E}(x_{1j^{'}}^{*}\mathbf{z}_{1j}t_{i}^{j})\right]t_{i^{'}}^{0}\\
 & \leq & 2b_{1}^{-1}\sqrt{\frac{k_{1}\log\max(p,\, d)}{n}}\max_{i^{'},\, i}t_{i^{'}}^{0T}\left[\frac{1}{n}x_{j^{'}}^{*T}\mathbf{z}_{j}t_{i}^{j}-\mathbb{E}(x_{1j^{'}}^{*}\mathbf{z}_{1j}t_{i}^{j})\right]t_{i^{'}}^{0}.
\end{eqnarray*}
Under Assumptions 3.3 and 3.4, $x_{j^{'}}^{*}$ is a sub-Gaussian
vector with parameter at most $\sigma_{X^{*}}$ for every $j^{'}=1,...,p$,
and $Z_{j}t_{i}^{j}:=\mathbf{w}_{j}$ is a sub-Gaussian vector with
parameter at most $\sigma_{W^{*}}$. An application of Lemma 6.8 and
a union bound yields 
\[
\mathbb{P}\left(\sup_{v\in S_{V}}\left|v^{0T}\left[\frac{1}{n}x_{j^{'}}^{*T}\mathbf{z}_{j}v^{j}\right]v^{0}-v^{0T}\left[\mathbb{E}(x_{1j^{'}}^{*}\mathbf{z}_{1j}v^{j})\right]v^{0}\right|\geq t\right)\leq81^{2sk_{1}}81^{2s}2\exp(-cn\min(\frac{t^{2}}{\sigma_{X^{*}}^{2}\sigma_{W}^{2}},\,\frac{t}{\sigma_{X^{*}}\sigma_{W}})),
\]
where the exponent $2sk_{1}$ in $81^{2sk_{1}}$ uses the fact that
there are at most $2s$ non-zero components in $v^{0}\in S_{V^{0}}$
and hence only $2s$ out of $p$ entries of $v^{1},...,v^{p}$ will
be multiplied by a non-zero scalar, which leads to a reduction of
dimensions. A second application of a union bound over the $\left(\begin{array}{c}
d_{j}\\
\left\lfloor k_{1}\right\rfloor 
\end{array}\right)\leq d^{k_{1}}$ choices of $V^{j}$ and respectively, the $\left(\begin{array}{c}
p\\
\left\lfloor 2s\right\rfloor 
\end{array}\right)\leq p^{2s}$ choices of $V^{0}$ yields 
\[
\mathbb{P}\left(\sup_{v\in\mathbb{K}(2s,\, p)\times\mathbb{K}(k_{1},\, d_{1})\times...\times\mathbb{K}(k_{1},\, d_{p})}\left|v^{0T}\left[\frac{1}{n}x_{j^{'}}^{*T}\mathbf{z}_{j}v^{j}\right]v^{0}-v^{0T}\left[\mathbb{E}(x_{1j^{'}}^{*}\mathbf{z}_{1j}v^{j})\right]v^{0}\right|\geq t\right)
\]
\begin{eqnarray*}
 & \leq & p^{2s}d^{2sk_{1}}\cdot2\exp(-cn\min(\frac{t^{2}}{\sigma_{X^{*}}^{2}\sigma_{W}^{2}},\,\frac{t}{\sigma_{X^{*}}\sigma_{W}}))\\
 & \leq & 2\exp(-cn\min(\frac{t^{2}}{\sigma_{X^{*}}^{2}\sigma_{W}^{2}},\,\frac{t}{\sigma_{X^{*}}\sigma_{W}})+2sk_{1}\log d+2s\log p).
\end{eqnarray*}
With the choice of $s=s(r):=\frac{1}{c}\frac{n}{k_{1}^{r}\log\max(p,\, d)}\min\left\{ \frac{\lambda_{\min}^{2}(\Sigma_{X^{*}})}{\sigma_{X^{*}}^{4}},\,1\right\} ,\; r\in[0,\,1]$
from the proof for Lemma 6.10 and $t=c^{'}k_{1}^{1-r}\sigma_{X^{*}}\sigma_{W}$
for some universal constant $c^{'}\geq1$, we have 
\[
\left|v^{0T}\frac{X^{*T}(\hat{X}-X^{*})}{n}v^{0}-\mathbb{E}[v^{0T}\frac{X^{*T}(\hat{X}-X^{*})}{n}v^{0}]\right|
\]

\begin{eqnarray*}
 & \leq & \left(\sup_{v\in\mathbb{K}(2s,\, p)\times\mathbb{K}(k_{1},\, d_{1})\times...\times\mathbb{K}(k_{1},\, d_{p})}\left|v^{0T}\left[\frac{1}{n}x_{j^{'}}^{*T}\mathbf{z}_{j}v^{j}-\mathbb{E}(x_{1j^{'}}^{*}\mathbf{z}_{1j}v^{j})\right]v^{0}\right|\right)b_{1}^{-1}\sqrt{\frac{k_{1}\log\max(p,\, d)}{n}}\\
 & \leq & c^{'}b_{1}^{-1}k_{1}^{1-r}\sqrt{\frac{k_{1}\log\max(p,\, d)}{n}}\sigma_{X^{*}}\sigma_{W}
\end{eqnarray*}
with probability at least $1-c_{1}^{'}\exp(-c_{2}^{'}nk_{1}^{1-r})-c_{1}^{''}\exp(-c_{2}^{''}\log\max(p,\, d))=1-c_{1}\exp(-c_{2}\log\max(p,\, d))$
(given $d>n$ and $p>n$ is the regime of our interests). Therefore,
we have 
\begin{eqnarray*}
\left|v^{0T}\frac{X^{*T}(\hat{X}-X^{*})}{n}v^{0}\right| & \leq & \left(\sup_{v\in\mathbb{K}(2s,\, p)\times\mathbb{K}(k_{1},\, d_{1})\times...\times\mathbb{K}(k_{1},\, d_{p})}\left|v^{0T}\left[\mathbb{E}(x_{1j^{'}}^{*}\mathbf{z}_{1j}v^{j})\right]v^{0}\right|\right)b_{1}^{-1}\sqrt{\frac{k_{1}\log\max(p,\, d)}{n}}\\
 & + & c^{'}b_{1}^{-1}k_{1}^{3/2-r}\sqrt{\frac{\log\max(p,\, d)}{n}}\sigma_{X^{*}}\sigma_{W}\\
 & \leq & cb_{2}b_{1}^{-1}k_{1}^{3/2-r}\sqrt{\frac{\log\max(p,\, d)}{n}},
\end{eqnarray*}
where $b_{2}=\max\left\{ \sigma_{X^{*}}\sigma_{W},\;\sup_{v\in\mathbb{K}(2s,\, p)\times\mathbb{K}(k_{1},\, d_{1})\times...\times\mathbb{K}(k_{1},\, d_{p})}\left|v^{0T}\left[\mathbb{E}(x_{1j^{'}}^{*}\mathbf{z}_{1j}v^{j})\right]v^{0}\right|\right\} $.
Notice that the term 
\[
\sup_{v\in\mathbb{K}(2s,\, p)\times\mathbb{K}(k_{1},\, d_{1})\times...\times\mathbb{K}(k_{1},\, d_{p})}\left|v^{0T}\left[\mathbb{E}(x_{1j^{'}}^{*}\mathbf{z}_{1j}v^{j})\right]v^{0}\right|
\]
is bounded above by the spectral norm of the matrix $\left[\mathbb{E}(x_{1j^{'}}^{*}\mathbf{z}_{1j}v^{j})\right]$
for some $v^{1}\times...\times v^{p}\in\mathbb{K}(k_{1},\, d_{1})\times...\times\mathbb{K}(k_{1},\, d_{p})$.

The term $\left|v^{0T}\frac{(\hat{X}-X^{*})^{T}(\hat{X}-X^{*})}{n}v^{0}\right|$
can be bounded using a similar argument. In particular, for the $(j^{'},\, j)$
element of the matrix\textbf{ }$\frac{(\hat{X}-X^{*})^{T}(\hat{X}-X^{*})}{n}$,
we have 
\begin{eqnarray*}
\frac{1}{n}(\hat{\mathbf{x}}_{j^{'}}-\mathbf{x}_{j^{'}}^{*})^{T}(\hat{\mathbf{x}}_{j}-\mathbf{x}_{j}^{*}) & = & (\hat{\pi}_{j^{'}}-\pi_{j^{'}}^{*})^{T}\left(\frac{1}{n}\sum_{i=1}^{n}\mathbf{z}_{ij^{'}}\mathbf{z}_{ij}\right)(\hat{\pi}_{j}-\pi_{j}^{*})\\
 & = & \frac{1}{n}(t_{i^{'}}^{j^{'}}+\bar{\pi}_{j^{'}}-t_{i^{'}}^{j^{'}})^{T}\mathbf{z}_{j^{'}}^{T}\mathbf{z}_{j}(t_{i}^{j}+\bar{\pi}_{j}-t_{i}^{j})b_{1}^{-2}\frac{k_{1}\log\max(p,\, d)}{n}\\
 & = & b_{1}^{-2}\frac{k_{1}\log\max(p,\, d)}{n}\{\frac{1}{n}t_{i^{'}}^{j^{'}T}\mathbf{z}_{j^{'}}^{T}\mathbf{z}_{j}t_{i}^{j}+\frac{1}{n}\triangle v^{j^{'}T}\mathbf{z}_{j^{'}}^{T}\mathbf{z}_{j}t_{i}^{j}\\
 & + & \frac{1}{n}t_{i^{'}}^{j^{'}T}\mathbf{z}_{j^{'}}^{T}\mathbf{z}_{j}\triangle v^{j}+\frac{1}{n}\triangle v^{j^{'}T}\mathbf{z}_{j^{'}}^{T}\mathbf{z}_{j}\triangle v^{j}\}
\end{eqnarray*}
Combining with \textbf{
\begin{eqnarray*}
v^{0T}\frac{(\hat{X}-X^{*})^{T}(\hat{X}-X^{*})}{n}v^{0} & = & t_{i^{''}}^{0T}\frac{(\hat{X}-X^{*})^{T}(\hat{X}-X^{*})}{n}t_{i^{''}}^{0}\\
 & + & 2\triangle v^{0T}\frac{(\hat{X}-X^{*})^{T}(\hat{X}-X^{*})}{n}t_{i^{''}}^{0}+\triangle v^{0T}\frac{(\hat{X}-X^{*})^{T}(\hat{X}-X^{*})}{n}\triangle v^{0},
\end{eqnarray*}
}Define $S_{V}:=S_{V^{0}}\times S_{V^{1}}^{2}\times...\times S_{V^{p}}^{2}$.
After some tedious algebra, we obtain 
\[
\left|v^{0T}\frac{(\hat{X}-X^{*})^{T}(\hat{X}-X^{*})}{n}v^{0}-\mathbb{E}(v^{0T}\frac{(\hat{X}-X^{*})^{T}(\hat{X}-X^{*})}{n}v^{0})\right|
\]
\begin{eqnarray*}
 & \leq & \sup_{v\in S_{V}}b_{1}^{-2}\frac{k_{1}\log\max(p,\, d)}{n}\left|v^{0T}\left[\frac{1}{n}v^{j^{'}}\mathbf{z}_{j^{'}}^{T}\mathbf{z}_{j}v^{j}-v^{j^{'}}\mathbb{E}(\mathbf{z}_{1j^{'}}^{T}\mathbf{z}_{1j})v^{j}\right]v^{0}\right|\\
 & \leq & b_{1}^{-2}\frac{k_{1}\log\max(p,\, d)}{n}\{\max_{i^{''},\, i^{'},\, i}\left|t_{i^{''}}^{0T}\left[\frac{1}{n}t_{i^{'}}^{j^{'}T}\mathbf{z}_{j^{'}}^{T}\mathbf{z}_{j}t_{i}^{j}-\mathbb{E}(t_{i^{'}}^{j^{'}}\mathbf{z}_{1j^{'}}^{T}\mathbf{z}_{1j}t_{i}^{j})\right]t_{i^{''}}^{0}\right|\\
 &  & +\frac{3439}{6561}\sup_{v\in S_{V}}\left|v^{0T}\left[\frac{1}{n}v^{j^{'}T}\mathbf{z}_{j^{'}}^{T}\mathbf{z}_{j}v^{j}-\mathbb{E}(v^{j^{'}}\mathbf{z}_{1j^{'}}^{T}\mathbf{z}_{1j}v^{j})\right]v^{0}\right|\}.
\end{eqnarray*}
Hence,\textbf{ 
\[
\sup_{v\in S_{V}}b_{1}^{-2}\frac{k_{1}\log\max(p,\, d)}{n}\left|v^{0T}\left[\frac{1}{n}v^{j^{'}}\mathbf{z}_{j^{'}}^{T}\mathbf{z}_{j}v^{j}-\mathbb{E}(v^{j^{'}}\mathbf{z}_{1j^{'}}^{T}\mathbf{z}_{1j}v^{j})\right]v^{0}\right|
\]
\begin{eqnarray*}
 & \leq & \frac{6561}{3122}b_{1}^{-2}\frac{k_{1}\log\max(p,\, d)}{n}\max_{i^{''},\, i^{'},\, i}\left|t_{i^{''}}^{0T}\left[\frac{1}{n}t_{i^{'}}^{j^{'}T}\mathbf{z}_{j^{'}}^{T}\mathbf{z}_{j}t_{i}^{j}-\mathbb{E}(t_{i^{'}}^{j^{'}}\mathbf{z}_{1j^{'}}^{T}\mathbf{z}_{1j}t_{i}^{j})\right]t_{i^{''}}^{0}\right|\\
 & \leq & 3b_{1}^{-2}\frac{k_{1}\log\max(p,\, d)}{n}\max_{i^{''},\, i^{'},\, i}\left|t_{i^{''}}^{0T}\left[\frac{1}{n}t_{i^{'}}^{j^{'}T}\mathbf{z}_{j^{'}}^{T}\mathbf{z}_{j}t_{i}^{j}-\mathbb{E}(t_{i^{'}}^{j^{'}}\mathbf{z}_{1j^{'}}^{T}\mathbf{z}_{1j}t_{i}^{j})\right]t_{i^{''}}^{0}\right|.
\end{eqnarray*}
}An application of Lemma 6.8 and a sequence of union bounds yields\textbf{
\[
\mathbb{P}\left(\sup_{v\in\mathbb{K}(2s,\, p)\times\mathbb{K}^{2}(k_{1},\, d_{1})\times...\times\mathbb{K}^{2}(k_{1},\, d_{p})}\left|v^{0T}\left[\frac{1}{n}v^{j^{'}}\mathbf{z}_{j^{'}}^{T}\mathbf{z}_{j}v^{j}\right]v^{0}-v^{0T}\left[\mathbb{E}(v^{j^{'}}\mathbf{z}_{1j^{'}}^{T}\mathbf{z}_{1j}v^{j})\right]v^{0}\right|\geq t\right)
\]
\[
\leq2\exp(-cn\min(\frac{t^{2}}{\sigma_{W}^{4}},\,\frac{t}{\sigma_{W}^{2}})+4sk_{1}\log d+2s\log p).
\]
}Under the choice of $s=s(r):=\frac{1}{c}\frac{n}{k_{1}^{r}\log\max(p,\, d)}\min\left\{ \frac{\lambda_{\min}^{2}(\Sigma_{X^{*}})}{\sigma_{X^{*}}^{4}},\,1\right\} ,\; r\in[0,\,1]$
from the proof for Lemma 6.10 and $t=c^{''}k_{1}^{1-r}\sigma_{W}^{2}$
for some universal constant $c^{''}\geq1$, we have, 
\[
\left|v^{0T}\frac{(\hat{X}-X^{*})^{T}(\hat{X}-X^{*})}{n}v^{0}-\mathbb{E}[v^{0T}\frac{(\hat{X}-X^{*})^{T}(\hat{X}-X^{*})}{n}v^{0}]\right|
\]
\begin{eqnarray*}
 & \leq & \left(\sup_{v\in\mathbb{K}(2s,\, p)\times\mathbb{K}^{2}(k_{1},\, d_{1})\times...\times\mathbb{K}^{2}(k_{1},\, d_{p})}\left|v^{0T}\left[\frac{1}{n}v^{j^{'}}\mathbf{z}_{j^{'}}^{T}\mathbf{z}_{j}v^{j}\right]v^{0}-v^{0T}\left[\mathbb{E}(v^{j^{'}}\mathbf{z}_{1j^{'}}^{T}\mathbf{z}_{1j}v^{j})\right]v^{0}\right|\right)b_{1}^{-2}\frac{k_{1}\log\max(p,\, d)}{n}\\
 & \leq & c^{''}b_{1}^{-2}\frac{k_{1}^{2-r}\log\max(p,\, d)}{n}\sigma_{W}^{2}
\end{eqnarray*}
with probability at least $1-c_{1}^{'}\exp(-c_{2}^{'}nk_{1}^{1-r})-c_{1}^{''}\exp(-c_{2}^{''}\log\max(p,\, d))=1-c_{1}\exp(-c_{2}\log\max(p,\, d))$
(given $d>n$ and $p>n$ is the regime of our interests). Therefore,
we have 
\[
\left|v^{0T}\frac{(\hat{X}-X^{*})^{T}(\hat{X}-X^{*})}{n}v^{0}\right|\leq
\]
\[
\left(\sup_{v\in\mathbb{K}(2s,\, p)\times\mathbb{K}^{2}(k_{1},\, d_{1})\times...\times\mathbb{K}^{2}(k_{1},\, d_{p})}\left|v^{0T}\left[\mathbb{E}(v^{j^{'}}\mathbf{z}_{1j^{'}}^{T}\mathbf{z}_{1j}v^{j})\right]v^{0}\right|\right)b_{1}^{-2}\frac{k_{1}\log\max(p,\, d)}{n}+c^{''}b_{1}^{-2}\frac{k_{1}^{2-r}\log\max(p,\, d)}{n}\sigma_{W}^{2}
\]
\[
\leq cb_{3}b_{1}^{-2}\frac{k_{1}^{2-r}\log\max(p,\, d)}{n},
\]
where $b_{3}=\max\left\{ \sigma_{W}^{2},\;\sup_{v\in\mathbb{K}(2s,\, p)\times\mathbb{K}^{2}(k_{1},\, d_{1})\times...\times\mathbb{K}^{2}(k_{1},\, d_{p})}\left|v^{0T}\left[\mathbb{E}(v^{j^{'}}\mathbf{z}_{1j^{'}}^{T}\mathbf{z}_{1j}v^{j})\right]v^{0}\right|\right\} $.
Notice that the term 
\[
\sup_{v\in\mathbb{K}(2s,\, p)\times\mathbb{K}^{2}(k_{1},\, d_{1})\times...\times\mathbb{K}^{2}(k_{1},\, d_{p})}\left|v^{0T}\left[\mathbb{E}(v^{j^{'}}\mathbf{z}_{1j^{'}}^{T}\mathbf{z}_{1j}v^{j})\right]v^{0}\right|
\]
is bounded above by the spectral norm of the matrix $\left[\mathbb{E}(v^{j^{'}}\mathbf{z}_{1j^{'}}^{T}\mathbf{z}_{1j}v^{j})\right]$
for some $(v^{1}\times...\times v^{p})\times(v^{1}\times...\times v^{p})\in\mathbb{K}^{2}(k_{1},\, d_{1})\times...\times\mathbb{K}^{2}(k_{1},\, d_{p})$.

By Lemma 6.9, the bound 
\[
\left|v^{0T}\frac{X^{*T}(\hat{X}-X^{*})}{n}v^{0}\right|\leq cb_{2}b_{1}^{-1}k_{1}^{3/2-r}\sqrt{\frac{\log\max(p,\, d)}{n}}\qquad\forall v^{0}\in\mathbb{K}(2s,\, p)
\]
implies 
\begin{equation}
\left|v^{0T}\frac{X^{*T}(\hat{X}-X^{*})}{n}v^{0}\right|\leq27cb_{2}b_{1}^{-1}k_{1}^{3/2-r}\sqrt{\frac{\log\max(p,\, d)}{n}}(|v^{0}|_{2}^{2}+\frac{1}{s}|v^{0}|_{1}^{2})\qquad\forall v^{0}\in\mathbb{R}^{p}.
\end{equation}
Similarly, the bound

\textbf{
\[
\left|v^{0T}\frac{(\hat{X}-X^{*})^{T}(\hat{X}-X^{*})}{n}v^{0}\right|\leq c^{''}b_{3}b_{1}^{-2}\frac{k_{1}^{2-r}\log\max(p,\, d)}{n}\qquad\forall v^{0}\in\mathbb{K}(2s,\, p)
\]
}implies 
\begin{equation}
\left|v^{0T}\frac{(\hat{X}-X^{*})^{T}(\hat{X}-X^{*})}{n}v^{0}\right|\leq27c^{''}b_{3}b_{1}^{-2}\frac{k_{1}^{2-r}\log\max(p,\, d)}{n}(|v^{0}|_{2}^{2}+\frac{1}{s}|v^{0}|_{1}^{2})\qquad\forall v^{0}\in\mathbb{R}^{p}.
\end{equation}
Therefore, applying Lemma 6.10 by choosing $s=s(r):=\frac{1}{c}\frac{n}{k_{1}^{r}\log\max(p,\, d)}\min\left\{ \frac{\lambda_{\min}^{2}(\Sigma_{X^{*}})}{\sigma_{X^{*}}^{4}},\,1\right\} ,\; r\in[0,\,1]$,
under the condition $n\succsim k_{1}^{3-2r}\log\max(p,\, d),$ we
have 
\begin{eqnarray*}
\left|v^{0T}\frac{\hat{X}^{T}\hat{X}}{n}v^{0}\right| & \geq & \left|v^{0T}\frac{X^{*T}X^{*}}{n}v^{0}\right|-c\max(b_{2}b_{1}^{-1},\, b_{3}b_{1}^{-2})k_{1}^{3/2-r}\sqrt{\frac{\log\max(p,\, d)}{n}}(|v^{0}|_{2}^{2}+\frac{1}{s}|v^{0}|_{1}^{2})\\
 & \geq & \left(\frac{\lambda_{\min}(\Sigma_{X^{*}})}{2}-c\max(b_{2}b_{1}^{-1},\, b_{3}b_{1}^{-2})k_{1}^{3/2-r}\sqrt{\frac{\log\max(p,\, d)}{n}}\right)\left|v^{0}\right|_{2}^{2}\\
 &  & -c^{'}\max\left\{ \lambda_{\min}(\Sigma_{X^{*}})\max\left\{ \frac{\sigma_{X^{*}}^{4}}{\lambda_{\min}^{2}(\Sigma_{X^{*}})},\,1\right\} ,\, b_{2}b_{1}^{-1},\, b_{3}b_{1}^{-2})\right\} \frac{k_{1}^{r}\log\max(p,d)}{n}\left|v^{0}\right|_{1}^{2},
\end{eqnarray*}
which can be written in the form 
\[
\frac{|\hat{X}v^{0}|_{2}^{2}}{n}\geq\left(\kappa_{1}-c\kappa_{2}k_{1}^{3/2-r}\sqrt{\frac{\log\max(p,\, d)}{n}}\right)|v^{0}|_{2}^{2}-c^{'}\kappa_{3}\frac{k_{1}^{r}\log\max(p,\, d)}{n}|v^{0}|_{1}^{2},\quad\textrm{for all }v^{0}\in\mathbb{R}^{p},
\]
with probability at least $1-c_{1}\exp(-c_{2}\log\max(p,\, d))$,
where $\kappa_{1}$, $\kappa_{2}$, and $\kappa_{2}$ are defined
in the statement of Lemma 6.3.$\square$

Again, recalling in proving Lemma 3.1, upon our choice $\lambda_{n}$,
we have shown\textbf{ 
\[
\hat{v}=\hat{\beta}_{H2SLS}-\beta^{*}\in\mathbb{C}(J(\beta^{*}),\,3),
\]
}and $|\hat{v}^{0}|_{1}^{2}\leq16|\hat{v}_{J(\beta^{*})}^{0}|_{1}^{2}\leq16k_{2}|\hat{v}_{J(\beta^{*})}^{0}|_{2}^{2}$.
Therefore, if we have the scaling 
\[
\frac{\min_{r\in[0,\,1]}\max\left\{ k_{1}^{3-2r}\log d,\, k_{1}^{3-2r}\log p,\, k_{1}^{r}k_{2}\log d,\, k_{1}^{r}k_{2}\log p\right\} }{n}=O(1),
\]
so that 
\[
c\kappa_{2}k_{1}^{3/2-r}\sqrt{\frac{\log\max(p,\, d)}{n}}+c^{'}\kappa_{3}\frac{k_{2}k_{1}^{r}\log\max(p,\, d)}{n}<\kappa_{1},
\]
then, \textbf{
\[
\left|\hat{v}^{0T}\frac{\hat{X}^{T}\hat{X}}{n}\hat{v}^{0}\right|\geq c_{0}\lambda_{\min}(\Sigma_{X^{*}})\left|\hat{v}^{0}\right|_{2}^{2},
\]
}provided $\sigma_{\eta}$, $\sigma_{W}$, $\sigma_{X^{*}}$, $\max_{j^{'},\, j}|\mathbb{E}(x_{ij^{'}}^{*},\,\mathbf{z}_{ij})|_{\infty}$,
and $\max_{j^{'},\, j}|\mathbb{E}(\mathbf{z}_{ij^{'}},\,\mathbf{z}_{ij})|_{\infty}$
are bounded from above while $\lambda_{\min}(\Sigma_{Z})$ and $\lambda_{\min}(\Sigma_{X^{*}})$
are bounded away from $0$. The above inequality implies RE (3). Because
the argument for showing Lemma 6.1 and that it implies RE (3) also
works under the assumptions of Lemma 6.3, we can combine the scaling
$\frac{k_{1}^{2}k_{2}^{2}\log\max(p,\, d)}{n}=O(1)$ from the proof
for Lemma 6.1 with the scaling $\frac{\min_{r\in[0,\,1]}\max\left\{ k_{1}^{3-2r}\log d,\, k_{1}^{3-2r}\log p,\, k_{1}^{r}k_{2}\log d,\, k_{1}^{r}k_{2}\log p\right\} }{n}=O(1)$
from above to obtain a more optimal scaling of the required sample
size 
\[
\frac{1}{n}\min\left\{ k_{1}^{2}k_{2}^{2}\log\max(p,\, d),\,\min_{r\in[0,\,1]}\max\left\{ k_{1}^{3-2r}\log d,\, k_{1}^{3-2r}\log p,\, k_{1}^{r}k_{2}\log d,\, k_{1}^{r}k_{2}\log p\right\} \right\} =O(1).
\]
\textbf{}\\
\textbf{Lemma 6.4} (Upper bound on $|\frac{1}{n}\hat{X}^{T}e|_{\infty}$):
Under Assumptions 1.1, 3.1-3.4, 3.5b, 3.6, and the condition $\frac{\max(k_{1}\log d,\, k_{1}\log p)}{n}=o(1)$,
we have 
\[
|\frac{1}{n}\hat{X}^{T}e|_{\infty}\precsim\max\left\{ \psi_{1}\sqrt{\frac{k_{1}\log\max(p,\, d)}{n}},\,\psi_{2}\sqrt{\frac{\log p}{n}}\right\} ,
\]
with probability at least $1-c_{1}\exp(-c_{2}\log\min(p,\, d))$ for
some universal constants $c_{1}$ and $c_{2}$, where $\psi_{1}$
and $\psi_{2}$ are defined in Lemma 6.2.\\
\\
\textbf{Proof}. Recall (5) from the proof for Lemma 6.2. Let us first
bound $|\frac{1}{n}X^{*T}(\hat{X}-X^{*})\beta^{*}|_{\infty}$. For
any $j^{'}=1,...,p$, we have
\begin{eqnarray*}
|\sum_{j=1}^{p}\beta_{j}^{*}\frac{1}{n}\sum_{i=1}^{n}x_{ij^{'}}^{*}(\hat{x}_{ij}-x_{ij}^{*})| & \leq & \max_{j^{'},\, j}|\frac{1}{n}\sum_{i=1}^{n}x_{ij^{'}}^{*}(\hat{x}_{ij}-x_{ij}^{*})||\beta^{*}|_{1}\\
 & = & \left|\frac{X^{*T}(\hat{X}-X^{*})}{n}\right|_{\infty}|\beta^{*}|_{1}
\end{eqnarray*}
In proving Lemma 6.3, we have shown that with a covering subset argument,
the $(j^{'},\, j)$ element of the matrix\textbf{ }$\frac{X^{*T}(\hat{X}-X^{*})}{n}$
can be rewritten as follows. 
\begin{eqnarray*}
\frac{1}{n}\mathbf{x}_{j^{'}}^{*T}(\hat{\mathbf{x}}_{j}-\mathbf{x}_{j}^{*}) & = & \left(\frac{1}{n}\sum_{i=1}^{n}x_{ij^{'}}^{*}\mathbf{z}_{ij}\right)(\hat{\pi}_{j}-\pi_{j}^{*})\\
 & = & b_{1}^{-1}\sqrt{\frac{k_{1}\log\max(p,\, d)}{n}}\left(\frac{1}{n}x_{j^{'}}^{*T}\mathbf{z}_{j}t_{i}^{j}+\frac{1}{n}x_{j^{'}}^{*T}\mathbf{z}_{j}\triangle v^{j}\right).
\end{eqnarray*}
Hence, 
\[
\left|\frac{1}{n}\mathbf{x}_{j^{'}}^{*T}(\hat{\mathbf{x}}_{j}-\mathbf{x}_{j}^{*})-\mathbb{E}(\frac{1}{n}\mathbf{x}_{j^{'}}^{*T}(\hat{\mathbf{x}}_{j}-\mathbf{x}_{j}^{*}))\right|
\]
\begin{eqnarray*}
 & \leq & \sup_{v^{j}\in S_{V^{j}}}b_{1}^{-1}\sqrt{\frac{k_{1}\log\max(p,\, d)}{n}}\left|\frac{1}{n}x_{j^{'}}^{*T}\mathbf{z}_{j}v^{j}-\mathbb{E}(x_{1j^{'}}^{*}\mathbf{z}_{1j}v^{j})\right|\\
 & \leq & b_{1}^{-1}\sqrt{\frac{k_{1}\log\max(p,\, d)}{n}}\{\max_{i}\left|\frac{1}{n}x_{j^{'}}^{*T}\mathbf{z}_{j}t_{i}^{j}-\mathbb{E}(x_{1j^{'}}^{*}\mathbf{z}_{1j}t_{i}^{j})\right|+\sup_{v^{j}\in S_{V^{j}}}\left|\frac{1}{n}x_{j^{'}}^{*T}\mathbf{z}_{j}\triangle v^{j}-\mathbb{E}(x_{1j^{'}}^{*}\mathbf{z}_{1j}\triangle v^{j})\right|\}\\
 & \leq & b_{1}^{-1}\sqrt{\frac{k_{1}\log\max(p,\, d)}{n}}\{\max_{i}\left|\frac{1}{n}x_{j^{'}}^{*T}\mathbf{z}_{j}t_{i}^{j}-\mathbb{E}(x_{1j^{'}}^{*}\mathbf{z}_{1j}t_{i}^{j})\right|+\sup_{v^{j}\in S_{V^{j}}}\frac{1}{9}\left|\frac{1}{n}x_{j^{'}}^{*T}\mathbf{z}_{j}v^{j}-\mathbb{E}(x_{1j^{'}}^{*}\mathbf{z}_{1j}v^{j})\right|\}\\
 & \leq & \frac{9}{8}b_{1}^{-1}\sqrt{\frac{k_{1}\log\max(p,\, d)}{n}}\max_{i}\left|\frac{1}{n}x_{j^{'}}^{*T}\mathbf{z}_{j}t_{i}^{j}-\mathbb{E}(x_{1j^{'}}^{*}\mathbf{z}_{1j}t_{i}^{j})\right|.
\end{eqnarray*}
With a similar argument as in the proof for Lemma 6.3, we obtain 
\[
\mathbb{P}\left(\max_{j^{'},\, j}\sup_{v^{j}\in\mathbb{K}(k_{1},\, d_{j})}\left|\frac{1}{n}x_{j^{'}}^{*T}\mathbf{z}_{j}v^{j}-\mathbb{E}(x_{1j^{'}}^{*}\mathbf{z}_{1j}v^{j})\right|\geq t\right)
\]
\begin{eqnarray*}
 & \leq & p^{2}d^{k_{1}}\cdot2\exp(-cn\min(\frac{t^{2}}{\sigma_{X^{*}}^{2}\sigma_{W}^{2}},\,\frac{t}{\sigma_{X^{*}}\sigma_{W}}))\\
 & = & 2\exp(-cn\min(\frac{t^{2}}{\sigma_{X^{*}}^{2}\sigma_{W}^{2}},\,\frac{t}{\sigma_{X^{*}}\sigma_{W}})+k_{1}\log d+2\log p).
\end{eqnarray*}
Consequently, under the condition $\frac{\max(k_{1}\log d,\,\log p)}{n}=o(1)$,
we have 
\[
\left|\frac{X^{*T}(\hat{X}-X^{*})}{n}-\mathbb{E}[\frac{X^{*T}(\hat{X}-X^{*})}{n}]\right|_{\infty}
\]
\begin{eqnarray*}
 & \leq & \left(\max_{j^{'},\, j}\sup_{v^{j}\in\mathbb{K}(k_{1},\, d_{j})}\left|\frac{1}{n}x_{j^{'}}^{*T}\mathbf{z}_{j}v^{j}-\mathbb{E}(x_{1j^{'}}^{*}\mathbf{z}_{1j}v^{j})\right|\right)b_{1}^{-1}\sqrt{\frac{k_{1}\log\max(p,\, d)}{n}}\\
 & \leq & c^{'}b_{1}^{-1}\sqrt{\frac{k_{1}\log\max(p,\, d)}{n}}\sigma_{X^{*}}\sigma_{W}\max\left\{ \sqrt{\frac{k_{1}\log d}{n}},\,\sqrt{\frac{\log p}{n}}\right\} ,
\end{eqnarray*}
with probability at least $1-c_{1}\exp(-c_{2}\log\max(p,\, d))$.
This implies, 
\begin{eqnarray}
\left|\frac{X^{*T}(\hat{X}-X^{*})}{n}\right|_{\infty} & \leq & \left(\max_{j^{'},\, j}\sup_{v^{j}\in\mathbb{K}(k_{1},\, d_{j})}\left|\mathbb{E}(x_{1j^{'}}^{*}\mathbf{z}_{1j}v^{j})\right|\right)b_{1}^{-1}\sqrt{\frac{k_{1}\log\max(p,\, d)}{n}}\nonumber \\
 &  & +c^{'}b_{1}^{-1}\sqrt{\frac{k_{1}\log\max(p,\, d)}{n}}\sigma_{X^{*}}\sigma_{W}\max\left\{ \sqrt{\frac{k_{1}\log d}{n}},\,\sqrt{\frac{\log p}{n}}\right\} 
\end{eqnarray}
with probability at least $1-c_{1}\exp(-c_{2}\log\max(p,\, d))$.
Notice that by the definition of $\mathbb{K}(k_{1},\, d_{j})$, 
\[
\max_{j^{'},\, j}\sup_{v^{j}\in\mathbb{K}(k_{1},\, d_{j})}\left|\mathbb{E}(x_{1j^{'}}^{*}\mathbf{z}_{1j}v^{j})\right|=\max_{j,j^{'}}|\mathbb{E}(x_{1j^{'}}^{*},\,\mathbf{z}_{1j})|_{\infty}.
\]
To bound the term $\left|\frac{(\hat{X}-X^{*})^{T}(\hat{X}-X^{*})}{n}\right|_{\infty}$,
recalling from the proof for Lemma 6.3, again with a covering subset
argument, the $(j^{'},\, j)$ element of the matrix\textbf{ }$\frac{(\hat{X}-X^{*})^{T}(\hat{X}-X^{*})}{n}$
can be rewritten as follows\\
\\
$\frac{1}{n}(\hat{\mathbf{x}}_{j^{'}}-\mathbf{x}_{j^{'}}^{*})^{T}(\hat{\mathbf{x}}_{j}-\mathbf{x}_{j}^{*})=b_{1}^{-2}\frac{k_{1}\log\max(p,\, d)}{n}\{\frac{1}{n}t_{i^{'}}^{j^{'}T}\mathbf{z}_{j^{'}}^{T}\mathbf{z}_{j}t_{i}^{j}+\frac{1}{n}\triangle v^{j^{'}T}\mathbf{z}_{j^{'}}^{T}\mathbf{z}_{j}t_{i}^{j}+\frac{1}{n}t_{i^{'}}^{j^{'}T}\mathbf{z}_{j^{'}}^{T}\mathbf{z}_{j}\triangle v^{j}+\frac{1}{n}\triangle v^{j^{'}T}\mathbf{z}_{j^{'}}^{T}\mathbf{z}_{j}\triangle v^{j}\}.$
\\
\\
With a similar argument as in the proof for Lemma 6.3, we obtain 
\[
\mathbb{P}\left(\max_{j^{'},\, j}\sup_{v^{j^{'}}\in\mathbb{K}(k_{1},\, d_{j^{'}}),\, v^{j}\in\mathbb{K}(k_{1},\, d_{j})}\left|\frac{1}{n}v^{j^{'}}\mathbf{z}_{j^{'}}^{T}\mathbf{z}_{j}v^{j}-\mathbb{E}(v^{j^{'}}\mathbf{z}_{1j^{'}}^{T}\mathbf{z}_{1j}v^{j})\right|\geq t\right)
\]
\begin{eqnarray*}
 & \leq & p^{2}d^{2k_{1}}\cdot2\exp(-cn\min(\frac{t^{2}}{\sigma_{W}^{4}},\,\frac{t}{\sigma_{W}^{2}}))\\
 & = & 2\exp(-cn\min(\frac{t^{2}}{\sigma_{W}^{4}},\,\frac{t}{\sigma_{W}^{2}})+2k_{1}\log d+2\log p).
\end{eqnarray*}
Consequently, under the condition $\frac{\max(k_{1}\log d,\,\log p)}{n}=o(1)$,
\[
\left|\frac{(\hat{X}-X^{*})^{T}(\hat{X}-X^{*})}{n}-\mathbb{E}[\frac{(\hat{X}-X^{*})^{T}(\hat{X}-X^{*})}{n}]\right|_{\infty}
\]
\begin{eqnarray*}
 & \leq & \left(\max_{j^{'},\, j}\sup_{v^{j^{'}}\in\mathbb{K}(k_{1},\, d_{j^{'}}),\, v^{j}\in\mathbb{K}(k_{1},\, d_{j})}\left|\frac{1}{n}v^{j^{'}}\mathbf{z}_{j^{'}}^{T}\mathbf{z}_{j}v^{j}-\mathbb{E}(v^{j^{'}}\mathbf{z}_{1j^{'}}^{T}\mathbf{z}_{1j}v^{j})\right|\right)b_{1}^{-2}\frac{k_{1}\log\max(p,\, d)}{n}\\
 & \leq & c^{'}b_{1}^{-2}\frac{k_{1}\log\max(p,\, d)}{n}\sigma_{W}^{2}\max\left\{ \sqrt{\frac{k_{1}\log d}{n}},\,\sqrt{\frac{\log p}{n}}\right\} ,
\end{eqnarray*}
with probability at least $1-c_{1}\exp(-c_{2}\log\max(p,\, d))$.
This implies, 
\begin{eqnarray}
\left|\frac{(\hat{X}-X^{*})^{T}(\hat{X}-X^{*})}{n}\right|_{\infty} & \leq & \left(\max_{j^{'},\, j}\sup_{v^{j^{'}}\in\mathbb{K}(k_{1},\, d_{j^{'}}),\, v^{j}\in\mathbb{K}(k_{1},\, d_{j})}\left|\mathbb{E}(v^{j^{'}}\mathbf{z}_{1j^{'}}^{T}\mathbf{z}_{1j}v^{j})\right|\right)b_{1}^{-2}\frac{k_{1}\log\max(p,\, d)}{n}\nonumber \\
 &  & +c^{'}b_{1}^{-2}\frac{k_{1}\log\max(p,\, d)}{n}\sigma_{W}^{2}\max\left\{ \sqrt{\frac{k_{1}\log d}{n}},\,\sqrt{\frac{\log p}{n}}\right\} 
\end{eqnarray}
with probability at least $1-c_{1}\exp(-c_{2}\log\max(p,\, d))$.
Notice that by the definition of $\mathbb{K}(k_{1},\, d_{j})$, 
\[
\max_{j^{'},\, j}\sup_{v^{j^{'}}\in\mathbb{K}(k_{1},\, d_{j^{'}}),\, v^{j}\in\mathbb{K}(k_{1},\, d_{j})}\left|\mathbb{E}(v^{j^{'}}\mathbf{z}_{1j^{'}}^{T}\mathbf{z}_{1j}v^{j})\right|=\max_{j,j^{'}}|\mathbb{E}(\mathbf{z}_{1j^{'}},\,\mathbf{z}_{1j})|_{\infty}.
\]
With exactly the same discretization argument as above, we can show
that, with probability at least $1-c_{1}\exp(-c_{2}\log\max(p,\, d))$,
\begin{eqnarray*}
|\frac{1}{n}(X^{*}-\hat{X})^{T}\boldsymbol{\eta}|_{\infty} & \leq & c^{'}b_{1}^{-1}\sigma_{\eta}\sigma_{W}\sqrt{\frac{k_{1}\log\max(p,\, d)}{n}}\max\left\{ \sqrt{\frac{k_{1}\log d}{n}},\,\sqrt{\frac{\log p}{n}}\right\} ,\\
|\frac{1}{n}(X^{*}-\hat{X})^{T}\epsilon|_{\infty} & \leq & c^{''}b_{1}^{-1}\sigma_{\epsilon}\sigma_{W}\sqrt{\frac{k_{1}\log\max(p,\, d)}{n}}\max\left\{ \sqrt{\frac{k_{1}\log d}{n}},\,\sqrt{\frac{\log p}{n}}\right\} .
\end{eqnarray*}
For the rest of terms in (5), we can use the bounds provided in the
proof for Lemma 6.2. In particular, recall we have, with probability
at least $1-c_{1}\exp(-c_{2}\log p)$, 
\begin{eqnarray*}
|\frac{1}{n}X^{*T}\boldsymbol{\eta}|_{\infty} & \leq & c^{'}\sigma_{X^{*}}\sigma_{\eta}\sqrt{\frac{\log p}{n}},\\
|\frac{1}{n}X^{*T}\epsilon|_{\infty} & \leq & c^{''}\sigma_{X^{*}}\sigma_{\epsilon}\sqrt{\frac{\log p}{n}}.
\end{eqnarray*}
Under the condition $\frac{\max(k_{1}\log d,\, k_{1}\log p)}{n}=o(1)$,
putting everything together yields the claim in Lemma 6.4. $\square$\\

Combining the bounds above with Lemmas 3.1 and 6.3, under the condition
\begin{eqnarray*}
\frac{1}{n}\min\left\{ k_{1}^{2}k_{2}^{2}\log\max(p,\, d),\,\min_{r\in[0,\,1]}\max\left\{ k_{1}^{3-2r}\log d,\, k_{1}^{3-2r}\log p,\, k_{1}^{r}k_{2}\log d,\, k_{1}^{r}k_{2}\log p\right\} \right\}  & = & O(1),\\
\frac{k_{1}\log\max(d,\, p)}{n} & = & o(1),
\end{eqnarray*}
and 
\[
\lambda_{n}\asymp k_{2}\sqrt{\frac{k_{1}\log\max(p,\, d)}{n}},
\]
we have 
\[
|\hat{\beta}_{H2SLS}-\beta^{*}|_{2}\precsim\max\{\varphi_{1}\sqrt{k_{2}}\sqrt{\frac{k_{1}\log\max(p,\, d)}{n}},\,\varphi_{2}\sqrt{\frac{k_{2}\log p}{n}}\},
\]
with probability at least $1-c_{1}\exp(-c_{2}\log\min(p,\, d))$ for
some universal positive constants $c_{1}$ and $c_{2}$, where $\varphi_{1}$
and $\varphi_{2}$ are defined in Theorem 3.2. This proves Theorem
3.3.$\square$

\subsection{Corollary 3.4, Theorems 3.5, and 3.6}

Corollary 3.4 is obvious from Theorem 3.3. The proof for Theorem 3.5 is completely identical to
that for Theorem 3.2 except we replace the inequality $\max_{j=1,...,p}|\hat{\pi}_{j}-\pi_{j}^{*}|_{1}\leq\frac{c\sigma_{\eta}}{\lambda_{\min}(\Sigma_{Z})}k_{1}\sqrt{\frac{\log\max(p,\, d)}{n}}$
by $\max_{j=1,...,p}|\hat{\pi}_{j}-\pi_{j}|_{1}\leq\sqrt{k_{1}}M(d,\, p,\, k_{1},\, n)$.
Also, the proof for Theorem 3.6 is completely identical to that for
Theorem 3.3 except we replace the inequality $\max_{j=1,...,p}|\hat{\pi}_{j}-\pi_{j}^{*}|_{2}\leq\frac{c\sigma_{\eta}}{\lambda_{\min}(\Sigma_{Z})}\sqrt{\frac{k_{1}\log\max(p,\, d)}{n}}$
by $\max_{j=1,...,p}|\hat{\pi}_{j}-\pi_{j}|_{2}\leq M(d,\, p,\, k_{1},\, n)$.

\subsection{Lemma 6.5}

\textbf{Lemma 6.5}: Suppose Assumptions 3.7 and 3.8 hold. Let $J(\beta^{*})=K$,
$\Sigma_{K^{c}K}:=\mathbb{E}\left[X_{1,K^{c}}^{*T}X_{1,K}^{*}\right]$,
$\hat{\Sigma}_{K^{c}K}:=\frac{1}{n}X_{K^{c}}^{*T}X_{K}^{*}$, and
$\tilde{\Sigma}_{K^{c}K}:=\frac{1}{n}\hat{X}_{K^{c}}^{T}\hat{X}_{K}$.
Similarly, let $\Sigma_{KK}:=\mathbb{E}\left[X_{1,K}^{*T}X_{1,K}^{*}\right]$,
$\hat{\Sigma}_{KK}:=\frac{1}{n}X_{K}^{*T}X_{K}^{*}$, and $\tilde{\Sigma}_{KK}:=\frac{1}{n}\hat{X}_{K}^{T}\hat{X}_{K}$.
(i) If the assumptions in Lemmas 6.1 and 6.2 hold, then under the
condition 
\begin{eqnarray*}
\frac{1}{n}\max\left\{ k_{1}k_{2}^{3/2}\log p,\, k_{2}^{3}\log p\right\}  & = & O(1),\\
\frac{1}{n}k_{1}^{2}k_{2}^{2}\log\max(p,\, d) & = & o(1),
\end{eqnarray*}
the sample matrix $\frac{1}{n}\hat{X}^{T}\hat{X}$ satisfies an analogous
version of the {}``mutual incoherence'' assumption with high probability,
i.e., 
\[
\mathbb{P}\left[\left\Vert \frac{1}{n}\hat{X}_{K^{c}}^{T}\hat{X}_{K}\left(\frac{1}{n}\hat{X}_{K}^{T}\hat{X}_{K}\right)^{-1}\right\Vert _{\infty}\geq1-\frac{\phi}{4}\right]\leq O\left(\frac{1}{\min(p,\, d)}\right).
\]
(ii) If the assumptions in Lemmas 6.3 and 6.4 hold, then under the
condition 
\begin{eqnarray*}
\frac{1}{n}\max\left\{ k_{1}^{\nicefrac{1}{2}}k_{2}^{\nicefrac{3}{2}}\log p,\; k_{2}^{3}\log p\right\}  & = & O(1),\\
\frac{1}{n}\min\left\{ k_{1}^{2}k_{2}^{2}\log\max(p,\, d),\,\min_{r\in[0,\,1]}\max\left\{ k_{1}^{3-2r}\log d,\, k_{1}^{3-2r}\log p,\, k_{1}^{r}k_{2}\log d,\, k_{1}^{r}k_{2}\log p\right\} \right\}  & = & o(1),\\
\frac{1}{n}k_{1}k_{2}^{2}\log\max(p,\, d) & = & o(1),
\end{eqnarray*}
the sample matrix $\frac{1}{n}\hat{X}^{T}\hat{X}$ satisfies an analogous
version of the {}``mutual incoherence'' assumption with high probability,
i.e., 
\[
\mathbb{P}\left[\left\Vert \frac{1}{n}\hat{X}_{K^{c}}^{T}\hat{X}_{K}\left(\frac{1}{n}\hat{X}_{K}^{T}\hat{X}_{K}\right)^{-1}\right\Vert _{\infty}\geq1-\frac{\phi}{4}\right]\leq O\left(\frac{1}{\min(p,\, d)}\right).
\]
\\
\textbf{Proof}. I use the following decomposition similar to the method
used in Ravikumar, et. al. (2010)\textbf{ 
\[
\tilde{\Sigma}_{K^{c}K}\tilde{\Sigma}_{KK}^{-1}-\Sigma_{K^{c}K}\Sigma_{KK}^{-1}=R_{1}+R_{2}+R_{3}+R_{4}+R_{5}+R_{6},
\]
}where\textbf{ 
\begin{eqnarray*}
R_{1} & = & \Sigma_{K^{c}K}[\hat{\Sigma}_{KK}^{-1}-\Sigma_{KK}^{-1}],\\
R_{2} & = & [\hat{\Sigma}_{K^{c}K}-\Sigma_{K^{c}K}]\Sigma_{KK}^{-1},\\
R_{3} & = & [\hat{\Sigma}_{K^{c}K}-\Sigma_{K^{c}K}][\hat{\Sigma}_{KK}^{-1}-\Sigma_{KK}^{-1}],\\
R_{4} & = & \hat{\Sigma}_{K^{c}K}[\tilde{\Sigma}_{KK}^{-1}-\hat{\Sigma}_{KK}^{-1}],\\
R_{5} & = & [\tilde{\Sigma}_{K^{c}K}-\hat{\Sigma}_{K^{c}K}]\hat{\Sigma}_{KK}^{-1},\\
R_{6} & = & [\tilde{\Sigma}_{K^{c}K}-\hat{\Sigma}_{K^{c}K}][\tilde{\Sigma}_{KK}^{-1}-\hat{\Sigma}_{KK}^{-1}].
\end{eqnarray*}
}By Assumption 3.7, we have 
\[
\left\Vert \Sigma_{K^{c}K}\Sigma_{KK}^{-1}\right\Vert _{\infty}\leq1-\phi.
\]
It suffices to show that $||R_{i}||_{\infty}\leq\frac{\phi}{6}$\textbf{
}for $i=1,...,3$ and $||R_{i}||_{\infty}\leq\frac{\phi}{12}$\textbf{
}for $i=4,...,6$. 

For the first term $R_{1}$, we have 
\[
R_{1}=-\Sigma_{K^{c}K}\Sigma_{KK}^{-1}[\hat{\Sigma}_{KK}-\Sigma_{KK}]\hat{\Sigma}_{KK}^{-1},
\]
Using the sub-multiplicative property $||AB||_{\infty}\leq||A||_{\infty}||B||_{\infty}$
and the elementary inequality $||A||_{\infty}\leq\sqrt{a}||A||_{2}$
for any symmetric matrix $A\in\mathbb{R}^{a\times a}$, we can bound
$R_{1}$ as follows: 
\begin{eqnarray*}
||R_{1}||_{\infty} & \leq & \left\Vert \Sigma_{K^{c}K}\Sigma_{KK}^{-1}\right\Vert _{\infty}\left\Vert \hat{\Sigma}_{KK}-\Sigma_{KK}\right\Vert _{\infty}\left\Vert \hat{\Sigma}_{KK}^{-1}\right\Vert _{\infty}\\
 & \leq & (1-\phi)\left\Vert \hat{\Sigma}_{KK}-\Sigma_{KK}\right\Vert _{\infty}\sqrt{k_{2}}\left\Vert \hat{\Sigma}_{KK}^{-1}\right\Vert _{2},
\end{eqnarray*}
where the last inequality follows from Assumption 3.7. Using bound
(16) from the proof for Lemma 6.11, we have 
\[
\left\Vert \hat{\Sigma}_{KK}^{-1}\right\Vert _{2}\leq\frac{2}{\lambda_{\min}(\Sigma_{KK})}
\]
with probability at least $1-c_{1}\exp(-c_{2}n)$. Next, applying
bound (11) from Lemma 6.11 with $\varepsilon=\frac{\phi\lambda_{\min}(\Sigma_{KK})}{12(1-\phi)\sqrt{k_{2}}}$,
we have 
\[
\mathbb{P}\left[\left\Vert \hat{\Sigma}_{KK}-\Sigma_{KK}\right\Vert _{\infty}\geq\frac{\phi\lambda_{\min}(\Sigma_{KK})}{12(1-\phi)\sqrt{k_{2}}}\right]\leq2\exp(-bn\min\{\frac{1}{k_{2}^{3}},\,\frac{1}{k_{2}^{3/2}}\}+2\log k_{2}).
\]
Then, we are guaranteed that 
\[
\mathbb{P}[||R_{1}||_{\infty}\geq\frac{\phi}{6}]\leq2\exp(-bn\min\{\frac{1}{k_{2}^{3}},\,\frac{1}{k_{2}^{3/2}}\}+2\log k_{2}).
\]
For the second term $R_{2}$, we first write 
\begin{eqnarray*}
||R_{2}||_{\infty} & \leq & \sqrt{k_{2}}\left\Vert \Sigma_{KK}^{-1}\right\Vert _{2}\left\Vert \hat{\Sigma}_{K^{c}K}-\Sigma_{K^{c}K}\right\Vert _{\infty}\\
 & \leq & \frac{\sqrt{k_{2}}}{\lambda_{\min}(\Sigma_{KK})}\left\Vert \hat{\Sigma}_{K^{c}K}-\Sigma_{K^{c}K}\right\Vert _{\infty}.
\end{eqnarray*}
An application of bound (10) from Lemma 6.11 with $\varepsilon=\frac{\phi}{6}\frac{\lambda_{\min}(\Sigma_{KK})}{\sqrt{k_{2}}}$
to bound the term $\left\Vert \hat{\Sigma}_{K^{c}K}-\Sigma_{K^{c}K}\right\Vert _{\infty}$
yields 
\[
\mathbb{P}[||R_{2}||_{\infty}\geq\frac{\phi}{6}]\leq2\exp(-bn\min\{\frac{1}{k_{2}^{3}},\,\frac{1}{k_{2}^{3/2}}\}+\log(p-k_{2})+\log k_{2}).
\]
For the third term $R_{3}$, by applying bounds (10) from Lemma 6.11
with $\varepsilon=\frac{\phi\lambda_{\min}(\Sigma_{KK})}{6}$ to bound
the term $\left\Vert \hat{\Sigma}_{K^{c}K}-\Sigma_{K^{c}K}\right\Vert _{\infty}$
and (12) from Lemma 6.11 to bound the term $\left\Vert \hat{\Sigma}_{KK}^{-1}-\Sigma_{KK}^{-1}\right\Vert _{\infty}$,
we have

\[
\mathbb{P}[||R_{3}||_{\infty}\geq\frac{\phi}{6}]\leq2\exp(-bn\min\{\frac{1}{k_{2}^{2}},\,\frac{1}{k_{2}}\}+\log(p-k_{2})+\log k_{2}).
\]
Putting everything together, we conclude that 
\[
\mathbb{P}[||\hat{\Sigma}_{K^{c}K}\hat{\Sigma}_{KK}^{-1}||_{\infty}\geq1-\frac{\phi}{2}]\leq O\left(\exp(-bn\min\{\frac{1}{k_{2}^{3}},\,\frac{1}{k_{2}^{3/2}}\}+2\log p)\right).
\]
For the fourth term $R_{4}$, we have, with probability at least $1-c\exp(-bn\min\{\frac{1}{k_{2}^{3}},\,\frac{1}{k_{2}^{3/2}}\}+2\log p)$,
\begin{eqnarray*}
||R_{4}||_{\infty} & \leq & \left\Vert \hat{\Sigma}_{K^{c}K}\hat{\Sigma}_{KK}^{-1}\right\Vert _{\infty}\left\Vert \tilde{\Sigma}_{KK}-\hat{\Sigma}_{KK}\right\Vert _{\infty}\left\Vert \tilde{\Sigma}_{KK}^{-1}\right\Vert _{\infty}\\
 & \leq & (1-\frac{\phi}{2})\left\Vert \tilde{\Sigma}_{KK}-\hat{\Sigma}_{KK}\right\Vert _{\infty}\sqrt{k_{2}}\left\Vert \tilde{\Sigma}_{KK}^{-1}\right\Vert _{2},
\end{eqnarray*}
where the last inequality follows from the bound on $||\hat{\Sigma}_{K^{c}K}\hat{\Sigma}_{KK}^{-1}||_{\infty}$
established previously. Using bounds (24) (or (26)) from the proof
for Lemma 6.12, we have 
\[
\left\Vert \tilde{\Sigma}_{KK}^{-1}\right\Vert _{2}\leq\frac{4}{\lambda_{\min}(\Sigma_{KK})}
\]
with probability at least $1-c_{1}\exp(-c_{2}\log\max(p,\, d))$.
Next, applying bound (18) (or (21)) from Lemma 6.12 with $\varepsilon=\frac{\phi\lambda_{\min}(\Sigma_{KK})}{48(1-\frac{\phi}{2})\sqrt{k_{2}}}$
to bound the term $\left\Vert \tilde{\Sigma}_{KK}-\hat{\Sigma}_{KK}\right\Vert _{\infty}$
yields, 
\begin{eqnarray*}
\mathbb{P}[||R_{4}||_{\infty} & \geq & \frac{\phi}{12}]\leq6\cdot\exp(-bn\min\{\frac{n}{k_{1}^{2}k_{2}^{3}\log\max(p,\, d)},\,\frac{\sqrt{n}}{k_{1}k_{2}^{3/2}\sqrt{\log\max(p,\, d)}}\}+\log d+2\log k_{2})\\
 &  & +c_{1}\exp(-c_{2}\log\max(p,\, d)),
\end{eqnarray*}
or, 
\begin{eqnarray*}
\mathbb{P}[||R_{4}||_{\infty} & \geq & \frac{\phi}{12}]\leq2\cdot\exp(-b^{'}n\min\{\frac{n}{k_{1}k_{2}^{3}\log\max(p,\, d)},\,\frac{\sqrt{n}}{\sqrt{k_{1}}k_{2}^{3/2}\sqrt{\log\max(p,\, d)}}\}+k_{1}\log d+2\log k_{2})\\
 &  & +c_{1}\exp(-c_{2}\log\max(p,\, d)).
\end{eqnarray*}
For the fifth term $R_{5}$, using bound (16) from the proof for Lemma
6.11, we have 
\begin{eqnarray*}
||R_{5}||_{\infty} & \leq & \sqrt{k_{2}}\left\Vert \hat{\Sigma}_{KK}^{-1}\right\Vert _{2}\left\Vert \tilde{\Sigma}_{K^{c}K}-\hat{\Sigma}_{K^{c}K}\right\Vert _{\infty}\\
 & \leq & \frac{2\sqrt{k_{2}}}{\lambda_{\min}(\Sigma_{KK})}\left\Vert \tilde{\Sigma}_{K^{c}K}-\hat{\Sigma}_{K^{c}K}\right\Vert _{\infty}.
\end{eqnarray*}
An application of bound (17) (or (20)) from Lemma 6.12 with $\varepsilon=\frac{\phi\lambda_{\min}(\Sigma_{KK})}{24\sqrt{k_{2}}}$
to bound the term $\left\Vert \tilde{\Sigma}_{K^{c}K}-\hat{\Sigma}_{K^{c}K}\right\Vert _{\infty}$
yields 
\begin{eqnarray*}
\mathbb{P}[||R_{5}||_{\infty} & \geq & \frac{\phi}{12}]\leq6\cdot\exp(-bn\min\{\frac{n}{k_{1}^{2}k_{2}^{3}\log\max(p,\, d)},\,\frac{\sqrt{n}}{k_{1}k_{2}^{3/2}\sqrt{\log\max(p,\, d)}}\}+\log d+2\log p)\\
 &  & +c_{1}\exp(-c_{2}\log\max(p,\, d)),
\end{eqnarray*}
or, 
\begin{eqnarray*}
\mathbb{P}[||R_{5}||_{\infty} & \geq & \frac{\phi}{12}]\leq2\cdot\exp(-b^{'}n\min(\frac{n}{k_{1}k_{2}^{3}\log\max(p,\, d)},\,\frac{\sqrt{n}}{\sqrt{k_{1}}k_{2}^{3/2}\sqrt{\log\max(p,\, d)}})+k_{1}\log d+2\log p)\\
 &  & +c_{1}\exp(-c_{2}\log\max(p,\, d)).
\end{eqnarray*}
For the sixth term $R_{6}$, by applying bounds (17) and (19) (or,
(20) and (22)) to bound the terms $\left\Vert \tilde{\Sigma}_{K^{c}K}-\hat{\Sigma}_{K^{c}K}\right\Vert _{\infty}$
and $\left\Vert \tilde{\Sigma}_{KK}^{-1}-\hat{\Sigma}_{KK}^{-1}\right\Vert _{\infty}$
respectively, with $\varepsilon=\frac{\phi}{12}\frac{\lambda_{\min}(\Sigma_{KK})}{8}$
for (17) (or (20)), we are guaranteed that 
\begin{eqnarray*}
\mathbb{P}[||R_{6}||_{\infty} & \geq & \frac{\phi}{12}]\leq6\cdot\exp(-bn\min\{\frac{n}{k_{1}^{2}k_{2}^{2}\log\max(p,\, d)},\,\frac{\sqrt{n}}{k_{1}k_{2}\sqrt{\log\max(p,\, d)}}\}+\log d+2\log p)\\
 &  & +c_{1}\exp(-c_{2}\log\max(p,\, d)),
\end{eqnarray*}
or, 
\begin{eqnarray*}
\mathbb{P}[||R_{6}||_{\infty} & \geq & \frac{\phi}{12}]\leq2\cdot\exp(-b^{'}n\min(\frac{n}{k_{1}k_{2}^{2}\log\max(p,\, d)},\,\frac{\sqrt{n}}{\sqrt{k_{1}}k_{2}\sqrt{\log\max(p,\, d)}})+k_{1}\log d+2\log p)\\
 &  & +c_{1}\exp(-c_{2}\log\max(p,\, d)).
\end{eqnarray*}
Under the assumptions in Lemmas 6.1 and 6.2 and the condition 
\begin{eqnarray*}
\frac{1}{n}\max\left\{ k_{1}k_{2}^{3/2}\log p,\, k_{2}^{3}\log p\right\}  & = & O(1),\\
\frac{1}{n}k_{1}^{2}k_{2}^{2}\log\max(p,\, d) & = & o(1),
\end{eqnarray*}
or, under the assumptions in Lemmas 6.3 and 6.4 and the condition
\begin{eqnarray*}
\frac{1}{n}\max\left\{ k_{1}^{\nicefrac{1}{2}}k_{2}^{\nicefrac{3}{2}}\log p,\; k_{2}^{3}\log p\right\}  & = & O(1),\\
\frac{1}{n}\min\left\{ k_{1}^{2}k_{2}^{2}\log\max(p,\, d),\,\min_{r\in[0,\,1]}\max\left\{ k_{1}^{3-2r}\log d,\, k_{1}^{3-2r}\log p,\, k_{1}^{r}k_{2}\log d,\, k_{1}^{r}k_{2}\log p\right\} \right\}  & = & o(1),\\
\frac{1}{n}k_{1}k_{2}^{2}\log\max(p,\, d) & = & o(1),
\end{eqnarray*}
putting the bounds on $R_{1}-R_{6}$ together, we conclude that 
\[
\mathbb{P}[||\tilde{\Sigma}_{K^{c}K}\tilde{\Sigma}_{KK}^{-1}||_{\infty}\geq1-\frac{\phi}{4}]\leq O\left(\frac{1}{\min(p,\, d)}\right).
\]
$\square$

\subsection{Theorems 3.7-3.8}

The proof for the first claim in Theorems 3.7 and 3.8 is established
in Lemma 6.6, which shows that $\hat{\beta}_{H2SLS}=(\hat{\beta}_{K},\,\mathbf{0})$
where $\hat{\beta}_{K}$ is the solution obtained in step 2 of the
PDW construction (recall we let $J(\beta^{*}):=K$ and $J(\beta^{*})^{c}:=K^{c}$
for notational convenience in Lemma 6.5). The second and third claims
are proved using Lemma 6.7. The last claim is a consequence of the
third claim.\\
\\
\textbf{Lemma 6.6}: If the PDW construction succeeds, then under Assumption
3.8, the vector $(\hat{\beta}_{K},\,\mathbf{0})\in\mathbb{R}^{p}$
is the unique optimal solution of the Lasso.\\
\\
\textbf{Proof}. The proof for Lemma 6.6 adopts the proof for Lemma
1 from Chapter 6.4.2 of Wainwright (2014). If the PDW construction
succeeds, then $\hat{\beta}=(\hat{\beta}_{K},\,\mathbf{0})$ is an
optimal solution with associated subgradient vector $\hat{\mu}\in\mathbb{R}^{p}$
satisfying $|\hat{\mu}_{K^{c}}|_{\infty}<1$, and $\left\langle \hat{\mu},\,\hat{\beta}\right\rangle =|\hat{\beta}|_{1}$.
Suppose $\tilde{\beta}$ is another optimal solution. Letting $F(\beta)=\frac{1}{2n}|y-\hat{X}\beta|_{2}^{2}$,
then \textbf{$F(\hat{\beta})+\lambda_{n}\left\langle \hat{\mu},\,\hat{\beta}\right\rangle =F(\tilde{\beta})+\lambda_{n}|\tilde{\beta}|_{1}$},
and hence\textbf{ $F(\hat{\beta})-\lambda_{n}\left\langle \hat{\mu},\,\tilde{\beta}-\hat{\beta}\right\rangle =F(\tilde{\beta})+\lambda_{n}\left(|\tilde{\beta}|_{1}-\left\langle \hat{\mu},\,\tilde{\beta}\right\rangle \right)$}.
However, by the zero-subgradient%
\footnote{Given a convex function $g:\,\mathbb{R}^{p}\mapsto\mathbb{R}$, $\mu\in\mathbb{R}^{p}$
is a subgradient at $\beta$, denoted by $\mu\in\partial g(\beta)$,
if $g(\beta+\triangle)\geq g(\beta)+\left\langle \mu,\,\triangle\right\rangle $
for all $\triangle\in\mathbb{R}^{p}$. When $g(\beta)=|\beta|_{1}$,
notice that $\mu\in\partial|\beta|_{1}$ if and only if $\mu_{j}=\textrm{sign}(\beta_{j})$
for all $j=1,...,p$, where $\textrm{sign}(0)$ is allowed to be any
number in $[-1,\,1]$.%
} conditions for optimality, we have \textbf{$\lambda_{n}\hat{\mu}=-\nabla F(\hat{\beta})$},
which implies that \textbf{$F(\hat{\beta})+\left\langle \nabla F(\hat{\beta}),\,\tilde{\beta}-\hat{\beta}\right\rangle -F(\tilde{\beta})=\lambda_{n}\left(|\tilde{\beta}|_{1}-\left\langle \hat{\mu},\,\tilde{\beta}\right\rangle \right)$}.\textbf{
}By convexity of \textbf{$F$}, the left-hand side is non-positive,
which implies that $|\tilde{\beta}|_{1}\leq\left\langle \hat{\mu},\,\tilde{\beta}\right\rangle $.
But since we also have \textbf{$\left\langle \hat{\mu},\,\tilde{\beta}\right\rangle \leq|\hat{\mu}|_{\infty}|\tilde{\beta}|_{1}$},
we must have\textbf{ $|\tilde{\beta}|_{1}=\left\langle \hat{\mu},\,\tilde{\beta}\right\rangle $}.
Since\textbf{ $|\hat{\mu}_{K^{c}}|_{\infty}<1$}, this equality can
only occur if $\tilde{\beta}_{j}=0$ for all\textbf{ $j\in K^{c}$}.
Thus, all optimal solutions are supported only on $K$, and hence
can be obtained by solving the oracle subproblem in the PDW procedure
described in Section 3.2. Given Assumption 3.8, this subproblem is
strictly convex, and hence it has a unique minimizer.\textbf{ $\square$}\\
\textbf{}\\
\textbf{Lemma 6.7}: Suppose Assumptions 1.1, 3.1-3.3, 3.5a, 3.7, and
3.8 hold. Let 
\begin{eqnarray*}
\varphi_{1} & = & \frac{\sigma_{\eta}\max_{j,j^{'}}|\textrm{cov}(x_{1j^{'}}^{*},\,\mathbf{z}_{1j})|_{\infty}|\beta^{*}|_{1}}{\lambda_{\min}(\Sigma_{Z})\lambda_{\min}(\Sigma_{KK})},\\
\varphi_{2} & = & \max\left\{ \frac{\sigma_{X^{*}}\sigma_{\eta}|\beta^{*}|_{1}}{\lambda_{\min}(\Sigma_{KK})},\;\frac{\sigma_{X^{*}}\sigma_{\epsilon}}{\lambda_{\min}(\Sigma_{KK})}\right\} .
\end{eqnarray*}
With the choice of the tuning parameter\\
\\
$\lambda_{n}\asymp\frac{48(2-\frac{\phi}{4})}{\phi}\max\left\{ \frac{\sigma_{\eta}\max_{j,j^{'}}|\textrm{cov}(x_{1j^{'}}^{*},\,\mathbf{z}_{1j})|_{\infty}|\beta^{*}|_{1}}{\lambda_{\min}(\Sigma_{Z})}k_{1}\sqrt{\frac{\log\max(p,\, d)}{n}},\:\sigma_{X^{*}}\sigma_{\eta}|\beta^{*}|_{1}\sqrt{\frac{\log p}{n}},\:\sigma_{X^{*}}\sigma_{\epsilon}\sqrt{\frac{\log p}{n}}\right\} $
\[
\asymp k_{2}k_{1}\sqrt{\frac{\log\max(p,\, d)}{n}},
\]
and under the condition $\frac{\max\{k_{1}^{2}\log d,\, k_{1}^{2}\log p\}}{n}=o(1)$,
we have $|\hat{\mu}_{K^{c}}|_{\infty}\leq1-\frac{\phi}{8}$ with probability
at least $1-c_{1}\exp(-c_{2}\log\min(p,\, d))$. Furthermore, 
\[
|\hat{\beta}_{K}-\beta_{K}^{*}|_{\infty}\leq c\max\left\{ \varphi_{1}k_{1}\sqrt{\frac{k_{2}\log\max(p,\, d)}{n}},\:\varphi_{2}\sqrt{\frac{k_{2}\log p}{n}}\right\} ,
\]
with probability at least $1-c_{1}\exp(-c_{2}\log\min(p,\, d))$.
If Assumptions 1.1, 3.1-3.4, 3.5b, 3.6-3.8 hold, then with the choice
of tuning parameter\\
\\
$\lambda_{n}\asymp\frac{48(2-\frac{\phi}{4})}{\phi}\max\left\{ \frac{\sigma_{\eta}\max_{j,j^{'}}|\textrm{cov}(x_{1j^{'}}^{*},\,\mathbf{z}_{1j})|_{\infty}|\beta^{*}|_{1}}{\lambda_{\min}(\Sigma_{Z})}\sqrt{\frac{k_{1}\log\max(p,\, d)}{n}},\:\sigma_{X^{*}}\sigma_{\eta}|\beta^{*}|_{1}\sqrt{\frac{\log p}{n}},\:\sigma_{X^{*}}\sigma_{\epsilon}\sqrt{\frac{\log p}{n}}\right\} $
\[
\asymp k_{2}\sqrt{\frac{k_{1}\log\max(p,\, d)}{n}},
\]
and under the condition $\frac{\max\{k_{1}\log d,\, k_{1}\log p\}}{n}=o(1)$,
we have $|\hat{\mu}_{K^{c}}|_{\infty}\leq1-\frac{\phi}{8}$ with probability
at least $1-c_{1}\exp(-c_{2}\log\min(p,\, d))$, and 
\[
|\hat{\beta}_{K}-\beta_{K}^{*}|_{\infty}\leq c^{'}\max\left\{ \varphi_{1}\sqrt{\frac{k_{1}k_{2}\log\max(p,\, d)}{n}},\:\varphi_{2}\sqrt{\frac{k_{2}\log p}{n}}\right\} ,
\]
with probability at least $1-c_{1}\exp(-c_{2}\log\min(p,\, d))$.
\\
\\
\textbf{Proof}. By construction, the sub-vectors $\hat{\beta}_{K}$,
$\hat{\mu}_{K}$, and $\hat{\mu}_{K^{c}}$ satisfy the zero-subgradient
condition in the PDW construction. Recall $e:=(X-\hat{X})\beta^{*}+\boldsymbol{\eta}\beta^{*}+\epsilon$
from Lemma 3.1. With the fact that $\hat{\beta}_{K^{c}}=\beta_{K^{c}}^{*}=0$,
we have 
\begin{eqnarray*}
\frac{1}{n}\hat{X}_{K}^{T}\hat{X}_{K}\left(\hat{\beta}_{K}-\beta_{K}^{*}\right)+\frac{1}{n}\hat{X}_{K}^{T}e+\lambda_{n}\hat{\mu}_{K} & = & 0,\\
\frac{1}{n}\hat{X}_{K^{c}}^{T}\hat{X}_{K}\left(\hat{\beta}_{K}-\beta_{K}^{*}\right)+\frac{1}{n}\hat{X}_{K^{c}}^{T}e+\lambda_{n}\hat{\mu}_{K^{c}} & = & 0.
\end{eqnarray*}
From the equations above, by solving for the vector $\hat{\mu}_{K^{c}}\in\mathbb{R}^{p-k_{2}}$,
we obtain 
\begin{eqnarray*}
\hat{\mu}_{K^{c}} & = & -\frac{1}{n\lambda_{n}}\hat{X}_{K^{c}}^{T}\hat{X}_{K}\left(\hat{\beta}_{K}-\beta_{K}^{*}\right)-\hat{X}_{K^{c}}^{T}\frac{e}{n\lambda_{n}},\\
\hat{\beta}_{K}-\beta_{K}^{*} & = & -\left(\frac{1}{n}\hat{X}_{K}^{T}\hat{X}_{K}\right)^{-1}\frac{\hat{X}_{K}^{T}e}{n}-\lambda_{n}\left(\frac{\hat{X}_{K}^{T}\hat{X}_{K}}{n}\right)^{-1}\hat{\mu}_{K},
\end{eqnarray*}
which yields 
\[
\hat{\mu}_{K^{c}}=\left(\tilde{\Sigma}_{K^{c}K}\tilde{\Sigma}_{KK}^{-1}\right)\hat{\mu}_{K}+\left(\hat{X}_{K^{c}}^{T}\frac{e}{n\lambda_{n}}\right)-\left(\tilde{\Sigma}_{K^{c}K}\tilde{\Sigma}_{KK}^{-1}\right)\hat{X}_{K}^{T}\frac{e}{n\lambda_{n}}.
\]
By the triangle inequality, we have 
\[
|\hat{\mu}_{K^{c}}|_{\infty}\leq\left\Vert \tilde{\Sigma}_{K^{c}K}\tilde{\Sigma}_{KK}^{-1}\right\Vert _{\infty}+\left|\hat{X}_{K^{c}}^{T}\frac{e}{n\lambda_{n}}\right|_{\infty}+\left\Vert \tilde{\Sigma}_{K^{c}K}\tilde{\Sigma}_{KK}^{-1}\right\Vert _{\infty}\left|\hat{X}_{K}^{T}\frac{e}{n\lambda_{n}}\right|_{\infty},
\]
where the fact that $|\hat{\mu}_{K}|_{\infty}\leq1$ is used in the
inequality above. By Lemma 6.5, we have $\left\Vert \tilde{\Sigma}_{K^{c}K}\tilde{\Sigma}_{KK}^{-1}\right\Vert _{\infty}\leq1-\frac{\phi}{4}$
with probability at least $1-c\exp(-\log\min(p,\, d))$. Hence, 
\begin{eqnarray*}
|\hat{\mu}_{K^{c}}|_{\infty} & \leq & 1-\frac{\phi}{4}+\left|\hat{X}_{K^{c}}^{T}\frac{e}{n\lambda_{n}}\right|_{\infty}+\left\Vert \tilde{\Sigma}_{K^{c}K}\tilde{\Sigma}_{KK}^{-1}\right\Vert _{\infty}\left|\hat{X}_{K}^{T}\frac{e}{n\lambda_{n}}\right|_{\infty}\\
 & \leq & 1-\frac{\phi}{4}+\left(2-\frac{\phi}{4}\right)\left|\hat{X}^{T}\frac{e}{n\lambda_{n}}\right|_{\infty}.
\end{eqnarray*}
Therefore, it suffices to show that $\left(2-\frac{\phi}{4}\right)\left|\hat{X}^{T}\frac{e}{n\lambda_{n}}\right|_{\infty}\leq\frac{\phi}{8}$
with high probability. This result is established in Lemma 6.13. Thus,
we have $|\hat{\mu}_{K^{c}}|_{\infty}\leq1-\frac{\phi}{8}$ with high
probability. 

It remains to establish a bound on the $l_{\infty}-$norm of the error
$\hat{\beta}_{K}-\beta_{K}^{*}$. By the triangle inequality, we have
\begin{eqnarray*}
|\hat{\beta}_{K}-\beta_{K}^{*}|_{\infty} & \leq & \left|\left(\frac{\hat{X}_{K}^{T}\hat{X}_{K}}{n}\right)^{-1}\frac{\hat{X}_{K}^{T}e}{n}\right|_{\infty}+\lambda_{n}\left\Vert \left(\frac{\hat{X}_{K}^{T}\hat{X}_{K}}{n}\right)^{-1}\right\Vert _{\infty}\\
 & \leq & \left\Vert \left(\frac{\hat{X}_{K}^{T}\hat{X}_{K}}{n}\right)^{-1}\right\Vert _{\infty}\left|\frac{\hat{X}_{K}^{T}e}{n}\right|_{\infty}+\lambda_{n}\left\Vert \left(\frac{\hat{X}_{K}^{T}\hat{X}_{K}}{n}\right)^{-1}\right\Vert _{\infty},
\end{eqnarray*}
Using bound (24) (or (26)) from Lemma 6.12, we have 
\[
\left\Vert \left(\frac{\hat{X}_{K}^{T}\hat{X}_{K}}{n}\right)^{-1}\right\Vert _{\infty}\leq\frac{2\sqrt{k_{2}}}{\lambda_{\min}(\hat{\Sigma}_{KK})}\leq\frac{4\sqrt{k_{2}}}{\lambda_{\min}(\Sigma_{KK})}.
\]
By Lemma 6.2, we have, with probability at least $1-c_{1}\exp(-c_{2}\log\min(p,\, d))$,
\[
|\frac{1}{n}\hat{X}^{T}e|_{\infty}\precsim\max\left\{ \frac{\sigma_{\eta}\max_{j,j^{'}}|\textrm{cov}(x_{1j^{'}}^{*},\,\mathbf{z}_{1j})|_{\infty}|\beta^{*}|_{1}}{\lambda_{\min}(\Sigma_{Z})}k_{1}\sqrt{\frac{\log\max(p,\, d)}{n}},\:\sigma_{X^{*}}\sigma_{\eta}|\beta^{*}|_{1}\sqrt{\frac{\log p}{n}},\:\sigma_{X^{*}}\sigma_{\epsilon}\sqrt{\frac{\log p}{n}}\right\} .
\]
By Lemma 6.4, we have, with probability at least $1-c_{1}\exp(-c_{2}\log\min(p,\, d))$),
\[
|\frac{1}{n}\hat{X}^{T}e|_{\infty}\precsim\max\left\{ \frac{\sigma_{\eta}\max_{j,j^{'}}|\textrm{cov}(x_{1j^{'}}^{*},\,\mathbf{z}_{1j})|_{\infty}|\beta^{*}|_{1}}{\lambda_{\min}(\Sigma_{Z})}\sqrt{\frac{k_{1}\log\max(p,\, d)}{n}},\:\sigma_{X^{*}}\sigma_{\eta}|\beta^{*}|_{1}\sqrt{\frac{\log p}{n}},\:\sigma_{X^{*}}\sigma_{\epsilon}\sqrt{\frac{\log p}{n}}\right\} .
\]
Putting everything together, with the choice of $\lambda_{n}$ given
in the statement of Lemma 6.7, we obtain 
\[
|\hat{\beta}_{K}-\beta_{K}^{*}|_{\infty}\leq c\max\left\{ \varphi_{1}k_{1}\sqrt{\frac{k_{2}\log\max(p,\, d)}{n}},\:\varphi_{2}\sqrt{\frac{k_{2}\log p}{n}}\right\} ,
\]
or, 
\[
|\hat{\beta}_{K}-\beta_{K}^{*}|_{\infty}\leq c^{'}\max\left\{ \varphi_{1}\sqrt{\frac{k_{1}k_{2}\log\max(p,\, d)}{n}},\:\varphi_{2}\sqrt{\frac{k_{2}\log p}{n}}\right\} ,
\]
with probability at least $1-c_{1}\exp(-c_{2}\log\min(p,\, d))$),
as claimed. $\square$

\subsection{Lemmas 6.8-6.13}

\textbf{Lemma 6.8}: If $X\in\mathbb{R}^{n\times p_{1}}$ is a zero-mean
sub-Gaussian matrix with parameters $(\Sigma_{X},\,\sigma_{X}^{2})$,
then for any fixed (unit) vector $v\in\mathbb{R}^{p_{1}}$, we have
\[
\mathbb{P}(\left||Xv|_{2}^{2}-\mathbb{E}[|Xv|_{2}^{2}]\right|\geq nt)\leq2\exp(-cn\min\{\frac{t^{2}}{\sigma_{X}^{4}},\,\frac{t}{\sigma_{X}^{2}}\}).
\]
Moreover, if $Y\in\mathbb{R}^{n\times p_{2}}$ is a zero-mean sub-Gaussian
matrix with parameters $(\Sigma_{Y},\,\sigma_{Y}^{2})$, then 
\[
\mathbb{P}(|\frac{Y^{T}X}{n}-\textrm{cov}(\mathbf{y}_{i},\,\mathbf{x}_{i})|_{\infty}\geq t)\leq6p_{1}p_{2}\exp(-cn\min\{\frac{t^{2}}{\sigma_{X}^{2}\sigma_{Y}^{2}},\,\frac{t}{\sigma_{X}\sigma_{Y}}\}),
\]
where $\mathbf{x}_{i}$ and $\mathbf{y}_{i}$ are the $i^{th}$ rows
of $X$ and $Y$, respectively. In particular, if $n\succsim\log p$,
then 
\[
\mathbb{P}(|\frac{Y^{T}X}{n}-\textrm{cov}(\mathbf{y}_{i},\,\mathbf{x}_{i})|_{\infty}\geq c_{0}\sigma_{X}\sigma_{Y}\sqrt{\frac{\log\left(\max\{p_{1},\, p_{2}\}\right)}{n}})\leq c_{1}\exp(-c_{2}\log\left(\max\{p_{1},\, p_{2}\}\right)).
\]
\textbf{Remark}. Lemma 6.8 is Lemma 14 in Loh and Wainwright (2012).
\\
\textbf{}\\
\textbf{Lemma 6.9}: For a fixed matrix $\Gamma\in\mathbb{R}^{p\times p}$,
parameter $s\geq1$, and tolerance $\tau>0$, suppose we have the
deviation condition 
\[
|v^{T}\Gamma v|\leq\tau\qquad\forall v\in\mathbb{K}(2s,\, p).
\]
Then, 
\[
|v^{T}\Gamma v|\leq27\tau\left(|v|_{2}^{2}+\frac{1}{s}|v|_{1}^{2}\right)\qquad\forall v\in\mathbb{R}^{p}.
\]
\\
\textbf{Remark}. Lemma 6.9 is Lemma 12 in Loh and Wainwright (2012).
\\
\textbf{}\\
\textbf{Lemma 6.10}: Under Assumptions 1.1 and 3.3, we have 
\[
\frac{|X^{*}v^{0}|_{2}^{2}}{n}\geq\kappa_{1}|v^{0}|_{2}^{2}-\kappa_{2}\frac{k_{1}^{r}\log\max(p,\, d)}{n}|v^{0}|_{1}^{2},\qquad\textrm{for all }v^{0}\in\mathbb{R}^{p},\, r\in[0,\,1]
\]
with probability at least $1-c_{1}\exp(-c_{2}n)$, where $\kappa_{1}=\frac{\lambda_{\min}(\Sigma_{X^{*}})}{2}$
and $\kappa_{2}=c_{0}\lambda_{\min}(\Sigma_{X^{*}})\max\left\{ \frac{\sigma_{X^{*}}^{4}}{\lambda_{\min}^{2}(\Sigma_{X^{*}})},\,1\right\} $.
\\
\\
\textbf{Proof}. First, we show 
\[
\sup_{v^{0}\in\mathbb{K}(2s,\, p)}\left|v^{0T}\left(\frac{X^{*T}X^{*}}{n}-\Sigma_{X^{*}}\right)v^{0}\right|\leq\frac{\lambda_{\min}(\Sigma_{X^{*}})}{54}
\]
with high probability, where $\Sigma_{X^{*}}=\mathbb{E}(X^{*T}X^{*})$.
Under Assumption 3.3, we have that $X^{*}$ is sub-Gaussian with parameters
$(\Sigma_{X^{*}},\,\sigma_{X^{*}})$. Therefore, by Lemma 6.8 and
a discretization argument similar to those in Lemma 6.3, we have 
\[
\mathbb{P}\left(\sup_{v^{0}\in\mathbb{K}(2s,\, p)}\left|v^{0T}\left(\frac{X^{*T}X^{*}}{n}-\Sigma_{X^{*}}\right)v^{0}\right|\geq t\right)\leq2\exp(-cn\min(\frac{t^{2}}{\sigma_{X^{*}}^{4}},\,\frac{t}{\sigma_{X^{*}}^{2}})+2s\log p),
\]
for some universal constants $c>0$. By choosing $t=\frac{\lambda_{\min}(\Sigma_{X^{*}})}{54}$
and 
\[
s=s(r):=\frac{1}{c^{'}}\frac{n}{k_{1}^{r}\log\max(p,\, d)}\min\left\{ \frac{\lambda_{\min}^{2}(\Sigma_{X^{*}})}{\sigma_{X^{*}}^{4}},\,1\right\} ,\quad r\in[0,\,1],
\]
where $c^{'}$ is chosen sufficiently small so that $s\geq1$, we
get 
\[
\mathbb{P}\left(\sup_{v^{0}\in\mathbb{K}(2s,\, p)}\left|v^{0T}\left(\frac{X^{*T}X^{*}}{n}-\Sigma_{X^{*}}\right)v^{0}\right|\geq\frac{\lambda_{\min}(\Sigma_{X^{*}})}{54}\right)\leq2\exp(-c_{2}n\min(\frac{\lambda_{\min}^{2}(\Sigma_{X^{*}})}{\sigma_{X^{*}}^{4}},\,1)).
\]
Now, by Lemma 6.9 and the following substitutions 
\[
\Gamma=\frac{X^{*T}X^{*}}{n}-\Sigma_{X^{*}},\quad\textrm{and}\quad\tau:=\frac{\lambda_{\min}(\Sigma_{X^{*}})}{54},
\]
we obtain 
\[
\left|v^{0T}\left(\frac{X^{*T}X^{*}}{n}-\Sigma_{X^{*}}\right)v^{0}\right|\leq\frac{\lambda_{\min}(\Sigma_{X^{*}})}{2}\left(|v^{0}|_{2}^{2}+\frac{1}{s}|v^{0}|_{1}^{2}\right),
\]
which implies 
\[
v^{0T}\frac{X^{*T}X^{*}}{n}v^{0}\geq v^{0T}\Sigma_{X^{*}}v^{0}-\frac{\lambda_{\min}(\Sigma_{X^{*}})}{2}\left(|v^{0}|_{2}^{2}+\frac{1}{s}|v^{0}|_{1}^{2}\right).
\]
Recalling the choice of 

\[
s=s(r):=\frac{1}{c^{'}}\frac{n}{k_{1}^{r}\log\max(p,\, d)}\min\left\{ \frac{\lambda_{\min}^{2}(\Sigma_{X^{*}})}{\sigma_{X^{*}}^{4}},\,1\right\} ,\quad r\in[0,\,1],
\]
where $c^{'}$ is chosen sufficiently small so $s\geq1$, the claim
follows. $\square$\\
\\
\textbf{Lemma 6.11}: Suppose Assumptions 1.1, 3.3, and 3.8 hold. For
any $\varepsilon>0$ and constant $c$, we have

\begin{equation}
\mathbb{P}\left\{ \left\Vert \hat{\Sigma}_{K^{c}K}-\Sigma_{K^{c}K}\right\Vert _{\infty}\geq\varepsilon\right\} \leq(p-k_{2})k_{2}\cdot2\exp(-cn\min\{\frac{\varepsilon^{2}}{4k_{2}^{2}\sigma_{X^{*}}^{4}},\,\frac{\varepsilon}{2k_{2}\sigma_{X^{*}}^{2}}\}),
\end{equation}
\begin{equation}
\mathbb{P}\left\{ \left\Vert \hat{\Sigma}_{KK}-\Sigma_{KK}\right\Vert _{\infty}\geq\varepsilon\right\} \leq k_{2}^{2}\cdot2\exp(-cn\min\{\frac{\varepsilon^{2}}{4k_{2}^{2}\sigma_{X^{*}}^{4}},\,\frac{\varepsilon}{2k_{2}\sigma_{X^{*}}^{2}}\}).
\end{equation}
Furthermore, under the scaling $n\succsim k_{2}\log p$, for constants
$b_{1}$ and $b_{2}$, we have

\begin{equation}
\left\Vert \hat{\Sigma}_{KK}^{-1}-\Sigma_{KK}^{-1}\right\Vert _{\infty}\leq\frac{1}{\lambda_{\min}(\Sigma_{KK})},
\end{equation}
with probability at least $1-c_{1}\exp(-c_{2}n\min\{\frac{\lambda_{\min}^{2}(\Sigma_{KK})}{4k_{2}\sigma_{X^{*}}^{4}},\,\frac{\lambda_{\min}(\Sigma_{KK})}{2\sqrt{k_{2}}\sigma_{X^{*}}^{2}}\})$.\textbf{}\\
\textbf{Proof}. Denote the element $(j^{'},\, j)$ of the matrix difference
$\hat{\Sigma}_{K^{c}K}-\Sigma_{K^{c}K}$ by $u_{j^{'}j}$. By the
definition of the $l_{\infty}$matrix norm, we have 
\begin{eqnarray*}
\mathbb{P}\left\{ \left\Vert \hat{\Sigma}_{K^{c}K}-\Sigma_{K^{c}K}\right\Vert _{\infty}\geq\varepsilon\right\}  & = & \mathbb{P}\left\{ \max_{j^{'}\in K^{c}}\sum_{j\in K}|u_{j^{'}j}|\geq\varepsilon\right\} \\
 & \leq & (p-k_{2})\mathbb{P}\left\{ \sum_{j\in K}|u_{j^{'}j}|\geq\varepsilon\right\} \\
 & \leq & (p-k_{2})\mathbb{P}\left\{ \exists j\in K\,\vert\,|u_{j^{'}j}|\geq\frac{\varepsilon}{k_{2}}\right\} \\
 & \leq & (p-k_{2})k_{2}\mathbb{P}\left\{ |u_{j^{'}j}|\geq\frac{\varepsilon}{k_{2}}\right\} \\
 & \leq & (p-k_{2})k_{2}\cdot2\exp(-cn\min\{\frac{\varepsilon^{2}}{k_{2}^{2}\sigma_{X^{*}}^{4}},\,\frac{\varepsilon}{k_{2}\sigma_{X^{*}}^{2}}\}),
\end{eqnarray*}
where the last inequality follows the deviation bound for sub-exponential
random variables, i.e., Lemma 6.8. Bound (11) can be obtained in a
similar way except that the pre-factor $(p-k_{2})$ is replaced by
$k_{2}$. To prove the last bound (12), write 
\begin{eqnarray}
\left\Vert \hat{\Sigma}_{KK}^{-1}-\Sigma_{KK}^{-1}\right\Vert _{\infty} & = & \left\Vert \Sigma_{KK}^{-1}\left[\Sigma_{KK}-\hat{\Sigma}_{KK}\right]\hat{\Sigma}_{KK}^{-1}\right\Vert _{\infty}\nonumber \\
 & = & \sqrt{k_{2}}\left\Vert \Sigma_{KK}^{-1}\left[\Sigma_{KK}-\hat{\Sigma}_{KK}\right]\hat{\Sigma}_{KK}^{-1}\right\Vert _{2}\nonumber \\
 & = & \sqrt{k_{2}}\left\Vert \Sigma_{KK}^{-1}\right\Vert _{2}\left\Vert \Sigma_{KK}-\hat{\Sigma}_{KK}\right\Vert _{2}\left\Vert \hat{\Sigma}_{KK}^{-1}\right\Vert _{2}\nonumber \\
 & \leq & \frac{\sqrt{k_{2}}}{\lambda_{\min}(\Sigma_{KK})}\left\Vert \Sigma_{KK}-\hat{\Sigma}_{KK}\right\Vert _{2}\left\Vert \hat{\Sigma}_{KK}^{-1}\right\Vert _{2}.
\end{eqnarray}
To bound the term $\left\Vert \Sigma_{KK}-\hat{\Sigma}_{KK}\right\Vert _{2}$
in (13), applying Lemma 6.8 with $X^{T}X=\hat{\Sigma}_{KK}$ and $t=\frac{\lambda_{\min}(\Sigma_{KK})}{2\sqrt{k_{2}}}$
yields 
\[
\left\Vert \hat{\Sigma}_{KK}-\Sigma_{KK}\right\Vert _{2}\leq\frac{\lambda_{\min}(\Sigma_{KK})}{2\sqrt{k_{2}}},
\]
with probability at least $1-c_{1}\exp(-c_{2}n\min\{\frac{\lambda_{\min}^{2}(\Sigma_{KK})}{4k_{2}\sigma_{X^{*}}^{4}},\,\frac{\lambda_{\min}(\Sigma_{KK})}{2\sqrt{k_{2}}\sigma_{X^{*}}^{2}}\})$. 

To bound the term\textbf{$\left\Vert \hat{\Sigma}_{KK}^{-1}\right\Vert _{2}$}
in (13), note that we can write 
\begin{eqnarray}
\lambda_{\min}(\Sigma_{KK}) & = & \min_{||h^{'}||_{2}=1}h^{'T}\Sigma_{KK}h^{'}\nonumber \\
 & = & \min_{||h^{'}||_{2}=1}\left[h^{'T}\hat{\Sigma}_{KK}h^{'}+h^{'T}(\Sigma_{KK}-\hat{\Sigma}_{KK})h^{'}\right]\nonumber \\
 & \leq & h^{T}\hat{\Sigma}_{KK}h+h^{T}(\Sigma_{KK}-\hat{\Sigma}_{KK})h
\end{eqnarray}
where $h\in\mathbb{R}^{k_{2}}$ is a unit-norm minimal eigenvector
of $\hat{\Sigma}_{KK}$. Applying Lemma 6.8 yields 
\[
\left|h^{T}\left(\Sigma_{KK}-\hat{\Sigma}_{KK}\right)h\right|\leq\frac{\lambda_{\min}(\Sigma_{KK})}{2}
\]
with probability at least $1-c_{1}\exp(-c_{2}n)$. Therefore, 
\begin{eqnarray}
\lambda_{\min}(\Sigma_{KK}) & \leq & \lambda_{\min}(\hat{\Sigma}_{KK})+\frac{\lambda_{\min}(\Sigma_{KK})}{2}\nonumber \\
\Longrightarrow\lambda_{\min}(\hat{\Sigma}_{KK}) & \geq & \frac{\lambda_{\min}(\Sigma_{KK})}{2},
\end{eqnarray}
and consequently,
\begin{equation}
\left\Vert \hat{\Sigma}_{KK}^{-1}\right\Vert _{2}\leq\frac{2}{\lambda_{\min}(\Sigma_{KK})}.
\end{equation}
Putting everything together, we have 
\[
\left\Vert \hat{\Sigma}_{KK}^{-1}-\Sigma_{KK}^{-1}\right\Vert _{\infty}\leq\frac{\sqrt{k_{2}}}{\lambda_{\min}(\Sigma_{KK})}\frac{\lambda_{\min}(\Sigma_{KK})}{2\sqrt{k_{2}}}\frac{2}{\lambda_{\min}(\Sigma_{KK})}=\frac{1}{\lambda_{\min}(\Sigma_{KK})}.
\]
with probability at least $1-c_{1}\exp(-c_{2}n\min\{\frac{\lambda_{\min}^{2}(\Sigma_{KK})}{4k_{2}\sigma_{X^{*}}^{4}},\,\frac{\lambda_{\min}(\Sigma_{KK})}{2\sqrt{k_{2}}\sigma_{X^{*}}^{2}}\})$.\textbf{
$\square$}\\
\textbf{}\\
\textbf{Lemma 6.12}:\textbf{ }(i) Suppose the assumptions in Lemmas
6.1 and 6.2 hold. For any $\varepsilon>0$, under the condition $\frac{k_{1}^{2}k_{2}^{2}\log\max(p,\, d)}{n}=o(1)$,
we have 
\[
\mathbb{P}\left\{ \left\Vert \tilde{\Sigma}_{K^{c}K}-\hat{\Sigma}_{K^{c}K}\right\Vert _{\infty}\geq\varepsilon\right\} \leq
\]
\[
6(p-k_{2})k_{2}\cdot\exp(-cn\min\{\frac{n\lambda_{\min}^{2}(\Sigma_{Z})\varepsilon^{2}}{\sigma_{\eta}^{2}\sigma_{X^{*}}^{2}\sigma_{Z}^{2}k_{1}^{2}k_{2}^{2}\log\max(p,\, d)},\,\frac{\sqrt{n}\lambda_{\min}(\Sigma_{Z})\varepsilon}{\sigma_{\eta}\sigma_{X^{*}}\sigma_{Z}k_{1}k_{2}\sqrt{\log\max(p,\, d)}}\}+\log d)
\]
\begin{equation}
+c_{1}\exp(-c_{2}\log\max(p,\, d)),
\end{equation}
\[
\mathbb{P}\left\{ \left\Vert \tilde{\Sigma}_{KK}-\hat{\Sigma}_{KK}\right\Vert _{\infty}\geq\varepsilon\right\} \leq
\]
\[
6k_{2}^{2}\cdot\exp(-cn\min\{\frac{n\lambda_{\min}^{2}(\Sigma_{Z})\varepsilon^{2}}{\sigma_{\eta}^{2}\sigma_{X^{*}}^{2}\sigma_{Z}^{2}k_{1}^{2}k_{2}^{2}\log\max(p,\, d)},\,\frac{\sqrt{n}\lambda_{\min}(\Sigma_{Z})\varepsilon}{\sigma_{\eta}\sigma_{X^{*}}\sigma_{Z}k_{1}k_{2}\sqrt{\log\max(p,\, d)}}\}+\log d)
\]
\begin{equation}
+c_{1}\exp(-c_{2}\log\max(p,\, d)).
\end{equation}
Furthermore, we have

\begin{equation}
\left\Vert \tilde{\Sigma}_{KK}^{-1}-\hat{\Sigma}_{KK}^{-1}\right\Vert _{\infty}\leq\frac{8}{\lambda_{\min}(\Sigma_{KK})}\quad\textrm{with probability at least }1-c_{1}\exp(-c_{2}\log\max(p,\, d)).
\end{equation}
(ii) Suppose the assumptions in Lemmas 6.3 and 6.4 hold. For any $\varepsilon>0$,
under the condition $\frac{k_{1}k_{2}^{2}\log\max(p,\, d)}{n}=o(1)$,
we have 
\[
\mathbb{P}\left\{ \left\Vert \tilde{\Sigma}_{K^{c}K}-\hat{\Sigma}_{K^{c}K}\right\Vert _{\infty}\geq\varepsilon\right\} \leq
\]
\[
2(p-k_{2})k_{2}\cdot\exp(-cn\min(\frac{n\lambda_{\min}^{2}(\Sigma_{Z})\varepsilon^{2}}{\sigma_{\eta}^{2}\sigma_{X^{*}}^{2}\sigma_{W}^{2}k_{1}k_{2}^{2}\log\max(p,\, d)},\,\frac{\sqrt{n}\lambda_{\min}(\Sigma_{Z})\varepsilon}{\sigma_{\eta}\sigma_{X^{*}}\sigma_{W}\sqrt{k_{1}}k_{2}\sqrt{\log\max(p,\, d)}})+k_{1}\log d)
\]
\begin{equation}
+c_{1}\exp(-c_{2}\log\max(p,\, d)),
\end{equation}
\[
\mathbb{P}\left\{ \left\Vert \tilde{\Sigma}_{KK}-\hat{\Sigma}_{KK}\right\Vert _{\infty}\geq\varepsilon\right\} \leq
\]
\[
2k_{2}^{2}\cdot\exp(-cn\min(\frac{n\lambda_{\min}^{2}(\Sigma_{Z})\varepsilon^{2}}{\sigma_{\eta}^{2}\sigma_{X^{*}}^{2}\sigma_{W}^{2}k_{1}k_{2}^{2}\log\max(p,\, d)},\,\frac{\sqrt{n}\lambda_{\min}(\Sigma_{Z})\varepsilon}{\sigma_{\eta}\sigma_{X^{*}}\sigma_{W}\sqrt{k_{1}}k_{2}\sqrt{\log\max(p,\, d)}})+k_{1}\log d)
\]
\begin{equation}
+c_{1}\exp(-c_{2}\log\max(p,\, d)).
\end{equation}
Furthermore, if 
\[
\frac{1}{n}\min\left\{ k_{1}^{2}k_{2}^{2}\log\max(p,\, d),\,\min_{r\in[0,\,1]}\max\left\{ k_{1}^{3-2r}\log d,\, k_{1}^{3-2r}\log p,\, k_{1}^{r}k_{2}\log d,\, k_{1}^{r}k_{2}\log p\right\} \right\} =o(1),
\]
we have

\begin{equation}
\left\Vert \tilde{\Sigma}_{KK}^{-1}-\hat{\Sigma}_{KK}^{-1}\right\Vert _{\infty}\leq\frac{8}{\lambda_{\min}(\Sigma_{KK})}\quad\textrm{with probability at least }1-c_{1}\exp(-c_{2}\log\max(p,\, d)).
\end{equation}
\textbf{Proof}. Denote the element $(j^{'},\, j)$ of the matrix difference
$\tilde{\Sigma}_{K^{c}K}-\hat{\Sigma}_{K^{c}K}$ by $w_{j^{'}j}$.
Using the same argument as in Lemma 6.11, we have 
\[
\mathbb{P}\left\{ \left\Vert \tilde{\Sigma}_{K^{c}K}-\hat{\Sigma}_{K^{c}K}\right\Vert _{\infty}\geq\varepsilon\right\} \leq(p-k_{2})k_{2}\mathbb{P}\left\{ |w_{j^{'}j}|\geq\frac{\varepsilon}{k_{2}}\right\} .
\]
Following the derivation of the upper bounds on $\left|\frac{(\hat{X}-X^{*})^{T}X^{*}}{n}\right|_{\infty}$
and $\left|\frac{(\hat{X}-X^{*})^{T}(\hat{X}-X^{*})}{n}\right|_{\infty}$
in the proof for Lemma 6.2 and the identity 
\[
\frac{1}{n}\left(\tilde{\Sigma}_{K^{c}K}-\hat{\Sigma}_{K^{c}K}\right)=\frac{1}{n}X_{K^{c}}^{*T}(\hat{X}_{K}-X_{K}^{*})+\frac{1}{n}(\hat{X}_{K^{c}}-X_{K^{c}}^{*})^{T}X_{K}^{*}+\frac{1}{n}(\hat{X}_{K^{c}}-X_{K^{c}}^{*})^{T}(\hat{X}_{K}-X_{K}^{*}),
\]
we notice that to upper bound $|w_{j^{'}j}|$, it suffices to upper
bound $3\cdot\left|\frac{1}{n}\mathbf{x}_{j^{'}}^{*T}(\hat{\mathbf{x}}_{j}-\mathbf{x}_{j}^{*})\right|$.
From the proof for Lemma 6.2, we have 
\[
\left|\frac{1}{n}\mathbf{x}_{j^{'}}^{*T}(\hat{\mathbf{x}}_{j}-\mathbf{x}_{j}^{*})\right|\leq\frac{c\sigma_{\eta}}{\lambda_{\min}(\Sigma_{Z})}k_{1}\sqrt{\frac{\log\max(p,\, d)}{n}}\left|\frac{1}{n}\sum_{i=1}^{n}x_{ij^{'}}^{*}\mathbf{z}_{ij}\right|_{\infty},
\]
with probability at least $1-c_{1}\exp(-c_{2}\log\max(p,\, d))$ and
\[
\mathbb{P}\left[|\frac{1}{n}\mathbf{x}_{j^{'}}^{*T}Z_{j}-\mathbb{E}(x_{ij^{'}}^{*},\,\mathbf{z}_{ij})|{}_{\infty}\geq t\right]\leq6\exp(-cn\min\{\frac{t^{2}}{\sigma_{X^{*}}^{2}\sigma_{Z}^{2}},\,\frac{t}{\sigma_{X^{*}}\sigma_{Z}}\}+\log d).
\]
Under the condition $\frac{k_{1}^{2}k_{2}^{2}\log\max(p,\, d)}{n}=o(1)$,
setting $t=\frac{\varepsilon\lambda_{\min}(\Sigma_{Z})}{c\sigma_{\eta}}\sqrt{\frac{n}{k_{2}^{2}k_{1}^{2}\log\max(p,\, d)}}$
for any $\varepsilon>0$ yields 
\[
\mathbb{P}\left[|w_{j^{'}j}|\geq\frac{\varepsilon}{k_{2}}\right]\leq
\]
\[
6\exp(-cn\min\{\frac{n\lambda_{\min}^{2}(\Sigma_{Z})\varepsilon^{2}}{\sigma_{\eta}^{2}\sigma_{X^{*}}^{2}\sigma_{Z}^{2}k_{1}^{2}k_{2}^{2}\log\max(p,\, d)},\,\frac{\sqrt{n}\lambda_{\min}(\Sigma_{Z})\varepsilon}{\sigma_{\eta}\sigma_{X^{*}}\sigma_{Z}k_{1}k_{2}\sqrt{\log\max(p,\, d)}}\}+\log d)
\]
\[
+c_{1}\exp(-c_{2}\log\max(p,\, d)).
\]
Therefore, 
\[
\mathbb{P}\left\{ \left\Vert \tilde{\Sigma}_{K^{c}K}-\hat{\Sigma}_{K^{c}K}\right\Vert _{\infty}\geq\varepsilon\right\} \leq
\]
\[
6(p-k_{2})k_{2}\cdot\exp(-cn\min\{\frac{n\lambda_{\min}^{2}(\Sigma_{Z})\varepsilon^{2}}{\sigma_{\eta}^{2}\sigma_{X^{*}}^{2}\sigma_{Z}^{2}k_{1}^{2}k_{2}^{2}\log\max(p,\, d)},\,\frac{\sqrt{n}\lambda_{\min}(\Sigma_{Z})\varepsilon}{\sigma_{\eta}\sigma_{X^{*}}\sigma_{Z}k_{1}k_{2}\sqrt{\log\max(p,\, d)}}\}+\log d)
\]
\[
+c_{1}\exp(-c_{2}\log\max(p,\, d)).
\]
Bound (18) can be obtained in a similar way except that the pre-factor
$(p-k_{2})$ is replaced by $k_{2}$. 

To prove bound (19), by applying the same argument as in Lemma 6.11,
we have 
\begin{eqnarray*}
\left\Vert \tilde{\Sigma}_{KK}^{-1}-\hat{\Sigma}_{KK}^{-1}\right\Vert _{\infty} & \leq & \frac{\sqrt{k_{2}}}{\lambda_{\min}(\hat{\Sigma}_{KK})}\left\Vert \hat{\Sigma}_{KK}-\tilde{\Sigma}_{KK}\right\Vert _{2}\left\Vert \tilde{\Sigma}_{KK}^{-1}\right\Vert _{2}\\
 & \leq & \frac{2\sqrt{k_{2}}}{\lambda_{\min}(\Sigma_{KK})}\left\Vert \hat{\Sigma}_{KK}-\tilde{\Sigma}_{KK}\right\Vert _{2}\left\Vert \tilde{\Sigma}_{KK}^{-1}\right\Vert _{2},
\end{eqnarray*}
where the last inequality comes from bound (15). 

To bound the term $\left\Vert \hat{\Sigma}_{KK}-\tilde{\Sigma}_{KK}\right\Vert _{2}$,
applying bound (18) with $\varepsilon=\frac{\lambda_{\min}(\Sigma_{KK})}{\sqrt{k_{2}}}$
yields 
\begin{eqnarray}
\left\Vert \hat{\Sigma}_{KK}-\tilde{\Sigma}_{KK}\right\Vert _{2} & \leq & \left\Vert \hat{\Sigma}_{KK}-\tilde{\Sigma}_{KK}\right\Vert _{\infty}\nonumber \\
 & \leq & \frac{\lambda_{\min}(\Sigma_{KK})}{\sqrt{k_{2}}},
\end{eqnarray}
with probability at least
\[
1-6\cdot\exp(-cn\min\{\frac{n\lambda_{\min}^{4}(\Sigma_{KK})}{\sigma_{\eta}^{2}\sigma_{X^{*}}^{2}\sigma_{Z}^{2}k_{1}^{2}k_{2}^{3}\log\max(p,\, d)},\,\frac{\sqrt{n}\lambda_{\min}^{2}(\Sigma_{KK})}{\sigma_{\eta}\sigma_{X^{*}}\sigma_{Z}k_{1}k_{2}^{\nicefrac{3}{2}}\sqrt{\log\max(p,\, d)}}\}+\log d+2\log k_{2})
\]
\[
-c_{1}\exp(-c_{2}\log\max(p,\, d))\geq1-O\left(\frac{1}{\max(p,\, d)}\right)
\]
if $\frac{k_{1}k_{2}^{3/2}\log\max(p,\, d)}{n}=O(1)$.

To bound the term\textbf{ $\left\Vert \tilde{\Sigma}_{KK}^{-1}\right\Vert _{2}$},
again we have, 
\begin{eqnarray*}
\lambda_{\min}(\hat{\Sigma}_{KK}) & \leq & h^{T}\tilde{\Sigma}_{KK}h+h^{T}(\hat{\Sigma}_{KK}-\tilde{\Sigma}_{KK})h\\
 & \leq & h^{T}\tilde{\Sigma}_{KK}h+k_{2}\left|\hat{\Sigma}_{KK}-\tilde{\Sigma}_{KK}\right|_{\infty}\\
 & \leq & h^{T}\tilde{\Sigma}_{KK}h+bk_{1}k_{2}\sqrt{\frac{\log\max(p,\, d)}{n}},
\end{eqnarray*}
where $h\in\mathbb{R}^{k_{2}}$ is a unit-norm minimal eigenvector
of $\tilde{\Sigma}_{KK}$. The last inequality follows from the bounds
on $\left|\frac{(\hat{X}-X^{*})^{T}X^{*}}{n}\right|_{\infty}$ and
$\left|\frac{(\hat{X}-X^{*})^{T}(\hat{X}-X^{*})}{n}\right|_{\infty}$
from the proof for Lemma 6.1 with probability at least $1-c_{1}\exp(-c_{2}\log\max(p,\, d))$.
Therefore, if $\frac{k_{1}^{2}k_{2}^{2}\log\max(p,\, d)}{n}=o(1)$,
then we have 

\[
\lambda_{\min}(\tilde{\Sigma}_{KK})\geq\frac{\lambda_{\min}(\hat{\Sigma}_{KK})}{2}
\]

\begin{eqnarray}
\Longrightarrow\left\Vert \tilde{\Sigma}_{KK}^{-1}\right\Vert _{2} & \leq & \frac{2}{\lambda_{\min}(\hat{\Sigma}_{KK})}\nonumber \\
 & \leq & \frac{4}{\lambda_{\min}(\Sigma_{KK})},
\end{eqnarray}
where the last inequality follows from bound (15) from the proof for
Lemma 6.11. Putting everything together, we have 
\[
\left\Vert \hat{\Sigma}_{KK}^{-1}-\tilde{\Sigma}_{KK}^{-1}\right\Vert _{\infty}\leq\frac{2\sqrt{k_{2}}}{\lambda_{\min}(\Sigma_{KK})}\frac{\lambda_{\min}(\Sigma_{KK})}{\sqrt{k_{2}}}\frac{4}{\lambda_{\min}(\Sigma_{KK})}=\frac{8}{\lambda_{\min}(\Sigma_{KK})}.
\]
with probability at least $1-c_{1}\exp(-c_{2}\log\max(p,\, d))$.\textbf{ }

For Part (ii) of Lemma 6.12, we can bound the terms using results
from Lemma 6.4 (bounds (8) and (9)) instead of Lemma 6.2. Denote the
element $(j^{'},\, j)$ of the matrix difference $\tilde{\Sigma}_{K^{c}K}-\hat{\Sigma}_{K^{c}K}$
by $w_{j^{'}j}$. From the proof for Lemma 6.4, we have 
\[
\left|\frac{1}{n}\mathbf{x}_{j^{'}}^{*T}(\hat{\mathbf{x}}_{j}-\mathbf{x}_{j}^{*})\right|\leq\frac{c\sigma_{\eta}}{\lambda_{\min}(\Sigma_{Z})}\sqrt{\frac{k_{1}\log\max(p,\, d)}{n}}\sup_{v^{j}\in\mathbb{K}(k_{1},d_{j})}\left|\frac{1}{n}\sum_{i=1}^{n}x_{ij^{'}}^{*}\mathbf{z}_{ij}v^{j}\right|.
\]
Following the discretization argument as in the proof for Lemma 6.4,
we have 
\[
\mathbb{P}\left(\sup_{v^{j}\in\mathbb{K}(k_{1},\, d_{j})}\left|\frac{1}{n}x_{j^{'}}^{*T}\mathbf{z}_{j}v^{j}-\mathbb{E}(x_{1j^{'}}^{*}\mathbf{z}_{1j}v^{j})\right|\geq t\right)\leq2\exp(-cn\min(\frac{t^{2}}{\sigma_{X^{*}}^{2}\sigma_{W}^{2}},\,\frac{t}{\sigma_{X^{*}}\sigma_{W}})+k_{1}\log d).
\]
Under the condition $\frac{k_{1}k_{2}^{2}\log\max(p,\, d)}{n}=o(1)$,
setting $t=\frac{\varepsilon\lambda_{\min}(\Sigma_{Z})}{c\sigma_{\eta}}\sqrt{\frac{n}{k_{2}^{2}k_{1}\log\max(p,\, d)}}$
for any $\varepsilon>0$ yields 
\[
\mathbb{P}\left\{ \left\Vert \tilde{\Sigma}_{K^{c}K}-\hat{\Sigma}_{K^{c}K}\right\Vert _{\infty}\geq\varepsilon\right\} 
\]
\begin{eqnarray*}
 & \leq & (p-k_{2})k_{2}\mathbb{P}\left\{ |w_{j^{'}j}|\geq\frac{\varepsilon}{k_{2}}\right\} \\
 & \leq & 2(p-k_{2})k_{2}\cdot\exp(-cn\min(\frac{n\lambda_{\min}^{2}(\Sigma_{Z})\varepsilon^{2}}{\sigma_{\eta}^{2}\sigma_{X^{*}}^{2}\sigma_{W}^{2}k_{1}k_{2}^{2}\log\max(p,\, d)},\,\frac{\sqrt{n}\lambda_{\min}(\Sigma_{Z})\varepsilon}{\sigma_{\eta}\sigma_{X^{*}}\sigma_{W}\sqrt{k_{1}}k_{2}\sqrt{\log\max(p,\, d)}})+k_{1}\log d)\\
 &  & +c_{1}\exp(-c_{2}\log\max(p,\, d)).
\end{eqnarray*}
Bound (21) can be obtained in a similar way except that the pre-factor
$(p-k_{2})$ is replaced by $k_{2}$.

To prove the last bound (22), applying the same argument as in Lemma
6.11, we have 
\begin{eqnarray*}
\left\Vert \hat{\Sigma}_{KK}^{-1}-\Sigma_{KK}^{-1}\right\Vert _{\infty} & \leq & \frac{\sqrt{k_{2}}}{\lambda_{\min}(\hat{\Sigma}_{KK})}\left\Vert \hat{\Sigma}_{KK}-\tilde{\Sigma}_{KK}\right\Vert _{2}\left\Vert \tilde{\Sigma}_{KK}^{-1}\right\Vert _{2}\\
 & \leq & \frac{2\sqrt{k_{2}}}{\lambda_{\min}(\Sigma_{KK})}\left\Vert \hat{\Sigma}_{KK}-\tilde{\Sigma}_{KK}\right\Vert _{2}\left\Vert \tilde{\Sigma}_{KK}^{-1}\right\Vert _{2},
\end{eqnarray*}
where the last inequality comes from bound (15). 

To bound the term $\left\Vert \hat{\Sigma}_{KK}-\tilde{\Sigma}_{KK}\right\Vert _{2}$,
applying bound (21) with $\varepsilon=\frac{\lambda_{\min}(\Sigma_{KK})}{\sqrt{k_{2}}}$
yields 
\begin{eqnarray}
\left\Vert \hat{\Sigma}_{KK}-\tilde{\Sigma}_{KK}\right\Vert _{2} & \leq & \left\Vert \hat{\Sigma}_{KK}-\tilde{\Sigma}_{KK}\right\Vert _{1}\nonumber \\
 & \leq & \frac{\lambda_{\min}(\Sigma_{KK})}{\sqrt{k_{2}}},
\end{eqnarray}
with probability at least
\[
1-2\cdot\exp(-cn\min(\frac{n\lambda_{\min}^{4}(\Sigma_{KK})}{\sigma_{\eta}^{2}\sigma_{X^{*}}^{2}\sigma_{W}^{2}k_{1}k_{2}^{3}\log\max(p,\, d)},\,\frac{\sqrt{n}\lambda_{\min}^{2}(\Sigma_{KK})}{\sigma_{\eta}\sigma_{X^{*}}\sigma_{W}\sqrt{k_{1}}k_{2}^{\nicefrac{3}{2}}\sqrt{\log\max(p,\, d)}})+k_{1}\log d+2\log k_{2})
\]
\[
-c_{1}\exp(-c_{2}\log\max(p,\, d))\geq1-O\left(\frac{1}{\max(p,\, d)}\right)
\]
if $\frac{\max(k_{1}k_{2}^{3/2}\log d,\, k_{1}^{1/2}k_{2}^{3/2}\log p)}{n}=O(1)$. 

To bound the term\textbf{ $\left\Vert \tilde{\Sigma}_{KK}^{-1}\right\Vert _{2}$},
we have, again
\[
\lambda_{\min}(\hat{\Sigma}_{KK})\leq h^{T}\tilde{\Sigma}_{KK}h+h^{T}(\hat{\Sigma}_{KK}-\tilde{\Sigma}_{KK})h
\]
where $h\in\mathbb{R}^{k_{2}}$ is a unit-norm minimal eigenvector
of $\tilde{\Sigma}_{KK}$. By bounds (6) and (7) in Lemma 6.3 and
choosing $s=s(r):=\frac{1}{c}\frac{n}{k_{1}^{r}\log\max(p,\, d)}\min\left\{ \frac{\lambda_{\min}^{2}(\Sigma_{X^{*}})}{\sigma_{X^{*}}^{4}},\,1\right\} $,
$r\in[0,\,1]$, under the condition $\frac{1}{n}k_{1}^{3-2r}\log\max(p,\, d)=o(1)$,
we have 
\begin{eqnarray*}
\left|h^{T}\left(\hat{\Sigma}_{KK}-\tilde{\Sigma}_{KK}\right)h\right| & \leq & bk_{1}^{3/2-r}\sqrt{\frac{\log\max(p,\, d)}{n}}(|h|_{2}^{2}+\frac{1}{s}|h|_{1}^{2})\\
 & \leq & b\max\{k_{1}^{3/2-r}\sqrt{\frac{\log\max(p,\, d)}{n}},\, k_{2}\frac{k_{1}^{r}\log\max(p,\, d)}{n}\}
\end{eqnarray*}
with probability at least $1-c_{1}\exp(-c_{2}\log\max(p,\, d))$.
Furthermore, if 
\[
\frac{1}{n}\min_{r\in[0,\,1]}\max\left\{ k_{1}^{3-2r}\log d,\, k_{1}^{3-2r}\log p,\, k_{1}^{r}k_{2}\log d,\, k_{1}^{r}k_{2}\log p\right\} =o(1),
\]
we have 
\[
\lambda_{\min}(\tilde{\Sigma}_{KK})\geq\frac{\lambda_{\min}(\hat{\Sigma}_{KK})}{2}
\]

\begin{eqnarray}
\Longrightarrow\left\Vert \tilde{\Sigma}_{KK}^{-1}\right\Vert _{2} & \leq & \frac{2}{\lambda_{\min}(\hat{\Sigma}_{KK})}\nonumber \\
 & \leq & \frac{4}{\lambda_{\min}(\Sigma_{KK})},
\end{eqnarray}
where the last inequality follows from bound (15) from the proof for
Lemma 6.11. Because the argument for showing Lemma 6.1 also works
under the assumptions of Lemma 6.3, we can combine the scaling $\frac{k_{1}^{2}k_{2}^{2}\log\max(p,\, d)}{n}=o(1)$
from the proof for bound (24) with the scaling 
\[
\frac{\min_{r\in[0,\,1]}\max\left\{ k_{1}^{3-2r}\log d,\, k_{1}^{3-2r}\log p,\, k_{1}^{r}k_{2}\log d,\, k_{1}^{r}k_{2}\log p\right\} }{n}=o(1)
\]
from above to obtain a more optimal scaling of the required sample
size 
\[
\frac{1}{n}\min\left\{ k_{1}^{2}k_{2}^{2}\log\max(p,\, d),\,\min_{r\in[0,\,1]}\max\left\{ k_{1}^{3-2r}\log d,\, k_{1}^{3-2r}\log p,\, k_{1}^{r}k_{2}\log d,\, k_{1}^{r}k_{2}\log p\right\} \right\} =o(1).
\]
Putting everything together, we have 
\[
\left\Vert \hat{\Sigma}_{KK}^{-1}-\tilde{\Sigma}_{KK}^{-1}\right\Vert _{\infty}\leq\frac{2\sqrt{k_{2}}}{\lambda_{\min}(\Sigma_{KK})}\frac{\lambda_{\min}(\Sigma_{KK})}{\sqrt{k_{2}}}\frac{4}{\lambda_{\min}(\Sigma_{KK})}=\frac{8}{\lambda_{\min}(\Sigma_{KK})}.
\]
with probability at least $1-c_{1}\exp(-c_{2}\log\max(p,\, d))$.\textbf{
$\square$}\\
\textbf{}\\
\textbf{Lemma 6.13}: (i) Suppose the conditions in Lemma 6.2 hold.
With the choice of the tuning parameter \\
\\
$\lambda_{n}\geq c\frac{48(2-\frac{\phi}{4})}{\phi}\max\left\{ \frac{\sigma_{\eta}\max_{j,j^{'}}|\textrm{cov}(x_{1j^{'}}^{*},\,\mathbf{z}_{1j})|_{\infty}|\beta^{*}|_{1}}{\lambda_{\min}(\Sigma_{Z})}k_{1}\sqrt{\frac{\log\max(p,\, d)}{n}},\:\sigma_{X^{*}}\sigma_{\eta}|\beta^{*}|_{1}\sqrt{\frac{\log p}{n}},\:\sigma_{X^{*}}\sigma_{\epsilon}\sqrt{\frac{\log p}{n}}\right\} $
\[
\asymp k_{2}k_{1}\sqrt{\frac{\log\max(p,\, d)}{n}},
\]
for some sufficiently large constant $c>0$, under the condition $\frac{\max\{k_{1}^{2}\log d,\, k_{1}^{2}\log p\}}{n}=o(1)$,
then, we have 
\[
\left(2-\frac{\phi}{4}\right)\left|\hat{X}^{T}\frac{e}{n\lambda_{n}}\right|_{\infty}\leq\frac{\phi}{8},
\]
with probability at least $1-c_{1}\exp(-c_{2}\log\min(p,\, d))$.
(ii) Suppose the conditions in Lemma 6.4 hold. Then the same result
can be obtained with the choice of tuning parameter\\
\\
$\lambda_{n}\geq c\frac{48(2-\frac{\phi}{4})}{\phi}\max\left\{ \frac{\sigma_{\eta}\max_{j,j^{'}}|\textrm{cov}(x_{1j^{'}}^{*},\,\mathbf{z}_{1j})|_{\infty}|\beta^{*}|_{1}}{\lambda_{\min}(\Sigma_{Z})}\sqrt{\frac{k_{1}\log\max(p,\, d)}{n}},\:\sigma_{X^{*}}\sigma_{\eta}|\beta^{*}|_{1}\sqrt{\frac{\log p}{n}},\:\sigma_{X^{*}}\sigma_{\epsilon}\sqrt{\frac{\log p}{n}}\right\} $
\[
\asymp k_{2}\sqrt{\frac{k_{1}\log\max(p,\, d)}{n}},
\]
and the condition $\frac{\max\{k_{1}\log d,\, k_{1}\log p\}}{n}=o(1)$.\textbf{}\\
\textbf{}\\
\textbf{Proof}. Recall from the proof for Lemma 6.2, 
\begin{eqnarray*}
\frac{1}{n}\hat{X}^{T}e & = & \frac{1}{n}\hat{X}^{T}\left[(X^{*}-\hat{X})\beta^{*}+\boldsymbol{\eta}\beta^{*}+\epsilon\right]\\
 & = & \frac{1}{n}X^{*T}\left[(X^{*}-\hat{X})\beta^{*}+\boldsymbol{\eta}\beta^{*}+\epsilon\right]+\frac{1}{n}(\hat{X}-X^{*})^{T}\left[(X^{*}-\hat{X})\beta^{*}+\boldsymbol{\eta}\beta^{*}+\epsilon\right].
\end{eqnarray*}
Hence, 
\begin{eqnarray}
|\frac{1}{n\lambda_{n}}\hat{X}^{T}e|_{\infty} & \leq & |\frac{1}{n\lambda_{n}}X^{*T}(\hat{X}-X^{*})\beta^{*}|_{\infty}+|\frac{1}{n\lambda_{n}}X^{*T}\boldsymbol{\eta}\beta^{*}|_{\infty}+|\frac{1}{n\lambda_{n}}X^{*T}\epsilon|_{\infty}\\
 & + & |\frac{1}{n\lambda_{n}}(\hat{X}-X^{*})^{T}(\hat{X}-X^{*})\beta^{*}|_{\infty}+|\frac{1}{n\lambda_{n}}(\hat{X}-X^{*})^{T}\boldsymbol{\eta}\beta^{*}|_{\infty}+|\frac{1}{n\lambda_{n}}(\hat{X}-X^{*})^{T}\epsilon|_{\infty}.\nonumber 
\end{eqnarray}
From the proof for Lemma 6.2, we have 
\begin{eqnarray*}
|\frac{1}{n}X^{*T}(\hat{X}-X^{*})\beta^{*}|_{\infty} & \leq & c\frac{\sigma_{\eta}\max_{j^{'},\, j}|\mathbb{E}(x_{ij^{'}}^{*},\,\mathbf{z}_{ij})|_{\infty}}{\lambda_{\min}(\Sigma_{Z})}|\beta^{*}|_{1}k_{1}\sqrt{\frac{\log\max(p,\, d)}{n}},\\
|\frac{1}{n}(\hat{X}-X^{*})^{T}(\hat{X}-X^{*})\beta^{*}|_{\infty} & \leq & c\frac{\sigma_{\eta}^{2}\max_{j^{'},\, j}|\mathbb{E}(\mathbf{z}_{ij^{'}},\,\mathbf{z}_{ij})|_{\infty}}{\lambda_{\min}^{2}(\Sigma_{Z})}|\beta^{*}|_{1}k_{1}^{2}\frac{\log\max(p,\, d)}{n},\\
|\frac{1}{n}X^{*T}\boldsymbol{\eta}\beta^{*}|_{\infty} & \leq & c\sigma_{X^{*}}\sigma_{\eta}|\beta^{*}|_{1}\sqrt{\frac{\log p}{n}},\\
|\frac{1}{n}(X^{*}-\hat{X})^{T}\boldsymbol{\eta}\beta^{*}|_{\infty} & \leq & c\frac{\sigma_{Z}\sigma_{\eta}^{2}}{\lambda_{\min}(\Sigma_{Z})}|\beta^{*}|_{1}k_{1}\frac{\log\max(p,\, d)}{n},\\
|\frac{1}{n}X^{*T}\epsilon|_{\infty} & \leq & c\sigma_{X^{*}}\sigma_{\epsilon}\sqrt{\frac{\log p}{n}},\\
|\frac{1}{n}(X^{*}-\hat{X})^{T}\epsilon|_{\infty} & \leq & c\frac{\sigma_{Z}\sigma_{\epsilon}\sigma_{\eta}}{\lambda_{\min}(\Sigma_{Z})}k_{1}\frac{\log\max(p,\, d)}{n},
\end{eqnarray*}
with probability at least $1-c_{1}\exp(-c_{2}\log\min(p,\, d))$.
Therefore as long as \\
\\
$\lambda_{n}\geq c^{'}\frac{48(2-\frac{\phi}{4})}{\phi}\max\left\{ \frac{\sigma_{\eta}\max_{j,j^{'}}|\textrm{cov}(x_{1j^{'}}^{*},\,\mathbf{z}_{1j})|_{\infty}|\beta^{*}|_{1}}{\lambda_{\min}(\Sigma_{Z})}k_{1}\sqrt{\frac{\log\max(p,\, d)}{n}},\:\sigma_{X^{*}}\sigma_{\eta}|\beta^{*}|_{1}\sqrt{\frac{\log p}{n}},\:\sigma_{X^{*}}\sigma_{\epsilon}\sqrt{\frac{\log p}{n}}\right\} $
\[
\asymp k_{2}k_{1}\sqrt{\frac{\log\max(p,\, d)}{n}},
\]
for some sufficiently large constant $c^{'}>0$, under the condition
$\frac{\max\{k_{1}^{2}\log d,\, k_{1}^{2}\log p\}}{n}=o(1)$, we have

\begin{eqnarray*}
|\frac{1}{n\lambda_{n}}X^{*T}(\hat{X}-X^{*})\beta^{*}|_{\infty} & \leq & \frac{\phi}{48(2-\frac{\phi}{4})},\\
|\frac{1}{n\lambda_{n}}X^{*T}\boldsymbol{\eta}\beta^{*}|_{\infty} & \leq & \frac{\phi}{48(2-\frac{\phi}{4})},\\
|\frac{1}{n}X^{*T}\epsilon|_{\infty} & \leq & \frac{\phi}{48(2-\frac{\phi}{4})},
\end{eqnarray*}
\begin{eqnarray*}
|\frac{1}{n\lambda_{n}}(\hat{X}-X^{*})^{T}(\hat{X}-X^{*})\beta^{*}|_{\infty}\precsim & k_{1}\sqrt{\frac{\log\max(p,\, d)}{n}}=o(1) & \leq\frac{\phi}{48(2-\frac{\phi}{4})},\\
|\frac{1}{n\lambda_{n}}(X^{*}-\hat{X})^{T}\boldsymbol{\eta}\beta^{*}|_{\infty}\precsim & \sqrt{\frac{\log\max(p,\, d)}{n}}=o(1) & \leq\frac{\phi}{48(2-\frac{\phi}{4})},\\
|\frac{1}{n}(X^{*}-\hat{X})^{T}\epsilon|_{\infty}\precsim & \sqrt{\frac{\log\max(p,\, d)}{n}}=o(1) & \leq\frac{\phi}{48(2-\frac{\phi}{4})}.
\end{eqnarray*}
with probability at least $1-c_{1}\exp(-c_{2}\log\min(p,\, d))$.
Putting everything together, we have 
\[
\left(2-\frac{\phi}{4}\right)\left|\hat{X}^{T}\frac{e}{n\lambda_{n}}\right|_{\infty}\leq\frac{\phi}{8},
\]
with probability at least $1-c_{1}\exp(-c_{2}\log\min(p,\, d))$.

The proof for Part (ii) of Lemma 6.13 follows from the similar argument
for proving Part (i) except that we bound the terms $|\frac{1}{n\lambda_{n}}X^{*T}(\hat{X}-X^{*})\beta^{*}|_{\infty}$,
$|\frac{1}{n\lambda_{n}}(\hat{X}-X^{*})^{T}(\hat{X}-X^{*})\beta^{*}|_{\infty}$,
$|\frac{1}{n}(X^{*}-\hat{X})^{T}\boldsymbol{\eta}\beta^{*}|_{\infty}$,
and $|\frac{1}{n}(X^{*}-\hat{X})^{T}\epsilon|_{\infty}$ using the
discretization argument as in the proof for Lemma 6.4. $\square$

\end{document}